\pdfoutput=1
\RequirePackage{ifpdf}
\ifpdf 
\documentclass[pdftex]{sigma}
\else
\documentclass{sigma}
\fi

\usepackage{youngtab}
\usepackage[all]{xy}

\def\bA{{\mathbb A}}

\def\CC{{\mathbb C}}

\def\GG{{\mathbb G}}

\def\ZZ{{\mathbb Z}}

\def\NN{{\mathbb N}}
\def\PP{{\mathbb P}}

\def\bH{{\mathbf F}}

\def\bH{{\mathbf H}}

\def\bby{{\mathbf y}}

\def\cA{{\mathcal A}}

\def\cI{{\mathcal I}}

\def\cM{{\mathcal M}}

\def\cO{{\mathcal O}}
\def\cP{{\mathcal P}}
\def\cQ{{\mathcal Q}}

\def\fB{{\mathfrak B}}
\def\fC{{\mathfrak C}}

\def\fE{{\mathfrak E}}

\def\fK{{\mathfrak K}}

\def\fY{{\mathfrak Y}}
\def\fZ{{\mathfrak Z}}
\def\fa{{\mathfrak a}}
\def\fb{{\mathfrak b}}
\def\fc{{\mathfrak c}}
\def\fd{{\mathfrak d}}
\def\fe{{\mathfrak e}}

\def\fh{{\mathfrak h}}

\def\fk{{\mathfrak k}}

\def\fm{{\mathfrak m}}

\def\fr{{\mathfrak r}}
\def\fs{{\mathfrak s}}

\def\fy{{\mathfrak y}}

\def\tA{{\widetilde A}}

\def\tR{{\widetilde R}}

\def\tQ{{\widetilde Q}}
\def\tR{{\widetilde R}}

\def\tb{{\widetilde b}}

\def\td{{\widetilde d}}

\def\th{{\widetilde h}}

\def\tj{{\widetilde j}}
\def\tk{{\widetilde k}}

\def\ty{{\widetilde y}}
\def\tw{{\widetilde w}}

\def\tdelta{{\widetilde \delta}}

\def\tvarphi{{\widetilde \varphi}}

\def\tfE{{\widetilde{\mathfrak E}}}

\def\tfe{{\widetilde{\mathfrak e}}}

\def\tfy{{\widetilde{\mathfrak y}}}

\def\BA{{\overline A}}

\def\BH{{\overline H}}
\def\BI{{\overline I}}

\def\BfE{{\overline{\mathfrak E}}}

\def\Bfe{{\overline{\mathfrak e}}}

\def\Bfr{{\overline{\mathfrak r}}}

\def\hc{{\widehat c}}

\def\hk{{\widehat k}}

\def\hH{{\widehat H}}

\def\hN{{\widehat N}}

\def\hR{{\widehat R}}
\def\hS{{\widehat S}}

\def\hfe{{\widehat{\mathfrak e}}}

\def\hfy{{\widehat{\mathfrak y}}}

\def\hby{\widehat {{\mathbf y}}}

\def\hlambda{{\widehat \lambda}}

\def\hphi{{\widehat \phi}}

\def\hzeta{{\widehat \zeta}}

\def\hUpsilon{{\widehat \Upsilon}}

\def\rfy{{\mathring{\mathfrak y}}}

\def\lH{{\overline H}}

\def\Hom{\mathrm{Hom}}

\def\diff{\mathrm{diff}}

\def\Ker{\operatorname{Ker}}

\def\Div{\operatorname{div}}

\def\Ann{\operatorname{Ann}}
\def\Aut{\operatorname{Aut}}
\def\Gal{\operatorname{Gal}}

\def\Spec{\operatorname{Spec}}

\def\wt{{\mathrm {wt}}}

\def\nuI#1{{\nu^{\mathrm{I}}_{#1}}}

\numberwithin{equation}{section}

\newtheorem{Theorem}{Theorem}[section]
\newtheorem*{Theorem*}{Theorem}
\newtheorem{Corollary}[Theorem]{Corollary}
\newtheorem{Lemma}[Theorem]{Lemma}
\newtheorem{Proposition}[Theorem]{Proposition}
 { \theoremstyle{definition}
\newtheorem{Definition}[Theorem]{Definition}

\newtheorem{Example}[Theorem]{Example}
\newtheorem{Remark}[Theorem]{Remark} }

\begin{document}
\allowdisplaybreaks

\newcommand{\arXivNumber}{2207.01905}

\renewcommand{\PaperNumber}{098}

\FirstPageHeading

\ShortArticleName{Complementary modules of Weierstrass canonical forms}

\ArticleName{Complementary Modules\\ of Weierstrass Canonical Forms}

\Author{Jiryo KOMEDA~$^{\rm a}$, Shigeki MATSUTANI~$^{\rm b}$ and Emma PREVIATO~$^{\rm c}$}

\AuthorNameForHeading{J.~Komeda, S.~Matsutani and E.~Previato}

\Address{$^{\rm a)}$~Department of Mathematics,
Center for Basic Education and Integrated Learning,\\
\hphantom{$^{\rm a)}$}~Kanagawa Institute of Technology, 1030 Shimo-Ogino, Atsugi, Kanagawa 243-0292, Japan}

\Address{$^{\rm b)}$~Faculty of Electrical, Information and Communication Engineering,\\
\hphantom{$^{\rm b)}$}~Kanazawa University, Kakuma Kanazawa, 920-1192, Japan}
\EmailD{\href{s-matsutani@se.kanazawa-u.ac.jp}{s-matsutani@se.kanazawa-u.ac.jp}}

\Address{$^{\rm c)}$~Department of Mathematics and Statistics, Boston University, Boston, MA 02215-2411, USA}

\ArticleDates{Received July 12, 2022, in final form December 07, 2022; Published online December 18, 2022}

\Abstract{The Weierstrass curve is a pointed curve $(X,\infty)$ with a numerical semigroup $H_X$, which is a normalization of the curve given by the Weierstrass canonical form, $y^r + A_{1}(x) y^{r-1} + A_{2}(x) y^{r-2} +\dots + A_{r-1}(x) y + A_{r}(x)=0$ where each $A_j$ is a polynomial in $x$ of degree $\leq j s/r$ for certain coprime positive integers $r$ and $s$, $r<s$, such that the generators of the Weierstrass non-gap sequence $H_X$ at $\infty$ include $r$ and $s$. The Weierstrass curve has the projection $\varpi_r\colon X \to {\mathbb P}$, $(x,y)\mapsto x$, as a covering space. Let $R_X := {\mathbf H}^0(X, {\mathcal O}_X(*\infty))$ and $R_{\mathbb P} := {\mathbf H}^0({\mathbb P}, {\mathcal O}_{\mathbb P}(*\infty))$ whose affine part is ${\mathbb C}[x]$. In this paper, for every Weierstrass curve~$X$, we show the explicit expression of the complementary module~$R_X^{\mathfrak c}$ of $R_{\mathbb P}$-module $R_X$ as an extension of the expression of the plane Weierstrass curves by Kunz. The extension naturally leads the explicit expressions of the holomorphic one form except~$\infty$, ${\mathbf H}^0({\mathbb P}, {\mathcal A}_{\mathbb P}(*\infty))$ in terms of $R_X$. Since for every compact Riemann surface, we find a Weierstrass curve that is bi-rational to the surface, we also comment that the explicit expression of $R_X^{\mathfrak c}$ naturally leads the algebraic construction of generalized Weierstrass' sigma functions for every compact Riemann surface and is also connected with the data on how the Riemann surface is embedded into the universal Grassmannian manifolds.}


\Classification{14H55; 14H50; 16S36; 13H10}

\section{Introduction}\label{introduction}

The Weierstrass $\sigma$ function is defined for an elliptic curve of Weierstrass' equation $y^2 = 4 x^3 - g_2 x - g_3$ in Weierstrass' elliptic function theory \cite{WeiWV,WhittakerWatson}.
Since $\bigr(\wp(u)=-\frac{{\rm d}^2}{{\rm d} u^2}\log\sigma(u), \frac{{\rm d}\wp(u)}{{\rm d} u}\bigr)$ is identical to a point $(x,y)$ of the curve, we can use the equivalence between the algebraic objects of the curve and the transcendental objects on its Jacobi variety.
The equivalence leads to the crucial relations among them and their algebraic and analytic properties associated with the elliptic curves~\cite{WeiWV}. These relations and properties affect the various fields of science and mathematics.\looseness=-1

We have studied the generalization of this picture to algebraic curves with higher genera in the series of the studies \cite{KM2020, KMP13, KMP16, KMP19} following Mumford's excellent studies for the hyperelliptic curves \cite{Mum81,Mum84}.

As the elliptic theta function was generalized by Riemann for an Abelian variety, its equivalent function Al was defined for any hyperelliptic curve by Weierstrass and was refined by Klein as a~generalization of the elliptic sigma function.
Since Klein defined his hyperelliptic sigma function by using only the data of the hyperelliptic Riemann surface and Jacobian transcendentally, Baker re-constructed Klein's sigma function by using only the data of the hyperelliptic curve itself from an algebraic viewpoint~\cite{Baker97}.
Buchstaber, Enolskii, and Leykin extend the sigma functions to certain plane curves, so-called $(n,s)$ curves, based on Baker's construction which we call EEL construction due to work by Eilbeck, Enolskii, and Leykin (see \cite{BEL20, EEL00} and references therein).
For the $(n,s)$ curves with the cyclic symmetry, the direct relations between the affine rings and the sigma functions were obtained as the Jacobi inversion formulae \cite{MP08, MP14}. Further, we generalized the sigma functions and the formulae to a particular class of the space curves using the EEL-construction \cite{KMP13, KMP19, MK13}.

Recently as a generalization of Klein's sigma functions, D.~Korotkin and V.~Shramchenko defined the sigma function of every compact Riemann surface transcendentally~\cite{KorotkinS}.
Every compact Riemann surface with a point $P$ has the Weierstrass non-gap sequence at $P$, which is described by a numerical semigroup $H$; we call the numerical semigroup \emph{Weierstrass semigroup}.
In \cite{Nak16}, Nakayashiki refined the sigma function for every compact Riemann surface with Weierstrass semigroup $H$ based on Sato's theory on the universal Grassmannian manifolds (UGM) \cite{SatoN84,SegalWilson}.

On the other hand,
it is well-known that for every compact Riemann surface $Y$ with the Weierstrass semigroup $H_Y$ at a point $P\in Y$,
there is an algebraic curve $X$ which is bi-rational to the surface $Y$ and is obtained by the normalization of the curve satisfying the Weierstrass canonical form, $y^r + A_{1}(x) y^{r-1} + A_{2}(x) y^{r-2} +\cdots + A_{r-1}(x) y + A_{r}(x)$ where each $A_j$ is a polynomial in $x$ of degree $\leq j s/r$ for certain coprime positive integers $r$ and $s$, $r<s$; the point $P\in Y$ corresponds to $\infty \in X$ and the Weierstrass semigroup $H_X$ at $\infty \in X$ is equal to $H_Y$ whose generators include $r$ and $s$
\cite{Baker97, CoppensKato, Kato80, Wei67}.
In this paper, we call such a curve \emph{Weierstrass curve} or \emph{W-curve}. The set of W-curves represents the set of compact Riemann surfaces.
The Weierstrass canonical form
 provides the projection $\varpi_r\colon X \to \PP$ (e.g., $(x,y) \mapsto x)$ as a covering space.
For $R_X = \bH^0(X, \cO_X(*\infty))$ and $R_\PP = \bH^0(\PP, \cO_\PP(*\infty))$, let their quotient fields be~$\cQ(R_X)$ and~$\cQ(R_\PP)$.
Then $\cQ(R_X)$ is a field extension of $\cQ(\PP)$, or $X$ is the holomorphic $r$-sheeted covering on~$\PP$, and~$R_X$ is regarded as an $R_\PP$-module.

The $(n,s)$ curves are identical to the plane W-curves, whereas the space curves studied in \cite{ KMP13, KMP19, MK13} are particular classes of the space W-curves.
In \cite{KMP19}, we implicitly showed that these studies on the sigma functions in \cite{BEL20,KMP13, MP08, MP14, Onishi18} are based on the algebraic structures of~$R_X$ as an $R_\PP$-module for such particular W-curves.
In \cite{KM2020}, we show that if we have the data of the holomorphic one-forms over $X$ except $\infty$, $\bH^0(X, \cA_X(*\infty))$,
we have the correspondence between $R_X$ and the Riemann theta function for the subvarieties of its Jacobi variety associated with~$S^k X$, $0\le k <g$, where $\cA_X$ is the sheaf of the holomorphic one-forms on $X$.

Our purpose of the studies series is to connect the algebra $R_X$ and the sigma function for every W-curve $X$ and show the equivalence between the algebraic objects of the W-curve $X$ and the transcendental objects on its Jacobi variety using the sigma functions as Baker did for the hyperelliptic sigma functions following Weierstrass' elliptic function theory~\cite{Baker97}.
The equivalence also leads to the crucial relations among them and their algebraic and analytic properties associated with the W-curve $X$, like Weierstrass' elliptic function theory.
In other words, we purpose to extend Weierstrass' elliptic function theory \cite{WeiWV} and Mumford's study \cite{Mum81,Mum84} to every W-curve.

It is a critical perspective that Weierstrass himself and related researchers had already obtained some, but not perfect, results \cite{Baker97,WeiWIV}.
Due to the difficulty, these attempts were given up and forgotten for a century.
However, the development of mathematics enables us to revive their approaches.
We also emphasize that the progress in the studies on the sigma function has provided several non-trivial results as a generalization of those in elliptic curves, which had been regarded as impossible for the century, e.g., \cite{BEL20,EEMOP07,EEMOP08, EEL00,FKMPA,MP08,MP14, Nak10b, Onishi18} and reference therein as Mumford did for the hyperelliptic curves using $\theta$ functions in \cite{Mum81,Mum84}.

In the general theory of the algebraic curves \cite{DedekindWeber,Kunz2004,Stichtenoth}, the algebraic curves with the Galois covering are studied well, and it turns out that the complementary module plays important roles and is connected with the K\"ahler differentials.
The origin of the complementary modules is in the study by Dedekind and Weber \cite{DedekindWeber}; Weierstrass already stated some results in \cite{WeiWIV}.
Kunz showed the explicit expression of the complementary module for every plane W-curve \cite[Theorem 15.1]{Kunz2004} though he did not call it W-curve.

To connect Nakayashiki's sigma functions in \cite{Nak16} for pointed compact Riemann surfaces with the Weierstrass semigroup $H_X$ and the algebraic properties of the W-curves $X$, in this paper, we study the algebraic properties of $R_X$ and the complementary module $R_X^\fc$ as $R_\PP$-modules.
This paper aims to show the explicit expression of $R_X^\fc$ since the expression of $R_X^\fc$ enables us to apply the algebraic construction of the sigma function following the EEL construction to every W-curve:
1) $R_X^\fc$ is directly related to the meromorphic one forms $\bH^0(X, \cA_X(*\infty))$, i.e., $\bH^0(X, \cA_X(*\infty))=R_X^\fc {\rm d} x$; we can find the Jaocobi inversion formula as mentioned in \cite{KM2020}.
2) As we show in a follow-up paper \cite{KMP2022b}, by using the structure of $R_X^\fc$, we can define the fundamental differential of the second kind $\Omega$ following the EEL-construction to obtain the Jacobi inversion formula of Nakayashiki's sigma function.

We remark that in Weierstrass' elliptic function theory, if we regard $x(u)=-\frac{{\rm d}^2}{{\rm d} u^2}\log \sigma$ with ${\rm d}u = {\rm d}x/y$ as a differential equation, it is known that $x$ can be integrated twice in $u$ to obtain the elliptic sigma function without theta functions, at least, formally.
Similarly, we regard that the Jacobi inversion formula, i.e., the relation between $R_X$ and the sigma function, for every W-curve also means an algebraic construction of the sigma function as Baker considered~\cite{Baker97}.
In the construction, the most critical and most complicated point is to obtain the explicit expression of the complementary module $R_X^\fc$.
Thus this paper is devoted to the expression.

By using the general theory of the algebraic curves, in this paper, we obtain the explicit expression of the complementary module $R_X^\fc$ for every W-curve~$X$, including a space curve as an extension of Kunz's results on the complementary modules of the plane W-curves~\cite{Kunz2004}.
Though the extension, especially to non-symmetric $H_X$ cases, is complicated, we have already obtained the expressions of particular space curves heuristically in \cite{KMP13, KMP19, MK13}.
By investigating them, we found the essentials of the general constructions of~$R_X^\fc$ in terms of data of~$R_X$ for general W-curves as follows.
Since the W-curves $X$ are characterized by the Weierstrass semigroups, $H_X$~\cite{Wei67}, the monomial curves~$X_H$ and the monomial rings $R_H^Z$ associated with $H_X$ determine the behavior of $R_X$ at $\infty\in X$.
Following such studies on the monomial curves by Kunz~\cite{Kunz1970}, Herzog~\cite{Herzog70}, and Pinkham~\cite{Pinkham74}, we use the multiple groups $\GG_\fm$ action on the monomial curve~$X_H$ and the monomial algebra~$R_H^Z$.
As we assume that the generators of the Weierstrass semigroup~$H_X$ include the minimal number~$r$, the cyclic group $\fC_r$ of order $r$ gives the standard basis of~$H_X$ with respect to~$r$.
The standard basis of $H_X$ governs the monomial algebra and induces the standard basis of~$R_X$ as the $R_\PP$-module globally.
By using the standard basis of~$R_X$, we investigate the $R_\PP$-module structure of~$R_X$.
As a key proposition, we find the characterization of~$R_X$ in Proposition~\ref{2pr:WCF2} which induces an embedding of $m_X$-dimensional space curve~$X$ into~$\PP^{2(m_X-1)}$ with projection to each singular plane curve in~$\PP^2$ in Proposition~\ref{2pr:X_to_PP^m-1}. This embedding leads the trace structure $R_H^Z$ in Lemma~\ref{2lm:h_H} and the structure determines the dual basis with the trace of $R_X$ in Proposition~\ref{2pr:h_RXe_ys}.
Finally we find the explicit expression of~$R_X^\fc$ in Theorem~\ref{2pr:R_Xfc}.
Using it, we describe $\bH^0(X, \cA_X(*\infty))=R_X^\fc {\rm d} x$ in terms of $R_X$ in Theorem~\ref{2lm:nuI_hX}.

Furthermore it is known that $\BH_X^\fc := \ZZ \setminus H_X$ provides the embedding of the algebraic systems associated with $X$ into the UGM due to Segal and Wilson \cite[p.~46]{SegalWilson}.
Since due to Riemann--Roch theorem, $R_X^\fc {\rm d} x$ has the weight sequence which is equal to $\BH_X^\fc-1$, $R_X^\fc {\rm d} x$ might lead the embedding as Nakayashiki did for $(n,s)$ curve in \cite{Nak10b} rather than the spin connection in~\cite{Nak16}.
Our results in this paper can be naturally applied to the generalization of the EEL-constructions~\cite{KMP19} as we mentioned above.
Our previous report on the Riemann constant on the theta function~\cite{KMP16} enables us to construct the sigma function of every W-curve algebraically and connect $R_X$ with the sigma function as we show in the follow-up paper~\cite{KMP2022b}. We mention it shortly in Re\-mark~\ref{rmk:Final}.\looseness=1

In \cite{FKMPA}, the explicit description of the Abelian functions in terms of the sigma function demonstrates the degenerating behavior of the sigma function for the degenerating family of curves given by the Weierstrass canonical form $f_X(x,y)=y^3-x(x-s)(x-b_1)(x-b_2)$ for $s \to 0$ for disjoint non-zero complex numbers $b_1$ and $b_2$ recently, which is much more precise than the known results~\cite{Igusa56}.
The results in this paper with this follow-up paper~\cite{KMP2022b} mean that 1)~as we handle the elliptic functions of an elliptic curve, we can handle the algebraic functions of any W-curve $X$ using the explicit connection between the sigma function for $X$ and the affine ring $R_X$, e.g., their additive structure, Jacobi inversion formula, and differential relations, 2)~as we did in~\cite{FKMPA}, we can basically express the degenerating behavior of sigma function (theta function) for any degenerating family of W-curves, and 3) in terms of them, we could have explicit expressions of the algebraic solutions of KP hierarchy more precisely:
Though it has not been a concern in the study of the integrable system, even for soliton solutions of KP hierarchy, there is no study on explicit description associated with the space curves except the recent interesting work by Kodama and Xie~\cite{KodamaX2021}.
Our results in this paper provide the bases.

Contents are as follows:
Section~\ref{sec:preliminary} consists of the two subsections:
Section~\ref{2sc:App_trace} reviews Dedekind's trace, and its related topics based on Kunz's book \cite{Kunz2004}.
Section~\ref{2sc:WSG} gives the summary of the numerical and Weierstrass semigroups.
We show the Weierstrass canonical forms and Weierstrass curves (W-curves) and their properties in Section~\ref{sec:WCF}.
Using their properties, we find the identities in $R_X$ in Proposition~\ref{2pr:WCF2}, as the first key proposition, to show the $R_\PP$-module structure of $R_X$ and an embedding of $X$ into $\PP^{2(m_X-1)}$, where $m_X$ is the minimal number of the generators of $H_X$.
Section~\ref{sec:CompM} is devoted to the construction of the complementary module~$R_X^\fc$.
After rewriting the tools in Section~\ref{2sc:App_trace} for W-curves~$X$ shortly, we start to review the explicit expression of the complementary module for every plane W-curve of $(m_X=2)$ \cite[Theorem~15.1]{Kunz2004} in Proposition~\ref{2pr:Kunz_m2}.
In order to extend it to a general W-curve $X$ including a~non-symmetric case, we consider the trace structure of the monomial ring $R_H^Z$ in Lemma~\ref{2lm:h_H} as the second key proposition.
Using it, we investigate the global structure of~$R_X$ as the $R_\PP$-module in Proposition~\ref{2pr:h_RXe_ys} and Lemma~\ref{2lm:ideal_hUpsilon} as the third key propositions.
Finally, we construct the complementary module in Theorem~\ref{2pr:R_Xfc} as the first main theorem in this paper.
In Section~\ref{2ssc:W-norm_nuI}, we consider the W-normalized Abelian differentials $\bH^0(X, \cA_X(*\infty))$ and show the second main theorem in Theorem \ref{2lm:nuI_hX}.
In Section~\ref{sec:Exmples}, we provide some examples of our results.

\section{Preliminary}\label{sec:preliminary}

\subsection{Trace and complementary modules}\label{2sc:App_trace}

We review Dedekind's trace and its related topics based on Kunz's book~\cite{Kunz2004} whose origin appeared in the paper by Dedekind and Weber~\cite{DedekindWeber}.

Let $R_P$ be an algebra over $\CC$ and $R_Y$ an $R_P$ algebra such that
$R_Y=\bigoplus_{i=0}^{r-1}R_P\bby_i$, $\bby_i\in R_Y$.

The dual of $R_Y$ is defined by
\[
\omega_{R_Y/R_P}=\operatorname{Hom}_{R_P}(R_Y, R_P)
\]
with the basis $\{\bby_i^*\}_{i\in \ZZ_r} \subset \omega_{R_Y/R_P}$ satisfying $\bby_i^*\bby_j=\delta_{i j}$.
We assume that for $x, y\in R_Y$ and $z\in \omega_{R_Y/R_P}$, $(x\circ z)(y)=z(xy)$, and we also regard $\omega_{R_Y/R_P}$ as an $R_Y$-module.

Let us introduce \emph{standard trace},
$\tau_{R_Y/R_P}:=\sum_{i=0}^{r-1}\bby_i\circ\bby_i^*$,
 and its complementary module $R_Y^\fc$ with respect to $\tau_{R_Y/R_P}$ given by
\[
R_Y^\fc:=\{z \in \cQ(R_Y) \, | \, \tau_{R_Y/R_P}(z R_Y) \subset R_P \}.
\]
We construct the complementary module $R_Y^\fc$ as follows because $R_Y^\fc$ is directly connected with the K\"ahler differentials.

\begin{Lemma}\label{2lm:trace_1.1}
There are structure coefficients $a_{i j k}$ in $R_P$ satisfying
$\bby_i \bby_j = \sum_{k} a_{i j k} \bby_k$ and ${a_{ijk}=a_{jik}}$, which determines the structure of the $R_P$-module $R_Y$.
Then the standard trace shows $\tau_{R_Y/R_P}(\bby_j)=\sum_i a_{i j i}$.
\end{Lemma}

\begin{proof}
The former statement is obvious and due to the definition of the standard trace, we have $\tau_{R_Y/R_P}(\bby_j)=\sum_{i}
\bby_i\circ \bby_i^*(\bby_j)=\sum_{i} \bby_i^*(\bby_i\bby_j) =\sum_{i} \bby_i^*(\sum_k a_{i j k} \bby_k)=\sum_i a_{i j i}$.
\end{proof}

The $R_Y$-action on the $R_P$-module $R_Y$, $x\colon R_Y \to R_Y$ for $x\in R_Y$, has the matrix expression $M_x$, i.e., $(M_x)_{ij}=\sum_k x_k a_{jki}$ since
$(xz)_i=\bby_i^*(\sum_{kj} x_k z_j\bby_k\bby_j)
=\bby_i^*(\sum_{kj\ell} x_k z_ja_{jk\ell}\bby_\ell)
=\sum_j(\sum_k x_k a_{jki}) z_j$.
Lemma~\ref{2lm:trace_1.1} asserts that the standard trace $\tau_{R_Y/R_P}(x)$ agrees with trace of~$M_x$, i.e, $\tau_{R_Y/R_P}(x)=\sum_{i} (M_x)_{ii}$.

\begin{Definition}\label{2df:trace_1.1}
If there is an element $\tau$ in $\omega_{R_Y/R_P}$ such that $\omega_{R_Y/R_P}=$ $R_Y\circ\tau$ as an $R_P$-module, i.e., $\omega_{R_Y/R_P}$ is a free $R_P$-module with the basis $\{\tau\}$, we say that $R_Y$ has the trace $\tau$.
\end{Definition}

We note that if $R_Y/R_P$ is separable, the standard trace $\tau_{R_Y/R_P}$ is a trace in this definition.

\begin{Lemma}\label{2lm:trace_1.2}
If $R_Y$ has a trace $\tau$, the following hold.
\begin{enumerate}\itemsep=0pt
\item[$1.$] If an element $a$ in $R_Y$ satisfies $a \circ\tau =0$, then $a=0$.
It means that $\omega_{R_Y/R_P}\cong R_Y\circ \tau$ as an $R_Y$-module.

\item[$2.$] For a basis $\{\bby_i\}$ of $R_Y$ as an $R_P$-module, there exists a subset $\{\hby_i\}\subset R_Y$ satisfying $\tau(\hby_i \bby_j)=\delta_{i j}$.
$($We call $\{\hby_i\}$ the dual basis with respect to the trace $\tau.)$
Then $\tau_{R_Y/R_P}=\big(\sum_{i=0}^{r-1} \hby_i \bby_i\big)\circ \tau$.
\end{enumerate}
\end{Lemma}

\begin{proof}(1) For any $x\in R_Y$, $0=a\circ \tau(x)=\tau(ax)=x\circ \tau(a)$. For every $b \in \omega_{R_Y/R_P}$, $b(a)=0$ implies that $a=0$.
The $R_P$-morphism $\omega_{R_Y/R_P}\to R_Y\circ \tau$ is injective.
(2) For the dual basis~$\{\bby_i^*\}$ of $\omega_{R_Y/R_P}$ such that $\bby_i^*(\bby_j)=\delta_{i j}$, we can find an element $\hby_i\in R_Y$ satisfying $\bby^*_i = \hby_i\circ \tau$.
Then $\tau(\hby_i \bby_j)=\bby_i^*(\bby_j)=\delta_{i j}$, and $\tau_{R_Y/R_P}=\sum_{i=0}^{r-1}\bby_i\circ \bby_i^*=\big(\sum_{i=0}^{r-1} \hby_i \bby_i\big)\circ \tau$.
\end{proof}

We construct $\tau$ and $\tau_{R_Y/R_P}$ in terms of the enveloping algebra
$R_Y^e:=R_Y \otimes_{R_P}R_Y$.
For~$R_Y^e$, the standard multiplication $\mu\colon R_Y^e\to R_Y$ is defined by $\mu(a\otimes b)=ab$.

\begin{Lemma}The kernel of $\mu$, $\Ker \mu$, is generated by $\{\bby_i \otimes 1 - 1 \otimes \bby_i\}_{i=0, \dots, r-1}$.
\end{Lemma}

\begin{proof}
Following \cite[Theorem~G.7]{Kunz2004}, we show it.
Let $I:=\langle \{a \otimes 1 - 1 \otimes a\}_{a\in R_Y}\rangle_{R_P}$.
Clearly $I \subset \Ker \mu$.
There are surjective $R_P$-homomorphisms,
\[
R_Y \otimes_{R_P} R_Y \xrightarrow{p_\mu} (R_Y \otimes_{R_P} R_Y)/I
\xrightarrow{p_\mu'}
(R_Y \otimes_{R_P} R_Y)/\Ker \mu \cong R_Y.
\]
For $a, b \in R_Y$, we have
$a\otimes b=(a\otimes 1)(1\otimes b)$
$=-(a\otimes 1)(b\otimes 1 -1\otimes b)+(ab\otimes 1)$.
Thus there is an injection $R_Y\xrightarrow{\iota_\mu}(R_Y \otimes_{R_P} R_Y)$
$\to (R_Y \otimes_{R_P} R_Y)/I$.
Hence $p_\mu'$ is the identity map as a set, and is bijective.
Further since
\[
(ab\otimes 1)-(1\otimes ab)=
(b\otimes 1)(a\otimes 1 -1\otimes a)+
(1\otimes a)(b\otimes 1 - 1\otimes b),
\]
every element in $\{a \otimes 1 - 1 \otimes a\}_{a\in R_Y}$ is generated by
$\{\bby_i \otimes 1 - 1 \otimes \bby_i\}_{i=0, \dots, r-1}$.
\end{proof}

We consider the annihilator of the kernel of $\mu$,
\[
\operatorname{Ann}_{R_Y^e}(\Ker \mu) :=\{ z\in R_Y^e\, |\, z\cdot \Ker \mu = 0\}.
\]
For an element $\sum a_i \otimes b_i \in\operatorname{Ann}_{R_Y^e}(\Ker \mu)$,
$(c \otimes 1-1 \otimes c)\cdot \sum a_i \otimes b_i =0$ or $\sum ca_i \otimes b_i =\sum a_i \otimes c b_i$.

\begin{Lemma}There is a natural embedding $\varphi\colon \operatorname{Ann}_{R_Y^e}(\Ker \mu)\to
\operatorname{Hom}_{R_Y}(\omega_{R_Y/R_P}, R_Y)$.
\end{Lemma}

\begin{proof}Since $\sum ca_i \otimes b_i =\sum a_i \otimes c b_i$ for an element $\sum a_i \otimes b_i \in\operatorname{Ann}_{R_Y^e}(\Ker \mu)$ and $c\in R_Y$, it should be defined
$\varphi(\sum_i a_i \otimes b_i)(\rho)$ $=\sum_i\rho(a_i) b_i\in R_Y$ for $\rho\in \omega_{R_Y/R_P}$;
for every $s\in R_Y$, we have $s \rho(a) b=\rho(a) s b=\rho(s a)b$.
\end{proof}

\begin{Proposition}[{\cite[Corollary H.20]{Kunz2004}}]\label{2pr:A_traceH20}
Suppose $R_Y/R_P$ has a trace.
Then $\varphi$ induces a bijection between the set of all traces of $R_Y/R_P$ and the set of all generators of the $R_Y$-module $\operatorname{Ann}_{R_Y^e}(\Ker \mu)$:
Each trace $\tau \in \omega_{R_Y/R_P}$ is mapped to the unique element $\Delta_\tau:=\sum_{i=0}^{r-1} \hby_i\otimes\bby_i \in \operatorname{Ann}_{R_Y^e}(\Ker \mu)$ associated with $\tau$ such that $\sum_{i=0}^{r-1}\tau( \hby_i)\bby_i =1$.
Furthermore

\begin{enumerate}\itemsep=0pt
\item[$1.$] $\Delta_\tau$ generates the $R_Y$-module $\operatorname{Ann}_{R_Y^e}(\Ker \mu)$, and $\{\hby_1, \dots, \hby_r\}$ is the dual basis of the basis $\{\bby_1, \dots, \bby_r\}$ of $R_X$ with respect to $\tau$; i.e.,
\[
\tau(\hby_i\bby_j)=\delta_{i,j}, \qquad i,j = 1, \dots, r.
\]
\item[$2.$] If $\sum_{i=0}^{r-1} \hby'_i \otimes \bby_i$ generates the $R_Y$-module $\operatorname{Ann}_{R_Y^e}(\Ker \mu)$, and if $\tau' \in \omega_{R_Y/R_P}$ is a linear form with $\sum_{i=0}^{r-1}\tau'( \hby_i')\bby_i =1$, then $\tau'$ is a trace of $R_Y/R_P$ and $\Delta_{\tau'}=\sum_{i=0}^{r-1} \hby'_i \otimes \bby_i$ is associated with the trace $\tau'$; hence $\{\hby'_1, \dots, \hby'_r\}$ is the dual basis of the basis $\{\bby_1, \dots, \bby_r\}$ of $R_X$ with respect to $\tau'$.

\item[$3.$] For each trace $\tau$ of $R_Y/R_P$,
\[
\tau_{R_Y/R_P} = \mu(\Delta_\tau)\circ \tau.
\]
\end{enumerate}
\end{Proposition}

\begin{proof}Let us consider $1 \in \operatorname{Hom}_{R_Y}(\omega_{R_Y/R_P}, R_Y)$, and consider the inverse image of $\{1\}$, $\varphi^{-1}(1)\in \operatorname{Ann}_{R_Y^e}(\Ker \mu)$, which is denoted by $\Delta_\tau$. By using the basis $\{\bby_i\} \subset R_Y$, we can express it as $\Delta_\tau=\sum_{i}\hby_i \otimes \bby_i$ for a certain $\{\hby_i\}$.
Thus $\varphi(\sum_i \hby_i \otimes \bby_i)(\tau)=\sum_i \tau(\hby_i) \bby_i=1$.
The fact that $\sum_i \tau(\hby_i) \bby_i=1$ implies that $\bby_j=\bby_j\varphi(\sum_i \hby_i \otimes \bby_i)(\tau)=\varphi(\sum_i \hby_i \otimes \bby_i)(\bby_j\circ \tau)=\sum_i \tau(\hby_i \bby_j)\bby_i$.
Thus $\tau(\hby_i \bby_j)=\delta_{i j}$.
Hence the set $\{\hby_i\}$ is a dual basis of $R_Y$ with respect to~$\tau$.

Let us consider the case $\operatorname{Ann}_{R_Y^e}(\Ker \mu)=R_Y\big(\sum_i \hby_i' \otimes \bby_i\big)$ with $\sum_{i=0}^{r-1}\tau'( \hby_i')\bby_i =1$ for ${\tau' \in \omega_{R_Y/R_P}}$. Then by the above arguments, obviously $\varphi(\sum_i \hby_i' \otimes \bby_i)=1$. we obtain~(2).

Lemma~\ref{2lm:trace_1.2} shows that $\tau_{R_Y/R_P}=\sum_i \bby_i\circ \bby_i^*=\sum_i \hby_i \bby_i\circ \tau$ $=\mu(\Delta_\tau)\circ \tau$.
\end{proof}

\subsection{Numerical and Weierstrass semigroup}\label{2sc:WSG}

This subsection is on numerical and Weierstrass semigroups based on \cite{ADGS2016, KMP13}.
An additive sub-monoid of the monoid, non-negative integers $\NN_0$ is called a~\emph{numerical semigroup} if its complement in $\NN_0$ is a finite set.
In this subsection, we review the numerical semigroups associated with algebraic curves.

Let $X$ be a smooth (complex projective) curve of genus g, and $\cM(X)$ be the set of the meromorphic functions on $X$.
For a point $P \in X$,
\[
 H(X,P):= \{n \in \NN_0\, |\, \text{there exists } f \in \cM(X)
 \text{ such that } (f)_\infty = n P \}
\]
forms a numerical semigroup by the Riemann--Roch theorem, which is called the Weierstrass semigroup of the point $P$.
If the Weierstrass gap sequence $H^\fc(X,P):=\NN_0 \setminus H(X, P)$ differs from the set $\{1, 2, \dots, g\}$, we say that $P$ is a Weierstrass point of~$X$~\cite{FarkasKra}.

In this paper, we consider a \emph{pointed curve}, a pair $(X,P)$ with $P$ a point of the curve $X$ with the Weierstrass semigroup $H(X,P)$.

In general, a numerical semigroup $H$ has a unique (finite) minimal set of generators, $M=M(H)$, ($H=\langle M\rangle$) and the finite cardinality $g$ of $H^\fc=\NN_0\setminus H$; $g$ is the genus of $H$ or $H^\fc$ and $H^\fc$ is called a~gap-sequence.
For example,
\begin{gather*}
H^\fc=\{1, 2, 4, 5\}, \qquad \text{for} \ H=\langle 3,7,8\rangle, \\
H^\fc=\{1,2,3,4,6,8,9,13\}, \qquad \text{for} \ H=\langle 5,7,11\rangle,\qquad
\text{and} \\
H^\fc=
\{1,2,3,4,5,7,8,9,10,11,17,23\}, \qquad \text{for}
\ H=\langle6,13,14,15,16\rangle.
\end{gather*}
We let $r_{\min}(H)$ be the smallest positive integer of $M(H)$, which is referred the multiplicity of $H$.
We call the semigroup $H$ \emph{an $r_{\min}(H)$-semigroup}, so that $\langle 3,7,8\rangle$ is a~$3$-semigroup and $\langle 6,13,14,15,16\rangle$ is a $6$-semigroup.
Let $N(i)$ and $N^\fc(i)$ be the $i$-th ordered element of $H=\{N(i)\, |\, i \in \NN_0\}$ and $H^\fc=\{N^\fc(i)\, |\, i = 0, 1, \dots, g-1\}$ satisfying $N(i) < N(i+1)$ and $N^\fc(i) < N^\fc(i+1)$, respectively.
The Schubert index of the set $H^\fc$ is
\begin{equation*}
\alpha(H) :=\{\alpha_0(H), \alpha_1(H), \dots, \alpha_{g-1}(H)\},
\end{equation*}
where $\alpha_i(H) := N^\fc(i) - i -1$.
Moreover,
\[
\alpha(\langle 3,7,8\rangle) =\big\{0^2, 1^2\big\}
\qquad\text{and}\qquad
\alpha(\langle 6,13,14,15,16\rangle) =\big\{0^5, 1^5, 6, 11\big\}.
\]
Further the conductor $c_H$ of $H$ is defined by the minimal natural number satisfying ${c_H + \NN_0 \subset H}$.
The number $c_H-1$ is known as the Frobenius number, which is the largest element of $H^\fc$.

By letting the row lengths be $\Lambda_i = \alpha_{g-i}(H)+1$, $i< g$, we have the Young diagram of the semigroup, $\Lambda:=(\Lambda_1, \dots, \Lambda_g)$, $\Lambda_i \ge \Lambda_{i+1}$.
The Young diagram $\Lambda$ is a partition of $\sum_i\Lambda_i$.
We say that such a Young diagram is associated with the numerical semigroup.
If for a given Young diagram $\Lambda$, we cannot find any numerical semigroup $H$ such that $\Lambda_i = \alpha_{g-i}(H)+1$, we say that $\Lambda$ is not associated with the numerical semigroup.
It is obvious that in general, the Young diagrams are not associated with the numerical semigroups.
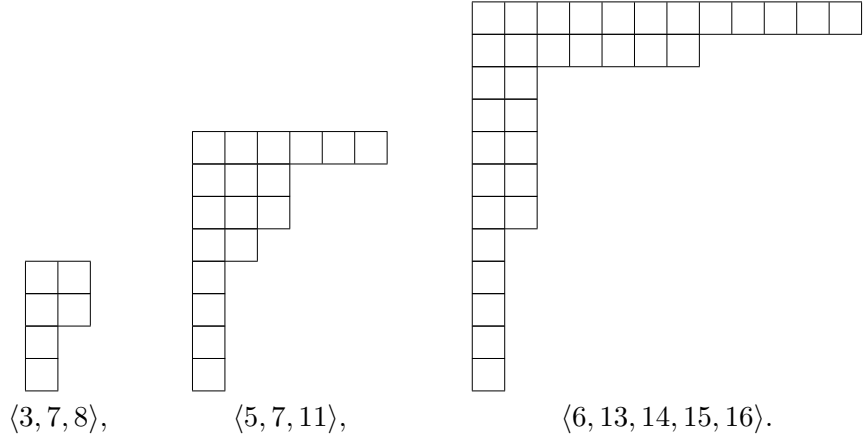
\begin{gather*}
\begin{matrix}
\yng(2,2,1,1)\qquad &
\yng(6,3,3,2,1,1,1,1)\qquad &
\yng(12,7,2,2,2,2,2,1,1,1,1,1)
\\
 \langle3, 7, 8\rangle, \qquad &
 \langle5,7,11\rangle, \qquad
 &\langle6,13,14,15,16\rangle.\\
\end{matrix}
\end{gather*}
The Young diagram and the associated numerical semigroup are called symmetric if the Young diagram is invariant under reflection across the main diagonal.
It is known that the numerical semigroup is symmetric if and only if $2g-1$ occurs in the gap sequence.
It means that if $c_H=2g$, $H$ is symmetric.

We obviously have the following proposition:
\begin{Proposition}\label{2pr:N(n)}
The following hold:
\begin{enumerate}\itemsep=0pt
\item[$(1)$] $N(n) - n \le g$ for every $n\in \NN_0$,

\item[$(2)$]
$N(n) - n = g$ for $N(n)\ge c_X=N(g)$,

\item[$(3)$] $N(n) - n < g$ for $0\le N(n)< c_X$,

\item[$(4)$]
$\#\{n \, |\, N^\fc(n) \ge g\}=\#\{n \, |\, N(n) < g\}$,

\item[$(5)$] for $N(i) < N^\fc(j)$, $N^\fc(j)-N(i)\in H^\fc$, and

\item[$(6)$] when $H$ is symmetric, $c_H=N(g)=2g$ and
$c_H - N(i)-1=N^\fc(g-i-1)$ for $0\le i\le g-1$.
\end{enumerate}
\end{Proposition}

\begin{proof} (1)--(3) and (5) are obvious.
By noting $\# H^\fc =g$, (4) means that what is missing must be filled later for~$H^\fc$. (6) is left to~\cite{ADGS2016}.
\end{proof}

In this paper, we mainly consider the $r$-numerical semigroup,~$H$.
We introduce the tools as follows:

\begin{Definition}\label{2df:NSG1}\quad
\begin{enumerate}\itemsep=0pt

\item
Let $\ZZ_r:= \{ 0, 1, 2, \dots, r-1\}$ and $\ZZ_r^\times := \ZZ_r \setminus\{0\}$.

\item
Let $\tfe_i :=\min\{h \in H \, |\, i \equiv h \ \text{mod}\ r\}$,
$i \in \ZZ_r$.

\item Let $\tfE_H:=\{\tfe_i \, |\, i \in \ZZ_r\}$ be the standard basis of $H$.
Further we define the ordered set $\fE_H:=\{\fe_i \in \tfE_H \, | \, \fe_i < \fe_{i+1}\}$, and $\fE_H^\times := \fE_H\setminus \{0\}$, e.g., $\fe_{0}=\tfe_0=0$.

\item Let $\fe_{\ell,i}^*$ be the element in $\fE_H$ such that $\fe_{\ell,i}^*= \fe_\ell - \fe_i$ modulo $r$.

\item Let $\lH^\fc:= H^\fc \bigcup (-\NN)$, where $\NN:=\NN_0\setminus\{0\}$.

\item The Ap\'ery set $\operatorname{Ap}(H, n)$ for a positive integer $n$ is defined by
\[
\operatorname{Ap}(H,n):=\{s \in H\, |\, s-n\not\in H\}.
\]
\end{enumerate}
\end{Definition}

Since it is obvious that $\operatorname{Ap}(H, r)=\fE_H=\tfE_H$ as a set, the standard basis is sometimes defined by the Ap\'ery set $\operatorname{Ap}(H, r)$.

We have the following elementary but essential results:

\begin{Lemma}\label{2lm:NSG1}
For $a \in \NN_0$, we define
\[
[a]_r:= \{ a + k r \, | \, k \in \NN_0\}, \qquad
\overline{[a]}_r^\fc:=
\{ a - k r \, | \, k \in \NN\}, \qquad
[a]_r^\fc:=\overline{[a]}_r^\fc\cap \NN.
\]
\begin{enumerate}\itemsep=0pt
\item[$1.$] We have the following decomposition:
\begin{enumerate}\itemsep=0pt
\item[$(a)$]
$ \NN_0 = \bigoplus_{i\in \ZZ_r} [i]_r$,

\item[$(b)$]
$H = \bigoplus_{i\in \ZZ_r} [\fe_i]_r$,

\item[$(c)$]
$\lH^\fc = \bigoplus_{i\in \ZZ_r} \overline{[\fe_i]}_r^\fc$, $H\bigcup \lH^\fc=\ZZ$,

\item[$(d)$]
$H^\fc = \bigoplus_{i\in \ZZ_r} [\fe_i]_r^\fc
= \bigoplus_{i\in \ZZ_r^\times} [\fe_i]_r^\fc$,\ $H\bigcup H^\fc=\NN_0$,
\end{enumerate}

\item[$2.$] for every $x_i \in [\fe_i]_r$ $(i\in \ZZ_r)$,
\[
\{x_i \text{ modulo } r\, |\, i \in \ZZ_r\}=\ZZ/r \ZZ,
\]
especially for $x \in [\tfe_i]_r$, $x = i$ modulo $r$, and

\item[$3.$] $\fe_{\ell,\ell}^*=\fe_0=0$, and $\fe_{\ell,0}^* =\fe_{\ell}$.
\end{enumerate}
\end{Lemma}

\begin{proof}
(1a), (2) and (3) are apparent.
From the definition of $\fE_H$, $H= \{ \fe_i + k r \, |\, i \in \ZZ_r,\, k \in \NN_0\}$.
For $i \neq j$, $[\fe_i]_r\cap [\fe_j]_r = \varnothing$ and thus we have the decomposition in (1b).
Since $H^\fc = \NN_0 \setminus H$, we have (1d) and (1c) noting (1a).
\end{proof}

The following is obvious:

\begin{Lemma}\label{2lm:rs_iris}
For the generators $r$ and $s$ in the numerical semigroup $H$, there are positive integers $i_s$ and $i_r$ such that
$ i_s s- i_r r = 1$.
\end{Lemma}

We remark that $\lH^\fc$ determines the structure of the differentials on a certain curve $X$ in Theorem~\ref{2lm:nuI_hX} and the embedding of the curve into the universal Grassmannian manifold as in \cite[p.~46]{SegalWilson}.

A numerical semigroup $H$ is said to be Weierstrass if there exists a pointed curve $(X, P)$ such that $H=H(X,P)$.
Hurwitz posed the problem of whether any numerical semigroup $H$ is Weierstrass.
The question was revived in the 1980s, viewed as the question of deformations of a reduced complex curve singularity $(X_0,\infty )$.
Buchweitz and Greuel showed a counterexample.

\section[Weierstrass canonical form and Weierstrass curves (W-curves)]{Weierstrass canonical form and Weierstrass curves\\ (W-curves)}\label{sec:WCF}

\subsection{Weierstrass canonical form}

In this subsection, we now review the ``Weierstrass canonical form'' (``Weierstrass normal form'') based on \cite{KM2020, KMP16, KMP19}, which is a generalization of Weierstrass' standard form for elliptic curves.
This form originated from Abel's insight, and Weierstrass investigated its primitive property \cite{Wei67, WeiWIV}.
Baker \cite[Chapter~V, Sections~60--79]{Baker97} gave a complete review, proof, and examples of the theory.
Here we refer to Kato \cite{CoppensKato, Kato80}, who also produces this representation with proof.

\begin{Proposition}[{\cite{CoppensKato, Kato80}}] \label{2pr:WCF}
For a pointed curve $(X,\infty )$ with Weierstrass semigroup $H_X:=H(X,\infty)$
for which $r_{\min}(H_X)=r$, and $\fe_i\in \fE_{H_X}$, $i \in \ZZ_r^\times$, in Definition~{\rm \ref{2df:NSG1}}, and we let
$s:=\min_{i\in \ZZ_r^\times}\{\fe_i \in \fE_{H_X}\, |\, (\fe_i, r)=1\}$ and $s=\fe_{\ell_s}$. $(X,\infty)$ is defined by an irreducible equation,
\begin{equation}
f_X(x,y)= 0,
\label{2eq:WCF1a}
\end{equation}
for a polynomial $f_X\in \CC[x,y]$ of type,
\begin{equation}
f_X(x,y):=y^r + A_{1}(x) y^{r-1}
+ A_{2}(x) y^{r-2}
+\cdots + A_{r-1}(x) y
 + A_{r}(x),
\label{2eq:WCF1b}
\end{equation}
where the $A_i(x)$'s are polynomials in $x$, $A_0=1$,
\begin{equation*}
A_i = \sum_{j=0}^{\lfloor i s/r\rfloor} \lambda_{i, j} x^j=:\prod_{j=1}^{\lfloor i s/r\rfloor}a_{i,0}(x-a_{i,j}),
\end{equation*}
and $\lambda_{i,j} \in \CC$, $\lambda_{r,s}=-1$.
\end{Proposition}

\begin{proof}We let $\BfE_{H_X} :=\{\fe_1, \fe_2, \dots, \fe_{r-1}\}\setminus \{s\}$, $\BfE_{H_X} =:\{\Bfe_2, \dots, \Bfe_{r-1}\}$, ($\Bfe_i<\Bfe_j$ for $i<j$), and $\Bfe_1=s$.
Let $\tfy_{\Bfe_i}$ be a meromorphic function on $X$ whose only pole is $\infty$ with order $\Bfe_i$, $i\in \ZZ_r^\times$, taking $y_s:=\tfy_{s}=\tfy_{\Bfe_{1}}$ and $1=\tfy_{\Bfe_{0}}$.
From the definition of $X$, we have, as $\CC$-vector spaces,
\begin{equation}
\bH^0(X,\cO_X(*\infty))
= \sum_{i=0}^{r-1}\sum_{j=0} \CC x^j \tfy_{\Bfe_i}.
\label{2eq:R_Xdecompose1}
\end{equation}
Let the affine ring $R_X$ of $X$ be defined by $\bH^0(X,\cO_X(*\infty))$, i.e., $R_X:=\bH^0(X,\cO_X(*\infty))$.
Thus for every $\Bfe_j\in \BfE_{H_X}$, $j=2,3,\dots,r-1$, we obtain the following equations,
\begin{gather}
\begin{cases}
y_{s} \tfy_{\Bfe_2} = A_{2,0} +
 A_{2,1} y_{s}+
 A_{2,2} \tfy_{\Bfe_{2}}+
 \cdots +
 A_{2,r-2} \tfy_{\Bfe_{r-2}}+
 A_{2,r-1} \tfy_{\Bfe_{r-1}}, \\
y_{s} \tfy_{\Bfe_3}= A_{3,0} +
 A_{3,1} y_{s}+
 A_{3,2} \tfy_{\Bfe_{2}}+
 \cdots +
 A_{3,r-2} \tfy_{\Bfe_{r-2}}+
 A_{3,r-1} \tfy_{\Bfe_{r-1}}, \\
\cdots\cdots\cdots\cdots\cdots\cdots\cdots\cdots\cdots\cdots\cdots\cdots\cdots\cdots\cdots\cdots\cdots\cdots\cdots\cdots\cdots\\
y_{s} \tfy_{\Bfe_{r-1}} = A_{{r-1},0} +
 A_{r-1,1} y_{s}+
 A_{r-1,2} \tfy_{\Bfe_{2}}+
 \cdots +
 A_{r-1,r-2} \tfy_{\Bfe_{r-2}}+
 A_{r-1,r-1} \tfy_{\Bfe_{r-1}}, \\
\end{cases}
\label{2eq:WCF2}
\\
y_{s}^2 = A_{1,0} +
 A_{1,1} y_{s}+
 A_{1,2} \tfy_{\Bfe_{2}}+
 \cdots +
 A_{1,r-2} \tfy_{\Bfe_{r-2}}+
 A_{1,r-1} \tfy_{\Bfe_{r-1}},
\label{2eq:WCF3}
\end{gather}
where $A_{i,j} \in \CC[x]$.

When $r = 2$, (\ref{2eq:WCF1a}) equals (\ref{2eq:WCF3}).
We assume that $r>2$ and then (\ref{2eq:WCF2}) is reduced to $(r-2)$ linear equations,
\begin{gather}
\begin{pmatrix}
A_{2,1} - y_s & A_{2,2} & \cdots & A_{1,r-1}\\
A_{3,1} & A_{3,2} - y_s& \cdots & A_{3,r-1}\\
\vdots & \vdots & \ddots & \vdots \\
A_{r-1,1} & A_{r-1,2} & \cdots & A_{r-1,r-1} - y_s
\end{pmatrix}
\begin{pmatrix}
\tfy_{\Bfe_2}\\
\tfy_{\Bfe_3}\\
\vdots\\
\tfy_{\Bfe_{r-1}}
\end{pmatrix}
=-
\begin{pmatrix}
A_{2,0} + A_{2, 1} y_s\\
A_{3,0} + A_{3, 1} y_s\\
\vdots\\
A_{r-1,0} + A_{r-1, 1} y_s
\end{pmatrix}.\label{2eq:WCF4}
\end{gather}

One can check that the determinant of the matrix on the left-hand side of~(\ref{2eq:WCF4}) is not equal to zero by computing the order of pole at $\infty$ of
the monomials $B_{i} y_s^{r-2-i}$ in the expression,
\begin{align*}
P(x, y_s):={}& \left|
\begin{matrix}
A_{2,1} - y_s & A_{2,2} & \cdots & A_{2,r-1}\\
A_{3,1} & A_{3,2} - y_s& \cdots & A_{3,r-1}\\
\vdots & \vdots & \ddots & \vdots \\
A_{r-1,1} & A_{r-1,2} & \cdots & A_{r-1,r-1} - y_s\\
\end{matrix}
\right|\\
={}& y_s^{r-2} + B_1 y_s^{r-3} +\cdots
+ B_{r-3} y_s + B_{r-4},
\end{align*}
which is $s (r -2 -i) + r \cdot\deg_x B_i$ by letting $\deg_x h(x)$ be the degree of $h$ with respect to $x$.
The fact that $(r,s)=1$ shows that $s (r -2 -i) + r \cdot\deg_x B_i\neq s (r -2 -j) + r \cdot\deg_x B_j$ for $i\neq j$.

Hence by solving equation (\ref{2eq:WCF4}), we have
\begin{equation}
\tfy_{\Bfe_i} = \frac{Q_{i}(x, y_s)}{P(x, y_s)},
\label{2eq:WCF5}
\end{equation}
where $Q_{i}(x, y_s) \in \CC[x,y_s]$ and a polynomial of order at most $r-2$ in $y_s$.
Note that the equations~(\ref{2eq:WCF5}) are not algebraically independent in general but in any case the function field of the curve can be generated by
some of these $\tfy_{\Bfe_j}$'s, and its affine ring $R_X$ can be given by a quotient ring of $\CC[x,y_{r_2}, y_{r_3}, \dots, y_{r_{m_X}}]$ for $i_j \in M_X:=M(H_X)$, where $M_X =\{ r_1, r_2, \dots, r_{m_X}\} \subset \NN^{{m_X}}$ with the conditions that mutually coprime, $(r_1, \dots, r_{m_X}) = 1$, $r_1=r$, $r_2=s$, and $r_i<r_j$ for $2<i<j$, is a minimal set of generators for~$H_X$. Here the cardinality of the generator~$M_X$ of~$H_X$ is $m_X (< r)$.

By putting (\ref{2eq:WCF5}) into (\ref{2eq:WCF3}), we obtain (\ref{2eq:WCF1b}) if it is irreducible.
If it is reducible, $f_X(x,y_s)$ is decomposed to polynomials whose degree is less than $r$ with respect to $y_s$. However the relation $(r,s)=1$ shows that there does not exist such monic polynomials due to the order of the singularity at $\infty$.

Further the order of the singularity of $A_i y_s^{r-i}$ is given by
$s(r-i)+r\deg_x A_i$.
The cases satisfying that $s(r-i)+r\deg_x A_i= s(r-j)+r\deg_x A_j$ mean that $i=j$ or $(i,j)=(0,r), (r,0)$ due to $(r,s)=1$. Hence $r\deg_x A_r = s$. For $i=1, \dots, r-1$, $s(r-i)+r\deg_x A_i<rs$ leads that $\deg_x A_i<s i/r$.
\end{proof}

\begin{Remark}\label{2rm:2.9}
Let us call the curve in Proposition \ref{2pr:WCF} a {\emph{Weierstrass curve}} or a~{\emph{W-curve}} in this paper.
The Weierstrass canonical form characterizes the W-curve, which has only one infinity point~$\infty$. The infinity point $\infty$ is a Weierstrass point if $H_X^\fc=H^\fc(X,\infty)=\{N^\fc(i)\}$ differs from $\{1, 2, \dots, g\}$.
Since every compact Riemann surface of the genus, $g(>1)$, has a Weierstrass point whose Weierstrass gap sequence has genus~$g$ \cite{ACGH85, FarkasKra}, it characterizes the behavior of the meromorphic functions at the point, and thus there is a W-curve which is bi-rationally equivalent to the compact Riemann surface.

Further Proposition \ref{2pr:WCF} is also applicable to a pointed compact Riemann surface $(Y, P)$ of genus $g$ whose point $P$ is an ordinary point rather than the Weierstrass point; its Weierstrass gap sequence at $P$ is $H^\fc(Y, P)=\{1, 2, \dots, g\}$.
Even for the case, we find the Weierstrass canonical form $f_X$ and the W-curve $X$ with $H_X^\fc=\{1, 2, \dots, g\}$ which is bi-rational to $Y$.
\end{Remark}

\begin{Remark}\label{2rm:g_Xc_X}
Let $R_{X^\circ}^\circ := \CC[x,y]/(f_X(x,y))$ for (\ref{2eq:WCF1a}) and its normalized ring be $R_X^\circ$ if $X^\circ:=\Spec R_{X^\circ}^\circ$ is singular.
$R_X^\circ$ is the coordinate ring of the affine part of $X\setminus\{\infty\}$ and we identify~$R_X^\circ$ with $R_X=\bH^0(X,\cO_X(*\infty))$.
Then the quotient field $\CC(X):=\cQ(R_X)$ of $R_X$ is considered as an algebraic function field on $X$ over $\CC$.

By introducing $R_\PP:=\bH^0(\PP, \cO_\PP(*\infty))=\CC[x]$ and its quotient field $\CC(x):=\cQ(R_\PP)$, $\cQ(R_X)$ is considered a finite extension of $\cQ(R_\PP)$.
We regard~$R_X$ as a finite extended ring of $R_\PP$ of rank~$r$, e.g., $R_{X^\circ}^\circ=R_\PP[y]/(f_X(x,y))$ as mentioned in Section~\ref{2sc:Covering} \cite{Kunz2004}.
Further as we will mention in Section~\ref{2ssc:R_Pmodule_R_X}, $\tfy_{\Bfe_{i}}$ in the proof of Proposition \ref{2pr:WCF} is the standard basis of $R_X$ as an $R_\PP$-module, and the matrix in~(\ref{2eq:WCF4}) are naturally interpreted as the $R_\PP$-module.

For the local ring $R_{X,P}$ of $R_X$ at $P \in X$, we have the ring homomorphism,
$\varphi_{P}\colon R_X\to R_{X,P}$.
We note that~$R_{X, \infty}$ plays crucial roles in the Weierstrass canonical form.
We let $M_X =\{ r_1, r_2, \dots, r_{m_X}\}$ be the minimal generator of the numerical semigroup $H_X=H(X, \infty)$ appearing in the proof of Proposition~\ref{2pr:WCF}.
The Weierstrass curve admits a local cyclic $\fC_r=\ZZ/r\ZZ$-action at $\infty$, cf.\ Section~\ref{2ssc:MCurve}.
The genus of~$X$ is denoted by~$g_X$, briefly $g$ and the conductor of~$H_X$ is denoted by $c_X:=c_{H_X}$; the Frobenius number $c_X-1$ is the maximal gap in~$H_X$, i.e., $c_H=N^\fc(g-1)+1$.
We let $\lH_X^\fc:=\ZZ\setminus H_X$.
\end{Remark}

\subsubsection[Projection from X to P]{Projection from $\boldsymbol{X}$ to $\boldsymbol{\PP}$}\label{2sssc:X_to_PP}

There is the natural projection,
\begin{equation*}
\varpi_{r}\colon \ X \to \PP, \qquad
\varpi_{r}(x,y_{r_2}, \dots, y_{r_{m_X}}) = x=y_r,
\end{equation*}
such that $\varpi_{r}(\infty) = \infty\in \PP$.

Let $\{y_{\bullet}\}:=\{y_{s}=y_{r_2}, y_{r_3}, \dots, y_{r_{m_X}}\}$ and $\CC[x, y_{\bullet}]:=\CC[x, y_{s}=y_{r_2}, y_{r_3}, \dots, y_{r_{m_X}}]$.

\subsubsection{Symmetric and non-symmetric Weierstrass curves (W-curves)}\label{2sssec:SymmWcurve}

We also investigate the W-curves whose Weierstrass semigroups $H_X$ are symmetric and non-symmetric, which are called {\emph{symmetric Weierstrass curve}} or {\emph{symmetric W-curve}}, and {\emph{non-symmetric Weierstrass curve}} or {\emph{non-symmetric W-curve}} respectively in this paper.

\subsection{The monomial curves and W-curves}\label{2ssc:MCurve}

This subsection shows the monomial curves and their relation to W-curves based on \cite{KM2020, KMP16, KMP19}.

For a given W-curve $X$ with the Weierstrass semigroup $H=H_X$, and its generator $M_X=\{r=r_1, r_2, \dots, r_{m_X}\}$, the behavior of singularities of the elements in $R_X$ at $\infty$ is described by a monomial curve $X_H^Z$.
For the numerical semigroup $H=\langle M_X\rangle$, the \emph{numerical semigroup ring}~$R_H$ is defined as $R_H:=\CC[ z^{r_1}, z^{r_2}, \dots, z^{r_{m_X}}]$.

Following a result of Herzog's \cite{Herzog70}, we recall the well-known proposition for a polynomial ring $\CC[Z]:=\CC[Z_{r_1}, Z_{r_2}, \dots, Z_{r_{m_X}}]$.
\begin{Proposition} \label{2pr:R_H}
For the $\CC$-algebra homomorphism $ \tvarphi^Z_H\colon \CC[Z] \to R_H$, $Z_a\mapsto z^a$, the kernel of $\tvarphi^Z_H$ is generated by certain binomials $f^H_i \in \CC[Z]$, $i=1,2,\dots,k_X$, for a positive integer $k_X$, $m_X-1\le k_X< \infty$, i.e., $\ker\tvarphi^Z_H=\big(f^H_1, f^H_2, \dots, f^H_{k_X}\big)$, and
\[
R_H \simeq \CC[Z]/ \ker\tvarphi^Z_H=:R^Z_H.
\]
\end{Proposition}

\begin{proof}This follows from a result of Herzog's \cite{Herzog70}.
There are the multiplicative group actions ($\GG_\fm$-actions) on $Z_{r_i}$'s, whereas~$R_H$ is invariant for the action.
It means that the number of generators of $\ker\tvarphi^Z_H$ is determined, i.e.,~$k_X$.
The relation in the $\ker\tvarphi^Z_H$ is reduced to a binary one.
\end{proof}

We call $R^Z_H=\CC[Z]/ \ker\tvarphi^Z_H$ a \emph{monomial ring}.
Sending $Z_r$ to $1/x$ and $Z_{r_i}$ to $1/y_{r_i}$, the monomial ring $R^Z_H$
determines the structure of gap sequence of $X$ at $\infty$ \cite{Herzog70, Pinkham74}.
Bresinsky showed that $k_X$ can be any finitely large number if $m_X>3$ \cite{Bresinsky}.

Let $X_H:=\Spec R_H^Z$,which we call a \emph{monomial curve}.
We also define the ring isomorphism on $R_H^Z$ induced from $\tvarphi^Z_H$, which is denoted by $\varphi^Z_H$,
\begin{equation*}
\varphi^Z_H\colon \ \CC[Z]/ \ker\tvarphi^Z_H=R^Z_H \to R_H.
\end{equation*}

Further, we let $\{Z_{\bullet}\}:=\{Z_{r_2}, Z_{r_3}, \dots, Z_{r_{m_X}}\}$, and $\CC[Z_r, Z_{\bullet}]:=\CC[Z_r, Z_{r_2}, Z_{r_3}, \dots, Z_{r_{m_X}}]$.
A~monomial curve is an irreducible affine curve with $\GG_\fm$-action, where $\GG_\fm$ is the multiplicative group of the complex numbers; $Z_a \mapsto g^{a} Z_a$ for $g\in \GG_\fm$, and it induces the action on the monomial ring $R_H^Z$.

The following action of the cyclic group of order $r$ plays a crucial role in this paper.

\begin{Lemma} \label{2lm:fC_r_action}
The cyclic group $\fC_r$ of order $r$ acts on the monomial ring $R_H^Z$;
the action of the generator $\hzeta_r \in \fC_r$ on $Z_a$ is defined by sending $Z_a$ to $\zeta_r^{a}Z_a$, where $\zeta_r$ is a primitive $r$-th root of unity.
By letting $\fr_i^*:=(r, r_i)$, $\fr_i := r/\fr_i^*$, and $\Bfr_i:=r_i/\fr_i^*$, the orbit of $Z_{r_i}$ forms $\fC_{\fr_i}$; especially for the case that $(r, r_i) = 1$, it recovers $\fC_r$.

Thus in $R_H^Z$, we have the following identities:
\begin{equation}
f_H^{(j)}(Z_r, Z_{r_j})= 0, \qquad f_H^{(j)}:=Z_{r_j}^{\fr_j}-Z_r^{\Bfr_j},
\qquad j = 2, \dots, m_X.\label{2eq:f_H^j}
\end{equation}
\end{Lemma}

For example, the case $M_X=\{3,7,8\}$ provides these elements $\big\{f^H_1,f^H_2, f^H_3\big\}$ are given by the $2 \times 2$ minors of
$\left|\begin{smallmatrix}
 Z_3^2 & Z_7 & Z_8 \\
 Z_7 & Z_8 & Z_3^3\\
\end{smallmatrix} \right|$.
There is a cyclic group $\fC_3:=\{\zeta_3^a\, | \, a=0,1,2\}$ action on $R_H$ as a $\GG_\fm$ action. Due to the relation, there are other possibilities which are given by the $2 \times 2$ minors of
$
\left|\begin{smallmatrix}
 Z_3^2 & \zeta_3^a Z_7 & \zeta_3^{2a}Z_8 \\
\zeta_3^a Z_7 & \zeta_3^{2a}Z_8 & Z_3^3\\
\end{smallmatrix} \right|$ for $a=0,1,2$.
It means that $f^H_i$ is unique up to the $\GG_\fm$ action.
On the other hand, $f_H^{(2)}=Z_7^3-Z_3^7$ and $f_H^{(3)}=Z_8^3-Z_3^8$ for~(\ref{2eq:f_H^j}).

There are non-negative integers $h_j^{(i\pm)}$ such that
\begin{equation}
f^H_i =
\bigg(
\prod_{j=2} Z_{r_j}^{h_j^{(i+)}}
\bigg)
-\bigg(
\prod_{j=1} Z_{r_j}^{h_j^{(i-)}}
\bigg),
\label{2eq:f^H_iZ}
\end{equation}
where, in the first term, $Z_{r_1}$ does not exist because $(r_1, r_2, \dots, r_{m_X})=1$.

Corresponding to the standard basis of $H_X$ in Definition \ref{2df:NSG1}, we find the monic monomial $\fZ_{\fe_i} \in R_H^Z$ such that $\varphi^Z_H(\fZ_{\fe_i}) = z^{\fe_i}$, and the standard basis $\{\fZ_{\fe_i}\, |\, i \in \ZZ_r\}$; $\fZ_{\fe_0}=1$.

\begin{Lemma}\label{2lm:Z_standardbasis}
The $\CC[Z_r]$-module $R_H^Z$ is given by
\[
R_H^Z=\CC[Z_r]\oplus \CC[Z_r]\fZ_{\fe_1}\oplus \cdots
\oplus \CC[Z_r]\fZ_{\fe_{r-1}},
\]
and thus $\fZ_H:=\{\fZ_{\fe_0}=1, \fZ_{\fe_1}, \dots, \fZ_{\fe_{r-1}}\}$ is the basis of the $\CC[Z_r]$-module~$R_H^Z$.
Then there is a~monomial $b_{i j k}\in \CC[Z_r]$ such that
\[
\fZ_{\fe_i}\fZ_{\fe_j}
=\sum_{k\in \ZZ_r}b_{i j k}\fZ_{\fe_k}.
\]

Further, there are elements $\fZ\in R_H^Z$ and $\fZ_{\fe_i}^*\in R_H^Z$, $(i\in \ZZ_r^\times)$ satisfying
\[
 \fZ_{\fe_i}^* \fZ_{\fe_i}=\fZ, \qquad \text{for} \quad
 i \in \ZZ_r^\times,
\]
and $\{Z_{r_j} \, |\, j=2, \dots, m_X\}\subset \fZ_H$.

Moreover, the cyclic group of order $r$ acts on these elements;
the action of the generator $\hzeta_r \in \fC_r$ on $\fZ_{\fe_i}$, $i=0, \dots, r-1$, is defined by sending $\fZ_{\fe_i}$ to $\zeta_r^{\fe_i}\fZ_{\fe_i}$.
For $f, g \in R_H^Z$, by letting $\hzeta_r (f g) =\hzeta_r (f)\hzeta_r (g)$, $\hzeta_r(\fZ_{\fe_1}\fZ_{\fe_2}\cdots\fZ_{\fe_{r-1}})=\fZ_{\fe_1}\fZ_{\fe_2}\cdots\fZ_{\fe_{r-1}}$.
\end{Lemma}

\begin{proof}
They are obtained from Definition \ref{2df:NSG1} and Lemma \ref{2lm:NSG1}.
From the definition of $M_X$, every $Z_{r_j}$, $j=2, \dots, m_X$, belongs to $\fZ_H$. The group action is obvious from Lemma \ref{2lm:fC_r_action}.
\end{proof}

To construct our curve $X$ from $R_H$ or $\Spec R_H$, we could follow Pinkham's strategy \cite{Pinkham74} with an irreducible curve singularity with the $\GG_\fm$ action, though we will not mention it in this paper.
For the coefficients $\lambda_{ij}$ in $R_X$, we may consider the coefficient ring $\CC[\lambda_{ij}]$, and then we also consider the case $\CC[\lambda_{ij}]/\fm_{\bA}$ divided by its maximal ideal $\fm_{\bA}=(\lambda_{ij})$.
Pinkham's investigations provide the following proposition \cite{Pinkham74}:
\begin{Proposition} \label{2pr:R_H^Z}
For a given W-curve $X$ and its associated monomial ring,
$R_H^Z = \CC[Z]/ \big(f^H_1, \allowbreak f^H_2, \dots, f^H_{k_X}\big)$, there is a surjective ring-homomorphism {\rm \cite[p.~80]{ADGS2016}}
\[
\varphi^X_H\colon\ R_X \to R_H^Z
\]
such that $R_X/\fm_{\bA}$ is isomorphic to $R_H^Z$, where $\fm_{\bA}$ is the maximal ideal $(\lambda_{ij})$ in the coefficient ring $\CC[\lambda_{ij}]$, and $\varphi^X_H(y_{r_i})=Z_{r_i}$, and there is a set of polynomials $\{f^{X}_i\}_{i = 1, \dots, k_X}\in \CC[x,y_\bullet]$ satisfying
\begin{enumerate}\itemsep=0pt
\item[$(1)$] $\big(f^X_i$ modulo $\fm_{\bA}\big)=f^H_i$ for $i=1,\dots,k_X$,

\item[$(2)$]
the affine part of $R_X$ is given by $R_X=\CC[x,y_{\bullet}]/\big(f^X_1, f^X_2, \dots, f^X_{k_X}\big)$, and

\item[$(3)$] the rank of the matrix
$\big(\frac{\partial f^X_i}{\partial y_{r_j}}\big)_{
i=1, 2, \dots, k_X, j = 1, 2, \dots, m_X}$ is $m_X-1$ for every point $P$ in $X$.
\end{enumerate}
\end{Proposition}

\begin{proof}As we showed in the proof of Proposition~\ref{2pr:WCF}, $R_X$ is a quotient ring of $\CC[x,y_{\bullet}]$.
Since $\varphi_H^X$ must be a surjective, there is a prime ideal $\big(f^X_1, f^X_2, \dots, f^X_{k_X}\big)\subset \CC[x,y_\bullet]$ generated by $\big\{f^{X}_i\big\}_{i = 1, \dots, k_X}\subset \CC[x,y_\bullet]$ satisfying~(1).
Thus we prove (2).
By noting $k_X\ge m_X$, Nagata's Jacobi criterion \cite[Theorem~30.10]{Matsumura} shows (3).
\end{proof}

\begin{Definition} \label{2df:wt_RX}\quad
\begin{enumerate}\itemsep=0pt
\item
Recalling Lemma \ref{2lm:rs_iris}, we define the \emph{arithmetic local parameter} at $\infty$ by \cite{Onishi18}
\[
t = \frac{x^{i_r}}{y^{i_s}}.
\]

\item The degree at $\cQ(R_{X,\infty})$ as the order of the zero or singularity with respect to $t$ is naturally defined by
\begin{equation*}
\wt=\deg_\infty\colon \ \cQ(R_X) \to \ZZ,
\end{equation*}
which is called \emph{Sato--Weierstrass weight} \cite{Wei67}.

\item In the ring of the formal power series $\CC[[t_1, \dots, t_\ell]]$,
we define the symbols $d_{>n}(t_1, \dots, t_\ell)$ and $d_{\ge n}(t_1, \dots, t_\ell)$, which express that they belong to the ideals
\begin{gather*}
d_{>n}(t_1, \dots, t_\ell)\in \Big\{\sum a_{i_1,\dots, i_\ell} t_1^{i_1}
\cdots t_\ell^{i_\ell}\, |\,a_{i_1,\dots, i_\ell}=0
\text{ for } i_1+\cdots+i_\ell \leq n \Big\},
\\
d_{\ge n}(t_1, \dots, t_\ell)\in \Big\{\sum a_{i_1,\dots, i_\ell} t_1^{i_1} \cdots t_\ell^{i_\ell}\, |\,a_{i_1,\dots, i_\ell}=0
\text{ for }i_1+\cdots+i_\ell < n \Big\}.
\end{gather*}
\end{enumerate}
\end{Definition}

The weight of $y_{r_i}$ is given by
\[
\wt(y_{r_i}) = -r_i, \qquad i=1,2,\dots,{m_X},
\qquad
y_{r_i}=\frac{1}{t^{r_i}}(1+d_{>0}(t)).
\]

Since $R_X$ or $R_X^\circ$ is given by a quotient ring of $\CC[x,y_{r_2},y_{r_3},\dots, y_{r_{m_X}}]$ divided by the relations, $\big\{f_i^X\big\}_{i=1, \dots, k_X}$, we have the decomposition of $R_X$ as a $\CC$-vector space,
\begin{equation}
R_X = \bigoplus_{i=0} \CC \phi_i,\label{2eq:RXphi}
\end{equation}
where $\phi_i$ is a monomial in $R_X$ satisfying the inequalities $-\wt\,\phi_i < -\wt\,\phi_j$ for $i<j$, i.e., $\phi_0 = 1$, $\phi_1 = x$, $\dots$.

Further by assigning a certain weight on each coefficient $\lambda_{i,j}$ in (\ref{2eq:WCF1a}) so that (\ref{2eq:WCF1a}) is a homogeneous equation of weight $rs$, we also define another weight,
\begin{equation*}
\wt_\lambda\colon \ R_X \to \ZZ.
\end{equation*}

\begin{Definition}\label{2df:t_SR_N}
We define $S_X:=\{\phi_i\, |\,i = 0,1,2, \dots\}$ by the basis of $R_X$ as in~(\ref{2eq:RXphi}).
\end{Definition}

Then $N(i)=-\wt(\phi_i)$, for $\{N(i)\, |\,i \in \NN_0\}=H_X$.

\begin{Lemma}\label{2lm:varphiinf}
Let $t$ be the arithmetic local parameter at $\infty$ of $R_X$.
\begin{enumerate}\itemsep=0pt
\item[$1.$] By the isomorphism $\varphi_{\mathrm{inv}}\colon z \mapsto \frac{1}{t}$, $\varphi_{\mathrm{inv}}(R_H)\big({\cong} R_H^Z\big)$ is a subring of $\CC\big[\frac{1}{t}\big]$; for $g(z)\in R_H$, $g\big(\frac{1}{t}\big)\in \CC\big[\frac{1}{t}\big]$.

\item[$2.$] There is a surjection of ring $\varphi_\infty\colon R_{X}\to R_H$ $\big({\cong} R_H^Z\big)$; for $f\in R_X$, there is $g(z)\in R_H$ such that
\[
 (f)_\infty = g
\left(\frac{1}{t}\right)(1+d_{>0}(t)) \in R_{X,\infty},
\]
where $(f)_\infty$ means the germ at $\infty$ or $(f)_\infty \in R_{X, \infty}$ via $\varphi_H^X$ in Proposition~{\rm \ref{2pr:R_H^Z}}.
It induces the surjection $R_{X,\infty}\to R_H$ $\big({\cong} R_H^Z\big)$.
\end{enumerate}
\end{Lemma}
\begin{proof}By letting $g=\varphi_H^Z \circ\varphi_H^X(f)$, the existence of $g$ is obvious.
\end{proof}

\subsection[R\_P-module R\_X]{$\boldsymbol{R_\PP}$-module $\boldsymbol{R_X}$}\label{2ssc:R_Pmodule_R_X}

$R_X$ is an $R_\PP$-module, and its affine part is given by the quotient ring of $R_\PP[y_{\bullet}]$.
We recall Definition \ref{2df:NSG1} and Lemma \ref{2lm:Z_standardbasis}, and apply them to W-curves:

\begin{Proposition}\label{2pr:RP-moduleRX}
For $\fe_i \in \fE_{H_X}$, we find $\fy_{\fe_i}$ such that it is the monic monomial in $R_X$ whose weight is~$-\fe_i$, $(\fy_{\fe_0}=1)$ and satisfies
\[
R_X = R_\PP \oplus \bigoplus_{i=1}^{r-1} R_\PP \fy_{\fe_i}
=\bigoplus_{i=0}^{r-1} R_\PP \fy_{\fe_i}
=\langle \fy_{\fe_0}, \fy_{\fe_1}, \dots, \fy_{\fe_{r-1}}\rangle_{R_\PP}
\]
with the relations,
\begin{equation}
\fy_{\fe_i} \fy_{\fe_j} =\sum_{k=0}^{r-1} \fa_{i j k} \fy_{\fe_k},
\label{2eq:fyfy=fy}
\end{equation}
where $\fa_{i j k} \in R_\PP$, $\fa_{i j k}=\fa_{j i k}$, especially $\fa_{0j k}=\fa_{j0k}=\delta_{j k}$.

Further we let $\hH_{\fe}:=\big\{{-}\wt(f)\, | \, f \in \bigoplus_{\fe_i \in \fE_H^\times} R_\PP \fy_{\fe_i}\big\}$ and then $\hH_{\fe}\subset H_X$.
\end{Proposition}

\begin{proof}The generating formula is directly obtained from~(\ref{2eq:R_Xdecompose1}) noting $\fy_{\fe_0}=1$ and Lemma~\ref{2lm:Z_standardbasis}.
 $\hH_{\fe}\subset H_X$ is obvious.
\end{proof}

The set $\fY_X:=\{\fy_{0}, \fy_{\fe_1}, \dots, \fy_{\fe_{r-1}}\}$ is called the standard basis of $R_X$ as an $R_\PP$-module, which is essentially the same as $\tfy_{\Bfe_i}$ in the proof of Proposition \ref{2pr:WCF}; $\fY_X=\{\tfy_{\Bfe_i}\}\cup\{1,y_s\}$, and thus we let $\ell_s$ be $y_s=\fy_{\fe_{\ell_s}}$.
$\fa_{\ell_s,k,\ell}$ in~(\ref{2eq:fyfy=fy}) corresponds to $A_{k\ell}$ in the proof of Proposition \ref{2pr:WCF}.

The following corollary is obvious:
\begin{Corollary}\label{2cr:Ideal_fyfy}
In $\CC[x, y_\bullet]$, the ideal generated by~\eqref{2eq:fyfy=fy} is a sub-ideal of $\cI_{R_X}:=\big(f_1^X, f_2^X, \dots,\allowbreak f_{k_X}^X\big)$ in Proposition~{\rm \ref{2pr:R_H^Z}}.
\end{Corollary}

As we regard an element in $R_X$ as an element in $R_\PP$-module, we introduce a polynomial in $R_X$ as an $R_\PP$-monomial such that it is given by $\delta(x) \prod_{i=2}^{m_X} y_{r_i}^{h_i}$ for certain non-negative integers $h_i$ and $\delta(x)\in R_\PP$, e.g., $\big(x^2-a x-b\big) y_{r_2}^2 y_{r_3}$.
Since we also define the weight, $\wt_\lambda$, on the $R_\PP$-monomials, we can consider homogeneous polynomials as elements in the $R_\PP$-module.

\begin{Remark}\label{2rk:RP-moduleRX}
Corresponding to (\ref{2eq:f^H_iZ}), the relation $f^X_i$ in Proposition \ref{2pr:R_H^Z} is decomposed as
\begin{gather*}
f^X_i =
\bigg(
\prod_{j=2} y_{r_j}^{h_j^{(i+)}}
\bigg)
-\delta_{f^X_i}(x)\bigg(
\prod_{j=2} y_{r_j}^{h_j^{(i-)}}
\bigg)+ \text{lower weight terms with respect to }{-}\wt
\end{gather*}
as an $R_\PP$-module, which is relevant to~(\ref{2eq:fyfy=fy}).
Here $\delta_{f^X_i}(x)$ is an element of $R_\PP$ whose degree is $h_1^{(i+)}$ in (\ref{2eq:f^H_iZ}).
The first and the second terms are homogeneous polynomials in the Sato--Weierstrass weight $\wt$.
\end{Remark}

\subsubsection[Embedding of X into P\^{}\{2(m\_X-1)\}]{Embedding of $\boldsymbol{X}$ into $\boldsymbol{\PP^{2(m_X-1)}}$}\label{2sssc:X_to_PP^m}

The projection from $X$ to $\PP$ in Section~\ref{2sssc:X_to_PP} with the $R_\PP$-module structure induces an embedding of~$X$ into $\PP^{2(m_X-1)}$ as follows.
Besides $\big\{f^X_i\, |\,i =1, \dots, k_X\big\}$, we introduce a subset of polynomials $\big\{f_X^{(i)} \, | \, i = 2, \dots, m_X\big\}$ of $\CC[x,y_\bullet]$, which are $R_X$-analog of $f_H^{(i)}$ in Lemma~\ref{2lm:fC_r_action} such that $f_X^{(i)}=0$ an identity in~$R_X$.

\begin{Proposition}\label{2pr:WCF2}
Let $X$ be the W-curve in Proposition~{\rm \ref{2pr:WCF}} with the affine ring $R_X=$\break
$\CC[x,y_{\bullet}]/\big(f^X_1, f^X_2, \dots, f^X_{k_X}\big)$.
There are polynomials $A^{(j)}_i\in \CC[x]$, $i=2, 3, \dots, \fr_j=r/(r,r_j$, $j=2, \dots, m_X$, satisfying $ A^{(j)}_i = \sum_{k=0}^{\lfloor i \Bfr_j/\fr_j\rfloor} \lambda^{(j)}_{i, k} x^k$, where $\Bfr_j = r_j/(r, r_j)$, $\lambda^{(j)}_{i,k} \in \CC$, $\lambda^{(j)}_{r,r_j}=1$, and an irreducible polynomial,
\begin{equation}
f^{(j)}_X(x,y_{r_j}):=y_{r_j}^{\fr_j} + A^{(j)}_{1} y_{r_j}^{\fr_j-1}
+ A^{(j)}_{2} y_{r_j}^{\fr_j-2}
+\cdots + A^{(j)}_{\fr_j-1} y_{r_j}
 + A^{(j)}_{\fr_j},
\label{2eq:WCF1b3}
\end{equation}
in $\CC[x,y_{r_j}]$, especially $f_X^{(2)}=f_X$ in~\eqref{2eq:WCF1b}, so that $f^{(j)}_X$ satisfies the identity $f^{(j)}_X(x,y_{r_j})=0$ in~$R_X$.
\end{Proposition}

\begin{proof}The $m_X=2$ case is trivial.
Let us consider $m_X>2$ case. Then $\Spec(\CC[x,y]/(f_X(x,y)))$ is singular and the commutative ring $\CC[x,y]/(f_X(x,y))$ is not normal.
We normalize it to obtain~$R_X$, in which every element in $\cQ(R_X)$ is expressed by a monic equation with coefficients in $R_X$.
It means that $y_{r_j}$ satisfies a certain relation $y_{r_j}^n+b_1 y_{r_j}^{n-1}+\cdots+b_{n-1} y_{r_j} + b_n=0$ for certain positive integer $n$ and $b_i \in R_X$.
We show that when $n=r_j$, it is irreducible and $b_j \in R_\PP$ as follows.

We remark that $y_{r_j}$ is equal to an element of the standard basis $\fY_X$ from the definition of~$y_{r_j}$ in the proof of Proposition~\ref{2pr:WCF};
we take $\tj$ such that $y_{r_j}=\fy_{\fe_{\tj}}$.
We apply the investigation of~$A_{k\ell}$ for $\fy_{\fe_{\ell_s}}=y_s$ in the proof of Proposition~\ref{2pr:WCF} and $\fa_{\ell_s i k}$ in~(\ref{2eq:fyfy=fy}) to this $\fy_{\fe_{\tj}}=y_{r_j}$ case.
We introduce
\begin{equation*}
\tA^{(j)}_{i,k}:=\fa_{\tj, i, k}-\fy_{\fe_\tj}\delta_{i,k}, \qquad
\tb_i:=-\fa_{\tj,i,0}-\fa_{\tj,i,\ell_\tj}\fy_{\fe_{\tj}}, \qquad
i, k \in \ZZ_r^\times\setminus\big\{\tj\big\},
\end{equation*}
and then we consider the $\fy_{\fe_{\tj}}=y_{r_j}$ action on $\fY_X$ in~(\ref{2eq:fyfy=fy}), which is described by{\samepage
\begin{gather}
\sum_{\ell \in \ZZ_r^\times\setminus\{\tj \}} \tA^{(j)}_{k\ell}\fy_{\fe_\ell}=\tb_k,
\qquad k \in \ZZ_r^\times \setminus\big\{\tj\big\},\label{2eq:Aeta_b02}
\\
y_{r_j}^2= \BA_{\tj,0}+\BA_{\tj\tj}y_{r_j}+
\sum_{k \in \ZZ_r^\times\setminus\{\tj \}} \BA_{\tj,k} \fy_{e_k},\label{2eq:Aeta_b02a}
\end{gather}
where $\BA_{\tj,\ell}:=\fa_{\tj, \tj, \ell}\in \CC[x]$.}

If $(r_j, r)=1$, the matrix $\tA^{(j)}_{k\ell}$ is regular and we obtain~(\ref{2eq:WCF1b3}) as in the proof of Proposition~\ref{2pr:WCF}.
Thus we assume that $r_j$ and $r$ are not coprime.
We recall $\fr_j^*=(r, r_j)$, $\fr_j=r/\fr_j^*$ and $\Bfr_j=r_j/\fr_j^*$.

If the determinant of the $(r-2)\times (r-2)$ matrix $\tA^{(j)}$ is not equal to zero, it is reduced to the above case.
Then we obtain the formula $f$ which is given by $y_{r_j}^r - x^{r_j} + \cdots$ and thus its image of $\varphi_H^X$ is reduced to $Z_{r_j}^r - Z_r^{r_j}$, which is decomposed into $\big(Z_{r_j}^{\fr_j} - Z_r^{\Bfr_j}\big)\big(Z_{r_j}^{\fr_j(\fr_j^*-1)}+\cdots+ Z_r^{\Bfr_j(\fr_j^*-1)}\big)$.
If the formula $f$ is irreducible, the infinity point $\infty$ in $X$ must not be unique, which contradicts the uniqueness of $\infty$ in W-curve as mentioned in Remark~\ref{2rm:2.9}.
Thus the formula must be reduced to two formulae; one of them must be the Sato--Weierstrass weight, $-r r_j/\fr_j^*$ and contain the terms $y_{\fr_j}^r - x^{\Bfr_j}$, which is equal to $f_X^{(j)}$~(\ref{2eq:WCF1b3}).
If the other formula $h_X^{(j)}$ is algebraically independent to~$f_X^{(j)}$, there are two proejctions, $\pi_{r, r_j}\colon X\to \Spec \big(\CC[x,y_{r_j}]/\big(f_X^{(j)}\big)\big)$ and $\pi_{r, r_j}'\colon X\to \Spec \big(\CC[x,y_{r_j}]/\big(h_X^{(j)}\big)\big)$.
It also contradicts the uniqueness of the infinity~$\infty$ at~$X$ in W-curves. Hence $f=\big(f_X^{(j)}\big)^{\fr_j^*}$ up to a constant factor. Hence we have~(\ref{2eq:WCF1b3}) and the identity $f_X^{(j)}=0$.

Hence we further assume that the matrix $\tA^{(j)}$ is singular, and its rank is $q(<r-2)$.

Now we consider the case that $q=0$.
Then in the relation~(\ref{2eq:Aeta_b02a}), we show that $\BA_{\tj,k}$, $k\in \ZZ_r\setminus\big\{\tj\big\}$, must vanish.
If some of $\BA_{\tj,k}$ does not vanish, there is an element $\ty \in \bigoplus_{i \in \ZZ_r\setminus\{\tj\}}R_\PP \fy_{\fe_i}$ which is expressed by a meromorphic function of~$x$ and~$y_{r_j}$ from the relation.
However, it means that there is a non-trivial relation in~(\ref{2eq:Aeta_b02}) and thus the rank $q$ must not be zero. It contradicts the assumption.
Hence~$y_{r_j}^2$ is expressed as $\BA_{\tj,0}+\BA_{\tj\tj}y_{r_j}$; the order of $\BA_{\tj,0}$ in $x$ is $\Bfr_j$.
Accordingly, if $q = 0$, there exists a $R_\PP$-submodule $\tR'=R_\PP \oplus R_\PP y_{r_j}$, which forms a subring of~$R_X$, $\CC[x, y_{r_j}]/$ $\big(y_{r_j}^2-\BA_{\tj,0}-\BA_{\tj\tj}y_{r_j}\big)\subset R_X$.
$X$ is a covering of a hyperelliptic (or elliptic) curve, and~$\wt(x)$ is divisible by two and~$2|r$.
It is obvious that $(y_{r_j}^2-\BA_{\tj,0}-\BA_{\tj\tj}y_{r_j})$ is irreducible as an $R_\PP$-module since $X$ is not a covering of decomposed curve, and we obtain~(\ref{2eq:WCF1b3}).

Hence we let $q\neq0$ or $0<q<r-2$.
We introduce a subset $I:=\{n_1, n_2, \dots, n_q\}$ of $\ZZ_r^\times$ and a submodule $\tR_X$ of $R_\PP$-module $R_X$ defined by
\[
\tR_X=R_\PP \oplus R_\PP y_{r_j}\oplus\bigoplus_{i \in I}R_\PP \fy_{\fe_i}.
\]

We assume that $\tR_X$ is closed for the $y_{r_j}$ action on $\tR_X$, i.e.,
\begin{equation}
y_{r_j}\tR_X\subset \tR_X.\label{2eq:Aeta_b03}
\end{equation}
The assumption enables us to find the regular submatrix $\tA^{I}:=\big(\tA^{(j)}_{i,k}\big)_{i, k\in I}$ of $\tA^{(j)}$ satisfying
\begin{equation}
\sum_{k \in I} \tA^{(j)}_{i k}\fy_{\fe_k} =\tb_i, \qquad i \in I.\label{2eq:Aeta_b020}
\end{equation}
We let $\BI:=I\cup\{0,r_j\}$, and $\fY_X^\BI:=\big\{ \fy_{\fe_i}\, |\,i \in \BI\big\}$.

By considering the image of (\ref{2eq:Aeta_b020}) under $\varphi_H^Z\circ\varphi_H^X$, the weight $-\wt$ of each component in $\tA^I$ obviously leads the fact that there is a sub-monoid $H':=\langle r, r_j, \fe_{n_1}, \dots, \fe_{n_q}\rangle$ such that $\langle r, r_j\rangle \subset H' \subset H_X$, and the set $\BI$ is characterized by
\begin{equation}
\BI = \{\ell \in \ZZ_r^\times \, |\,\fe_{\ell}\in H'\}\neq \varnothing, \qquad\text{and}\qquad
 \sum_{k\in \BI} [\fe_{k}]_r = H'.\label{2eq:Aeta_b022}
\end{equation}
Then we find an expression $\fy_{e_{k}}$ of $k \in I$ as a meromorphic function of $x$ and $y_{r_j}$,
\begin{equation}
\fy_{\fe_{k}}= \frac{\tQ_k(x, y_{r_j}) } {|\tA^{I}|(x, y_{r_j})},\label{2eq:Aeta_b021}
\end{equation}
where $\big|\tA^{I}\big|$ is a monic degree $q$ polynomial of $y_{r_j}$ with $R_\PP$ coefficients, whereas the degree of $\tQ_k$ in $y_{r_j}$ is $q-1$.
(\ref{2eq:Aeta_b021})~means that $\fy_{\fe_{\ell}}$ satisfies the relation $\big|\tA^{I}\big|\fy_{\fe_{k}}=\tQ_k(x, y_{r_j})$.

Then by substituting (\ref{2eq:Aeta_b021}) into $\fy_{\fe_{k}}$ in the relation
$y_{r_j}^2= \BA_{\tj,0}+\BA_{\tj\tj}y_{r_j}+ \sum_{k\in I} \BA_{\tj,k} \fy_{e_k}$, we obtain
\[
y_{r_j}^{q+2} + \BA_{1}y_{r_j}^{q+1}+\cdots+ \BA_{q+1}y_{r_j}+\BA_{q+2}=0,
\]
where $\BA_i$ is a certain element in $R_\PP$.
We state $(q+2)|r$ because $\varphi_H^X\big(y_{r_j}^{\fr_j}-x^{\Bfr_j}\big)=0$ and it must belong to $\varphi_H^{X-1}(\{0\})$. Thus we let $(q+2) = n\fr_j$ for an integer $n(\ge 1)$. However due to $\varphi_H^X\big(y_{r_j}^{\fr_j}-x^{\Bfr_j}\big)=0$ again, $y_{r_j}^{\fr_j}-x^{\Bfr_j}$ is equal to lower weight terms with respect to~$-\wt$ because of the uniqueness of the~$\infty$ in~$X$.
Thus the equation is reduced to multiply the same equation with~$\fr_j$ order in~$y_{r_j}$, due to the above arguments.
It means that $H'\subset \langle r, r_j\rangle$, and thus $H'=\langle r, r_j\rangle$, i.e., $q=n\fr_j-2$.
Hence we obtain~(\ref{2eq:WCF1b3}) as in the proof of Proposition~\ref{2pr:WCF}, which is irreducible.

We now show (\ref{2eq:Aeta_b03}) under (\ref{2eq:Aeta_b022}).
Assume that $y_{r_j}\tR_X\setminus \tR_X \neq \varnothing$ and let $\BI^\fc:=\ZZ_r^\times\setminus\big(\BI\big)$.
Obviously $y_{r_j} \fy_{\fe_i}$ belongs to $\tR_X$ for $i \in I$ because of~(\ref{2eq:Aeta_b020}). Hence the assumption means
\[
y_{r_j}^2= \BA_{\tj,0}+\BA_{\tj\tj}y_{r_j}
+ \sum_{i\in I} \BA_{\tj,i} \fy_{\fe_{i}}
+\sum_{k\in \BI^\fc} \BA_{\tj,k} \fy_{\fe_k},
\]
where there exists, at least, a non-vanishing $\BA_{\tj,\ell} \in R_\PP$ for a certain $\ell \in \BI^\fc$, and $\BA_{\tj,i} \in R_\PP$. By considering the Sato--Weierstrass weight of the both hand sides, we have $-\wt(\BA_{\tj,\ell}\fy_{\fe_{\ell}})=2r_j \in H'$ and $\fe_{\ell}$ belongs to $H'$. It means $\ell \in \BI$ and thus contradicts (\ref{2eq:Aeta_b022}). Hence we show $y_{r_j}\tR_X\subset \tR_X$ or (\ref{2eq:Aeta_b03}).
\end{proof}

Proposition \ref{2pr:WCF2} obviously leads the following observations.
Using~(\ref{2eq:WCF1b3}) in Proposition~\ref{2pr:WCF2}, we introduce $R_{X^{(i)}}:=\CC[x, y_{r_i}]/\big(f_X^{(i)}\big)$, $i=2, \dots, m_X$, which are $R_\PP$-modules and the unnormalized rings for $m_X>2$, and their associated singular curves $X^{(i)}$, $i=2, \dots, m_X$, with the projection $\varpi_{X^{(i)}}\colon X^{(i)} \to \PP$, $\varpi_{X^{(i)}}(x,y_i)=x$, $i = 2, \dots, m_X$.
Since $f_X^{(i)}$ is irreducible, $R_{X^{(i)}}$ is a~subring of $R_X$, and $R_X$ is the normalized ring of~$R_{X^{(i)}}$. There are injective ring-homomorphisms $R_\PP\xrightarrow{\iota_r^{(i)}} R_{X^{(i)}}\xrightarrow{\iota_{r,r_i}} R_X$; thus it induces the projections $\varpi_{r_i, r}\colon X \to X^{(i)}$ $((x, y_\bullet) \mapsto (x, y_{r_i}))$ and $\varpi_{r}^{(i)}\colon X^{(i)} \to \PP$ $((x, y_{r_i}) \mapsto x)$. They satisfy the commutative diagrams,
\begin{equation*}
\xymatrix{
 R_X& R_{X^{(j)}}\ar[l]_{\iota_{r_j,r}}\\
R_{X^{(i)}} \ar[u]^{\iota_{r_i,r}}&
R_\PP,\ar[l]^{\iota_{r}^{(i)}}
\ar[lu]^{\iota_r}
\ar[u]_{\iota_{r}^{(j)}}
}\qquad
\xymatrix{
 X \ar[dr]^{\varpi_r}\ar[r]^-{\varpi_{r,r_j}}
\ar[d]_{\varpi_{r,r_i}}& X^{(j)} \ar[d]^-{\varpi_{r}^{(j)}} \\
X^{(i)} \ar[r]_{\varpi_{r}^{(i)}} & \PP.
}
\end{equation*}

Further we also define the tensor product of these rings $R_{X^{(2)}}\otimes_{R_\PP} R_{X^{(3)}}\otimes_{R_\PP} \cdots \otimes_{R_\PP}R_{X^{(m_X)}}$, and its geometrical picture $X^{[m_X-1]}_\PP:=X^{(2)}\times_\PP X^{(3)} \times_\PP\cdots \times_\PP X^{(m_X)}$.
By identifying $\CC[x,y_{\bullet}]$ $/\big(f_X^{(2)},$ $f_X^{(3)},\dots, f_X^{(m_X)}\big)=R_X^{\otimes[m_X-1]}$ with a ring $R_{X^{(2)}}\otimes_{R_\PP} R_{X^{(3)}}\otimes_{R_\PP} \cdots \otimes_{R_\PP}R_{X^{(m_X)}}$, we have the natural projection $\varphi_{R_X^{\otimes[m_X-1]}}\colon R_X^{\otimes[m_X-1]}\to R_X$, i.e.,
$R_X=R_X^{\otimes[m_X-1]}/\big(f_1^{X}, \dots, f_{k_X}^X\big)$, and the injection $\iota_{R_X^{\otimes[m_X-1]}}\colon R_\PP\hookrightarrow R_X^{\otimes[m_X-1]}$.
It induces the injection $\iota_{ X^{[m_X-1]}_\PP}\colon X \to X^{[m_X-1]}_\PP$ and
the projection $\prod_i \varpi_{X^{(i)}}\colon X^{[m_X-1]}_\PP \to \PP$.

Moreover, we also define the direct product of these rings $R_X^{[m_X-1]}:=R_{X^{(2)}}\times R_{X^{(3)}}\times \cdots \times R_{X^{(m_X)}}$, and its geometrical picture $X^{[m_X-1]}:=\prod X^{(i)}\subset \PP^{2(m_X-1)}$.
Then we have the following proposition.

\begin{Proposition}\label{2pr:X_to_PP^m-1}
Expression \eqref{2eq:WCF1b3} in Proposition~{\rm \ref{2pr:WCF2}} provides surjective and injective ring homomorphisms,
\[
\varphi_{R_X^{\otimes[m_X-1]}}\colon \
R_X^{\otimes[m_X-1]}
\to R_X,\qquad
\iota_{R_X^{\otimes[m_X-1]}}\colon \
R_\PP\hookrightarrow R_X^{\otimes[m_X-1]},
\]
so that they satisfy the commutative diagrams
\begin{equation*}
\xymatrix{
R_X& R_X^{\otimes[m_X-1]} \ar[l]_{\varphi_{R_X^{\otimes[m_X-1]}}}, \\
R_\PP\ar[u]\ar[ur]_{\iota_{R_X^{\otimes[m_X-1]}}} &
}
\qquad
\xymatrix{
X \ar[d] \ar[r]^{\iota_{ X^{[m_X-1]}_\PP}}& X^{[m_X-1]}_\PP \ar[dl]^-{\prod\varpi_{X^{(i)}}}, \\
\PP & \\
}
\end{equation*}
which is consistent with the projection
$\varpi_{ X^{[m_X-1]}_\PP}\colon X^{[m_X-1]}_\PP\to X$.
They induce the homomorphism by the ideal $(x-x_2, x-x_3, \dots, x-x_{m_X})$
and an embedding,
\[
\varphi_{R_X^{[m_X-1]}}\colon \ R_X^{[m_X-1]}
\to R_X, \qquad \iota_{ X^{[m_X-1]}}\colon \ X \hookrightarrow X^{[m_X-1]} (\subset \PP^{2(m_X-1)}).
\]
\end{Proposition}

\begin{Remark}
Since the normalization of a ring is not unique in general, the surjective ring homomorphism $\varphi_{R_X^{\otimes[m_X-1]}}$ is not injective except $m_X=2$.
For example, in the $(3,4,5)$ case, there is an surjective ring homomorphisms whose $\varphi_H$ image is given by
\begin{gather*}
\CC[Z_3,Z_4,Z_5]/\big(Z_4^3-Z_3^4,Z_5^2-Z_3^5\big)\\
\qquad{} \to \CC[Z_3,Z_3,Z_4]/\big(Z_4^2- \zeta_3^a Z_3 Z_5, Z_4 Z_5 - Z_3^3, Z_5^2 - \zeta_3^{2a}Z_3^2 Z_4\big).
\end{gather*}
Thus $y$ and $y_{r_j}$ are related via $f^X_i$'s whose
image of $\varphi_H^X$ are binomial relations.
\end{Remark}

Direct computation gives the following relation based on~(\ref{2eq:WCF1b3}).

\begin{Lemma}\label{2lm:dy_rj/dx}
\begin{align*}
\big(\delta_{Y,y} f_X^{(j)}\big)(x,Y,y)
:={}& \frac{f_X^{(j)}(x,Y) - f_X^{(j)}(x,y)}{Y-y} \\
={}& \sum_{\ell=0}^{\Bfr_j-1}A_{\ell}^{(j)}(x)
\sum_{i=0}^{\fr_j-\ell-1}
Y^i y^{\fr_j-\ell-i-1} \in R_X \otimes_{R_\PP} R_X.
\end{align*}
\end{Lemma}

\subsection{The covering structures in W-curves}\label{2sc:Covering}

We follow \cite{Kunz2004, Stichtenoth} to investigate the covering structure in W-curves.

\subsubsection{Covering structure}
\label{2sc:GaloisCover}

As mentioned in Remark \ref{2rm:g_Xc_X}, let us consider the Riemann sphere $\PP$ and $R_\PP=\bH^0(\PP, \cO_{\PP}(*\infty))$. We identify $R_\PP$ with its affine part $R_\PP^\circ=\CC[x]$ and its quotient field is denoted by $\CC(x)=\cQ(R_\PP)$.
The quotient field $\cQ(R_X)=\CC(X)$ of $R_X$ is an extension of the field $\CC(x)$.

Following the above description, we consider the covering structure of the W-curve $X$.
The covering $\varpi_r\colon X \to \PP$ $((x,y_\bullet) \mapsto x)$ is obviously the holomorphic $r$-sheeted covering.
Further, when we have the Galois group $\Gal(\cQ(R_X)/\cQ(R_\PP))=\Aut(X/\PP)=\Aut(\varpi_r)$, it is denoted by $G_X$.
The $\varpi_{r}$ is a finite branched covering.
Each point in $\varpi_r^{-1}(x)$ for $x \in \PP$ except certain finite points is biholomorphic.
A ramification point of $\varpi_{r}$ is defined as a point of such that is not biholomorphic at the point.
The image $\varpi_{r}$ of the ramification point is called the branch point of $\varpi_{r}$. The number of the finite ramification points is denoted by $\ell_{\fB}$.

We basically focus on the holomorphic $r$-sheeted covering $\varpi_{r}\colon X \to \PP$. $G_x$ denotes the finite group action on $\varpi_r^{-1}(x)$
for $x\in \mathbb P$, refered to as group action at $x$ in this paper:
\begin{Definition}\label{2df:fBX}
Let $\fB_X:=\fB_{X,r}$ $=\{B_{i}\}_{i=0, \dots, \ell_{\fB}}$ and $\fB_{\PP}:=\varpi_r(\fB_X)=\{b_{i}\}_{i=0, \dots, \ell_{\fB}}$, where $\ell_{\fB}:= \# \fB_X-1$, $B_{0}=\infty\in X$ and $b_i :=\varpi_r(B_i)$.
\end{Definition}

We have the following results.

\begin{Lemma}\label{2lm:B_iVanish}
Every element $B_i$ in $\fB_X \setminus \{\infty\}$ is given by the point $(x, y_\bullet)$ at which there exists, at least, a certain $j$ in $\{2\le j\le m_X\}$ such that $\frac{\partial f_X^{(j)}(x,y_{r_j})}{\partial y_{r_j}}=0$
and $\frac{\partial f_X^{(j)}(x,y_{r_j})}{\partial x}\neq 0$.
\end{Lemma}

\begin{proof}
When $\frac{\partial f_X^{(i)}(x,y_{r_i})}{\partial x}=0$ and $\frac{\partial f_X^{(i)}(x,y_{r_i})}{\partial y_{r_j}}=0$ at a point $P\in X$, the point means the singular point as the plane curve given by $f_X^{(i)}(x,y_{r_i})=0$ at $\varpi_{r, r_i}(P)$.
However, since~$X$ is not singular, there exists~$j$ satisfying the condition.
Then at the point, ${\rm d} x$ is identically zero, and thus at the point, $x$ is not a local parameter of~$X$. Thus~$P$ must be an element in $\fB_X\setminus \{\infty\}$, $P=B_{i}$.
It means ${\rm d} x={\rm d}(t^{e_i})(1+d_{>0}(t))$ for a positive integer $e_i>1$ in terms of the local parameter $t$ such that $t(B_i) = 0$, and there exists $j$ such that ${\rm d} y_{r_j} = {\rm d} t$ and ${\rm d} y_{r_i} =t^{f_i} {\rm d} t(1+d_{>0}(t))$, $i \neq j$, $f_i \ge 1$.
\end{proof}

The $e_i$ appearing the proof for the ramification point $B_i$ in Lemma~\ref{2lm:B_iVanish} is called \emph{the ramification index}, and denoted by~$e_{B_i}$.

\subsubsection{Riemann--Hurwitz theorem}

Let us consider the behaviors of the covering $\varpi_r\colon X \to \PP$, including the ramification points.
The Riemann--Hurwitz theorem \cite{Kunz2004},
\begin{equation*}
2g-2 =-2r+ \sum_{i=1}^{\ell_{\fB}} (e_{B_i}-1)+(r-1)=
\sum_{i=1}^{\ell_{\fB}} (e_{B_i}-1)-(r+1),
\end{equation*}
 shows the following:
\begin{Corollary}\label{2cr:(dx)} The divisor of ${\rm d} x$ is given by
\[
\Div({\rm d} x) =\sum_{i=1}^{\ell_{\fB}} (e_{B_i}-1) B_i - (r+1)\infty.
\]
\end{Corollary}

\section[Complementary module R\_X\^{}c of R\_X]{Complementary module $\boldsymbol{R_X^\fc}$ of $\boldsymbol{R_X}$}\label{sec:CompM}

\subsection{Trace in the covering structure}\label{2ssc:TraceOperator}

It is known that the Riemann--Hurwitz relation and the divisor of ${\rm d} x$ are obtained via the Dedekind different \cite[Theorem~15.11]{Kunz2004}.
In this subsection, we consider the properties of $R_\PP$-module $R_X$, which are related to the trace, the complementary module, and the Dedekind different \cite[Chapter~15]{Kunz2004}.

\subsubsection[Trace in R\_X/R\_P]{Trace in $\boldsymbol{R_X/R_\PP}$}

We review the general results of $R_X$ as a ring extension of $R_\PP$ following \cite{Kunz2004} (see Section~\ref{2sc:App_trace}).
Let us consider the covering structures of $\varpi_r\colon X\to \PP$ to discriminate its lifted points, and the enveloping field $\cQ(R_X)^e=\CC(X)^e:=\CC(X)\otimes_{\CC(x)}\CC(X)=\cQ(R_X) \otimes_{\cQ(R_\PP)} \cQ(R_X)$.

The field extension $\CC(X)/\CC(x)$ induces the extension ring $R_X$ of $R_\PP$.
As mentioned in Section~\ref{2sc:App_trace}, we consider the dual of $R_X$ with respect to $R_\PP$,
\[
\omega_{R_X/R_\PP}:=\Hom_{R_\PP}(R_X, R_\PP),
\]
which is a free $R_\PP$-module.
For $R_X= \oplus_{i=0}^{r-1} R_\PP \bby_i$ (i.e., $\bby_i=\fy_{\fe_i}$), it is expressed by the dual basis $\{\bby_i^*\}_{i=0, \dots, r-1}$ such that
\begin{equation*}
\omega_{R_X/R_\PP}=\bigoplus_{i=0}^{r-1} R_\PP {\mathbf{y}_i^*},
\qquad\text{and}\qquad
\bby_i^* \bby_j = \delta_{i,j}.
\end{equation*}
Here $\bby_i^*$ and $\bby_j$ correspond to the point in fiber which lies over $x\in R_\PP$.
For the standard trace $\tau_{R_X/R_\PP}:=\sum_{i=0}^{r-1} \bby_i \circ \bby_i^*\in \omega_{R_X/R_\PP}$, we define the complementary module $R_{X}^{\fc}$ over $R_X/R_\PP$ with respect to the standard trace $\tau_{R_X/R_\PP}$ by
\begin{equation}
R_{X}^{\fc}:=\{z \in \cQ(R_X)\, | \, \tau_{R_X/R_\PP}(z a) \in R_\PP, \, \forall a \in R_X\}.\label{2eq:RXc_tau}
\end{equation}

It is obviously that the localization $(R_{X}^{\fc})_P$ at $P\in X$ is equal to $R_{X,P}^{\fc}$ from the definition.
Since for a point $P\in X$, $R_{X,P}$ is a~principal ideal domain, every ideal is generated by a~certain element in~$R_{X,P}$.
There is an element $\fh_{X,P} \in R_{X,P}$ such that $R_{X,P}^{\fc}=\fh_{X,P}R_{X,P}$ \cite[Proposition~3.4.2]{Stichtenoth}.
Following \cite[Definition~3.4.3]{Stichtenoth}, we define the different of~$R_X/R_\PP$.

\begin{Definition}\label{2df:2.32}
The different $\diff(R_X/R_\PP)$ is a divisor defined by
\[
\diff(R_X/R_\PP):=\sum_{P\in X\setminus\{\infty\}} {\rm d}_{P} P,
\]
where the different divisor ${\rm d}_{P}:= -\deg_P(\fh_{X,P})$ for $R_{X,P}^\fc =\fh_{X,P} R_{X,P}$.
\end{Definition}

By Dedekind's different theorem, we have the following:
\begin{Proposition}\label{2pr:DedekindDiff1}
${\rm d}_{B_i}=e_{B_i}-1$ for $P \in \fB_X\setminus\{\infty\}$, and the support of $\diff(R_X/R_\PP)$ equals $\fB_X\setminus \{\infty\}$.
\end{Proposition}

\begin{proof}See \cite[Theorem 15.11]{Kunz2004} and \cite[Theorem 3.5.1]{Stichtenoth}.
\end{proof}

\subsection[Trace operator for plane W-curves: R\_X=R\_\{X\^{}circ\} (m\_X=2) case]{Trace operator for plane W-curves: $\boldsymbol{R_X = R_{X^\circ}}$ ($\boldsymbol{m_X=2}$) case}\label{2sssec:m_X=2}

Following Kunz \cite[Theorem 15.1]{Kunz2004}, we review $R_X^\fc$ of the $m_X=2$ case:
\begin{Proposition}[{\cite[Theorem 15.1]{Kunz2004}}]\label{2pr:Kunz_m2}
For the plane W-curve $(m_X=2)$, we have the relation,%
\begin{equation}
 R_X^\fc = \frac{1}{f_{X,y}}\cdot R_X.\label{2eq:RXfc}
\end{equation}
\end{Proposition}

\begin{proof}Let us show the proof by Kunz.
We note that the extension of field $\cQ(R_X)$ of $\cQ(R_\PP)$ is separable and $f_X$ is monic,
\[
\cQ(R_X)=\cQ(R_\PP)[Y]/(f_X) = \bigoplus_{i=0}^{r-1} \cQ(R_\PP)y^i,
\]
and thus we note $y^i =\bby_i = \fy_{\fe_i}$ in Proposition \ref{2pr:RP-moduleRX} and in Section~\ref{2sc:App_trace}, and
\[
\cQ(R_X)^e = \cQ(R_X)\otimes_{\cQ(R_\PP)} \cQ(R_X) \cong \cQ(R_X)[Y]/(f_X).
\]
Using the ring-homomorphism $\mu\colon \cQ(R_X)^e\to \cQ(R_X)$ $(a\otimes b \mapsto ab)$, the standard trace $\tau_{R_X/R_\PP}:=\sum_{i=0}^{r-1} y^i \circ \bby_i^*$ is obtained by the extension of a trace $\tau$.
Let us find the basis $\{\hby_i\}$ of $R_X$ with respect to $\tau$, and $\Delta_\tau:=\bigoplus_{i=0}^{r-1} \hby_i\otimes y^i$ as an element of $R_\PP$-module $\Ann_{R_X^e}(\Ker \mu)$. If we find~$\{\hby_i\}$ and~$\tau$, using them, we obtain the standard trace $\tau_{R_X/R_\PP}=\mu(\Delta_\tau)\circ \tau$.
Accordingly, we construct the $\tau$ and $\{\hby_i\}$ as follows.

For an element $\th \in \cQ(R_X)[Y]$ such that $f_X(x,Y)=(Y-y) \cdot \th$,
we can identify $\Ann_{R_X^e}(\Ker \mu)$ with the principal ideal $(\th)/(f_X(x,Y))$.
Indeed, in $\cQ(R_X)[Y]/(f_X(x,Y))$,
\[
f_X(x,Y) = f_X(x,Y)-f_X(x,y) = (Y-y)\frac{f_X(x,Y)-f_X(x,y)}{Y-y}=(Y-y)\th=0.
\]
It means that $\th(x,Y,y):=\delta_{Y,y} f_X(x,Y,y)$ belongs to $R_X\otimes_{R_\PP} R_X$ noting Lemma~\ref{2lm:dy_rj/dx}.
Thus we have
\[
\Delta_{\th}=\sum_{\ell=0}^{r-1}A_\ell(x)\sum_{j=0}^{r-\ell-1}
y^{j}\otimes y^{r-\ell-j-1} \in R_X^e=R_X \otimes_{R_\PP} R_X,
\]
which generates the ideal $\Ann_{R_X^e}(\Ker \mu)$ in $R_X^e$, i.e.,
$(a\otimes 1-1\otimes a) \Delta_{\th}=0$ for every $a\in R_X$.
Further from Proposition~\ref{2pr:A_traceH20}, it corresponds to the trace $\tau_\th\in \omega_{R_X/R_\PP}$ by
\[
1=\sum_{\ell=0}^{r-1}A_\ell(x)
\sum_{j=0}^{r-\ell-1}\tau_\th\big(y^{j}\big)y^{r-\ell-j-1}.
\]
For example, $r=4$ case, it is
\[
1=\tau_\th\big(y^3 + A_1 y^2 + A_2 y+ A_3\big) +
\tau_\th\big(y^2 + A_1 y + A_2\big)y +
\tau_\th(y + A_1)y^2+\tau_\th(1)y^3,
\]
which should be interpreted as $\tau_\th(1)=0$, $\tau_\th(y)=0$, $\tau_\th\big(y^2\big)=0$, and thus $\tau_\th\big(y^3\big)=1$.
Similarly since $A_0=1$ and $y^{r-\ell-j-1}=1$ when $\ell=0$ and $j=r-1$, we compare the both sides in the equation and obtain
\[
\tau_\th(y^i) =
\begin{cases}
1 & \text{for } i=r-1,\\
0 & \text{otherwise}.
\end{cases}
\]
In other words, for $\bby_i = y^i$, $i\in \ZZ_r$, we have
\[
\hby_i = \sum_{\ell=0}^{r-1-i}A_\ell(x)y^{\ell}
\]
as the dual basis of the basis $\{y^i =\bby_i = \fy_{\fe_i}\}_{i\in \ZZ_r}$ with respect to $\tau_\th$ in Section~\ref{2sc:App_trace}.
Then $\tau_\th(R_X) \subset R_\PP$. It is obvious that $\mu(\Delta_{\th})= f_{X,y}$.
Thus the standard trace $\tau_{R_X/R_\PP}$ of $R_X/R_\PP$ is given by
\[
\tau_{R_X/R_\PP}=\mu(\Delta_\th)\circ \tau_\th=
 f_{X,y}\circ \tau_\th.
\]
We have
\[
\tau_{R_X/R_\PP}
\left(\frac{y^i}{ f_{X,y}(x,y)}\right)
= \begin{cases}
1 & \text{for } i=r-1,\\
0 & \text{otherwise}.
\end{cases}
\]
Let us consider an element in $R_{X}^{\fc}$ as in~(\ref{2eq:RXc_tau}).
We consider $u \in \cQ(R_X)$ which is expressed by
\[
u = \frac{1}{f_{X,y}(x,y)} \sum_{i=0}^{r-1} a_i y^i \in \cQ(R_X),
\]
and its conditions in (\ref{2eq:RXc_tau}).
For every element $v = \sum_{j=0}^{r-1} b_j y^j \in R_X$, $b_j \in R_\PP$, $u$ satisfies $\tau_{R_X/R_\PP}(u v)$ $ \in R_\PP$, i.e.,
\[
\tau_{R_X/R_\PP}(u v)= \sum_{i,j} \tau_\th\big(a_i b_j y^i y^j\big)= \sum_{i=0}^{r-1} a_i b_{r-1-i}\in R_\PP.
\]
It implies that $a_i \in R_\PP$ for every $i=0, 1, \dots, r-1$.
Thus we obtain the relation (\ref{2eq:RXfc}).
\end{proof}

\begin{Remark}\label{2rmk:hetam=2}
Here we note that instead of $\{\hby_i\}$ in the proof, we define a simpler dual basis of $\{y^i =\bby_i = \fy_{\fe_i}\}_{i\in \ZZ_r}$ with respect to $\tau_\th$ by
\[
\big\{\hfy_{\fe_i}:=y^{r-i-1}\big\}_{i \in \ZZ_r},
\]
because it is obvious that
$\tau_\th(\hfy_{\fe_i}\fy_{\fe_j})=\delta_{ij}$ from the definition of $\tau_\th$.
\end{Remark}

From (\ref{2eq:RXfc}), the Dedekind different of the plane W-curve is given by
\[
\diff(R_X/R_\PP)=-\sum_{P\in X\setminus\{\infty\}}\deg_P(f_{X,y}) P.
\]
By Dedekind's different theorem (Proposition \ref{2pr:DedekindDiff1}), we have the following \cite[Chapter 15]{Kunz2004}:
\begin{Proposition}\label{2pr:DedekindDiff2}
$\deg_{B_i}(f_{X,y})=e_{B_i}-1$.
\end{Proposition}

\subsubsection[Trace in the holomorphic r-sheeted covering at R\_X]{Trace in the holomorphic $\boldsymbol{r}$-sheeted covering at $\boldsymbol{R_X}$}\label{2sec:Trace_RX1}

We generalize the result of the plane curves to the space curves.
The surjective ring homomorphism in Proposition~\ref{2pr:X_to_PP^m-1} can be interpreted as follows.
We implicitly introduce the dual modules $\omega_{R_{X^{(i)}}/R_\PP}$ of $R_{X^{(i)}}$ and their tensor product $\omega_{R_X/R_\PP}^{\otimes[m_X-1]}:=\omega_{R_{X^{(2)}}/R_\PP}\otimes_{R_\PP}\cdots\otimes_{R_\PP}\omega_{R_{X^{(m_X)}}/R_\PP}$ to find~$R_X^\fc$.
More precisely, we implicitly construct the trace $\tau$ in $R_X$ by using the data of $(\tau_2, \dots, \tau_{m_X})$ of $R_X^{\otimes[m_X-1]}$ by regarding~$R_X$ as $R_X=R_X^{\otimes[m_X-1]}/\big(f^X_1, f^X_2, \dots, f^X_{k_X}\big)$.
By investigating them, we obtain $R_X^\fc$ for the $m_X>2$ case.

We first investigate the similar relations in $R_{H}$ following Lemma~\ref{2lm:fC_r_action}.

\subsubsection[Trace in the Galois covering at R\_H]{Trace in the Galois covering at $\boldsymbol{R_{H}}$}\label{2sec:Trace_RH}

In order to generalize the investigation of the complementary module for the plane curves $m_X=2$ to general W-curves, we investigate the trace structure at $R_{H}$ since the monomial curve in Section~\ref{2ssc:MCurve} is crucial in W-curves.
We assume $m_X\ge 2$.

We use the surjection $\varphi_\infty\colon R_X \to R_H$ in Lemma~\ref{2lm:varphiinf}, and consider the behavior of the trace in $R_H^Z$.
We investigate a {\lq\lq}covering{\rq\rq} structure in $\varpi_H\colon X_H=\Spec R_H^Z \to \PP=\Spec \CC[Z_r]$.
The cyclic group $\fC_r$ of order $r$ acts on $R_H^Z$ and $X_H$ as the $\GG_\fm$ action.
We regard it as the Galois covering and consider $\CC[Z_r]$-module $R_H^Z$ and $R_H$.

We introduce a meromorphic function on $\Spec R_H^Z \times_{\Spec \CC[Z_r]} \Spec R_H^Z$ or an element $p$ in the enveloping field $\cQ(R_H^Z)^e:=\cQ\big(R_H^Z \otimes_{\CC[Z_r]} R_H^Z\big)$, and an element $h$ in its associated enveloping ring ${R_H^{Z}}^e:=R_H^Z \otimes_{\CC[Z_r]} R_H^Z$, i.e., $Z_{r_i}^{r}=Z_r^{r_i}=Z_{r_i}^{\prime r}$ for $i=2, \dots, m_X$.
We extend the group action of $\fC_r$ to that on $\cQ\big(R_H^Z \otimes_{\CC[Z_r]} R_H^Z\big)$ such as $\zeta h(Z_r, Z_\bullet, Z'_\bullet) = h(Z_r, \zeta Z_\bullet, \zeta Z'_\bullet)$ for $\zeta \in \fC_r$.
We define an element $p_{H}$ in $\cQ\big(R_H^Z \otimes_{\CC[Z_r]} R_H^Z\big)$ by
\begin{gather}
p_{H}(Z_r, Z_\bullet, Z'_\bullet) :=\prod_{i=2}^{m_X}p_{H, r_i}(Z_r, Z_{r_i}, Z'_{r_i}),\nonumber\\
p_{H, r_i}(Z_r, Z_{r_i}, Z'_{r_i}) :=\frac{Z^{\fr_i-1}_{r_i}+Z^{\fr_i-2}_{r_i}Z^{\prime}_{r_i}
+Z^{\fr_i-3}_{r_i}Z^{\prime 2}_{r_i}+\cdots+Z^{\prime \fr_i-1}_{r_i}}{r Z_{r_i}^{\fr_i-1}}.\label{2eq:p_varphiH}
\end{gather}
Due to $Z_{r_i}^{\fr_i} = Z_r^{\fr_i}=Z_{r_i}^{\prime \Bfr_i}$ in Lemma~\ref{2lm:fC_r_action}, then each factor behaves like
\[
p_{H, r_i}(Z_r, Z_{r_i}, Z'_{r_i})
= \begin{cases}
1 & \text{for } Z^{\prime}_{r_i}=Z_{r_i},\\
0 & \text{otherwise}.\end{cases}
\]
We have the trace property,
\begin{equation}
p_{H}(Z_r, Z_\bullet, Z'_\bullet)
= \begin{cases}
1 & \text{for } Z^{\prime}_{\bullet}=Z_{\bullet},\\
0 & \text{otherwise}.\end{cases}
\label{2eq:6aa}
\end{equation}

\begin{Lemma}\label{2lm:h_H}
There are polynomials $\th_{H}(Z_r, Z_\bullet, Z'_\bullet)\in R_H^Z \otimes_{\CC[Z_r]} R_H^Z$ and $h_H(Z)=$ $h_H(Z_r, Z_\bullet)$ $:= \th_{H}(Z_r, Z_\bullet, Z_\bullet)\in R_H^Z$ such that
\begin{equation}
\varphi_H^Z\big(\th_{H}(Z_r, Z_\bullet, Z'_\bullet)\big)=z^{d_h}\sum_{i=0}^{r-1}
z^{\prime\,\fe_i }z^{-\fe_i}=\sum_{i=0}^{r-1}z^{\prime\,\fe_i }z^{\hfe_i}
\in R_H \otimes_{\CC[z^r]} R_H,\label{2eq:Sigma_X^H2}
\end{equation}
where
\begin{enumerate}\itemsep=0pt
\item[$(1)$] $z$ and $z'$ are given by $\varphi_H^Z(Z_{r_i})=z^{r_i}$, $\varphi_H^Z(Z_{r_i}')=z^{\prime\, r_i}$, for $i \in \ZZ_r$, and $z^r = z^{\prime\, r}$,
\item[$(2)$] $d_h$, $\hfe_i$, $\delta_i$, $i\in \ZZ_r$, and $\ell$ are non-negative integers satisfying the conditions that
\begin{enumerate}\itemsep=0pt
\item[$(a)$] $d_h \ge \fe_i$ for every $i\in\ZZ_r$,
\item[$(b)$] $d_h$ is determined by Lemma~{\rm \ref{2lm:h_RXe}},
\item[$(c)$] $d_h = \fe_\ell+ \delta_0 r$ such that $\fe_\ell = \sum_{i=2}^{m_X} (\fr_i-1)r_i$ modulo $r$, and
\item[$(d)$]
\begin{equation}
\hfe_i:=d_h - \fe_i=\fe_{\ell,i}^*+\delta_i r,\qquad
z^{\fe_i}z^{\hfe_i}=z^{\fe_0}z^{\hfe_0}=z^{d_h},\label{2eq:hfe}
\end{equation}
for every $i\in\ZZ_r$, especially
\begin{equation}
\hfe_0 = d_h = \fe_\ell + \delta_0 r, \qquad
\hfe_\ell = \delta_0 r, \qquad \delta_\ell =\delta_0,
\label{2eq:hfe0}
\end{equation}
\end{enumerate}
\item[$(3)$] $\th_{H}(Z_r, Z_\bullet, Z'_\bullet)$ consists of $r$ monomials corresponding to each term in \eqref{2eq:Sigma_X^H2} and satisfies
\begin{equation}
p_{H}(Z_r, Z_\bullet, Z'_\bullet)=\frac{\th_{H}(Z_r, Z_\bullet, Z'_\bullet)}{h_H(Z)}.\label{2eq:hfe2}
\end{equation}
\end{enumerate}
\end{Lemma}

\begin{proof}The cyclic group $\fC_r$ acts on $p_{H}(Z_r, Z_\bullet, Z'_\bullet)$ so that it is invariant.
We consider $R_H$ rather than $R_H^Z$. We note the following:
1) For $j=2, \dots, m_X$, $r_i$ and $r_j\in M_X$, $\{r_j i \text{ modulo } r\, |\,i\in \ZZ_r\}=\ZZ_r$,
2) the numerator in each $p_{H, r_i}$ in (\ref{2eq:p_varphiH}) is homogeneous, and
3) their product is also homogeneous.
Therefore we see that there are non-negative integers, $\td_\Delta$ and $\Delta_i$, $\Delta_{i}<\Delta_{i+1}$, such that
\[
\frac{z^{\td_\Delta}\sum_{i=0}^{r-1}
z^{\prime\, \Delta_i}z^{-\Delta_i}
}
{r z^{\td_\Delta}}
\]
has the property of the right-hand side of (\ref{2eq:6aa}) after acting $\varphi_H^Z$ both sides in~(\ref{2eq:6aa}), and $\{\Delta_i \text{ modulo } r\, |\,i\in \ZZ_r\}=\ZZ_r$.
It means that $\{\Delta_i \, |\,i\in \ZZ_r\}$ equals $\{\fe_i + n_i r \, |\,i\in \ZZ_r\}$ for a~certain non-negative number $n_i \in \NN_0$ from Lemma~\ref{2lm:NSG1}, and $\Delta_0 = 0$.

Due to the isomorphism $\varphi_H^Z$, for sufficiently large $n$ and $m$, we find an element in $R_H^Z \otimes_{\CC[Z_r]} R_H^Z$ whose image of $\varphi_H^Z$ is $z^{n}z^{\prime\, m}$.
For $\ell_i$ satisfying $\sum_i r_i \ell_i \equiv 0$ modulo $r$, $\prod Z_{r_i}^{\ell_i} = \prod Z_{r_i}^{\prime \ell_i}$, and thus we can find $d_h$ such that
$z^{\td_\Delta}\sum_{i=0}^{r-1} z^{\prime\, \Delta_i}z^{-\Delta_i}
=z^{d_h}\sum_{i=0}^{r-1}z^{\prime\, \fe_i}z^{-\fe_i}$ noting $z^r=z^{\prime\, r}$.
There is a preimage $\th_{H}(Z_r, Z_\bullet, Z'_\bullet)$ as an element in $R_H^Z \otimes_{\CC[Z_r]} R_H^Z$.
It is obvious that $\th_{H}(Z_r, Z_\bullet, Z'_\bullet)$ consists of $r$ monomials and satisfies~(\ref{2eq:Sigma_X^H2}) and~(\ref{2eq:hfe2}) for $p_{H}(Z_r, Z_\bullet, Z'_\bullet)$.

We note that the determination of $d_h$ has the ambiguity up to $r$, and thus we set it such that it satisfies Lemma \ref{2lm:h_RXe}.
Here $d_h = \td_h + n_h r$ so that $\td_h$ is the minimal element satisfying the above relations.

From the definition and Lemma~\ref{2lm:NSG1}\,(3), we obtain the relations (\ref{2eq:hfe}) of $i=0$ and $i=\ell$, and $\delta_0 = \delta_\ell$.

We consider the cases $d_h$ modulo $r$.
\begin{enumerate}\itemsep=0pt
\item $d_h = \fe_0 = 0$ modulo $r$ case:
We note that $\{\fe_i $ modulo $r \, |\,i\in \ZZ_r^\times\}= \ZZ_r^\times$, and $\fe_\ell < \fe_{r-1}$, whereas $\hfe_\ell$ must be non-negative.
Each $\fe_{i}$, $i \in \ZZ_r^\times$, cannot be divided by $r$ and thus $\fe_{i}$, especially, $\fe_{r-1}$ is not equal to $d_h$.
From Definition \ref{2df:NSG1}. we find a non-negative integer $\delta_i$ such that $\hfe_i =\fe_{0,i}^* + \delta_i r$, $i \in \ZZ_r^\times$. It is obvious that $d_h=\hfe_0=\delta_0r >\fe_{r-1}>0$.
We also note $\hfe_{r-1}>0$ because of $\fe_{0,r-1}^*>0$.

\item $d_h = \fe_\ell$ modulo $r$ case ($\ell = 1, \dots, r-2$):
Similarly, since $d_h = \fe_\ell + \delta_0 r$ satisfies $d_h - \fe_i\ge0$ for $i \in \ZZ_r$ (especially $d_h \ge \fe_{r-1}$), we find non-negative integers $\delta_i$ such that $\hfe_i =\fe_{\ell,i}^* + \delta_i r$ for $i \in \ZZ_r$, and $\hfe_\ell = \delta_\ell r=\delta_0 r>0$, noting $\fe_{\ell,\ell}^*=0$.
Then $\hfe_{r-1}>0$ because of $\fe_{\ell,r-1}^*>0$.

\item $d_h = \fe_{r-1}$ modulo $r$ case:
Similarly since $d_h = \fe_{r-1} + \delta_0 r$ satisfies $d_h - \fe_i\ge0$, especially $d_h \ge \fe_{r-1}$, $\delta_0$ is non-negative.
We find non-negative integers $\delta_i$ such that $\hfe_i =\fe_{r-1,i}^* + \delta_i r$ for $i \in \ZZ_r^\times$, and $\hfe_{r-1} = \delta_{r-1} r\ge 0$ or $\delta_{r-1}=\delta_0$ because of $\fe_{r-1,r-1}^*=0$.
\end{enumerate}
These show the statements in the proposition.
\end{proof}

\begin{Proposition}\label{2pr:dh_symmetric}
The $\delta_0$ in \eqref{2eq:hfe0} equals zero if and only if $d_h = \fe_{r-1}$.

The case $d_h = \fe_{r-1}$ or $\delta_0=0$, occurs if and only if~$H$ is symmetric whereas
the case $\delta_0\neq 0$ if and only if $H_X$ is not symmetric.

Thus we say that if $\delta_0=0$, $d_h$ is symmetric and otherwise, $d_h$~is not symmetric.
\end{Proposition}

\begin{proof}They are proved in Lemma~\ref{2lm:nuIo1}.
\end{proof}

\begin{Remark}\label{2rm:R_Hparameters}
We remark that $R_H$ and $R_H^Z$ are characterized by these parameters $(M_X=\{r_i\},\allowbreak m_X, k_X, \{\fe_i\},\{\hfe_i\},d_h, \{\delta_i\}, \ell)$.
Especially $\ell$ is a fixed number for a given $X$ in this paper.
\end{Remark}

\begin{Example}\quad
\begin{enumerate}\itemsep=0pt

\item $H=\langle 4, 6, 7, 9\rangle$ (non-symmetric) case $(d_h=4\delta_0+9,\, \delta_0=1,\, \ell=3)$:
$H = \{0, 4, 6, 7, 8, 9, 10, \allowbreak 11, \dots\}$, $H^\fc = \{1, 2, 3, 5\}$,
$R_H = \CC[Z_4, Z_6, Z_{7}, Z_{9}]/{\sim} $ and then
\begin{gather*}
h_{{R_Z^H}}(Z, Z_\bullet, Z_\bullet') =
Z_4Z_9 + Z_{7} Z_{6}'+Z_6 Z_{7}'+ Z_4 Z_9',
\qquad
h_H(Z) = 4 Z_4Z_9.
\end{gather*}
\begin{table}[htb]\centering\vspace{-5mm}
 \begin{tabular}{|c|cccc|}
\hline
$i$&0&1&2&3\\
\hline
$\fe_i $ &0 &6 & 7& 9 \\
$\hfe_i$ &$4+9$ &$7$ & $6$& 4 \\
\hline
 \end{tabular}
\end{table}\label{2tb:Ex50}

\item $H=\langle 5, 7, 11, 13\rangle$ (non-symmetric) case $(d_h=5 \delta_0,\, \delta_0=5, \, \ell=0)$:
$H = \{0, 5, 7, 10, 11, 12,\allowbreak 13, 14, \dots\}$, $H^\fc = \{1, 2, 3, 4, 6, 8, 9\}$,
$R_H = \CC[Z_5, Z_7, Z_{11}, Z_{13}]/ {\sim} $ and then
\begin{gather*}
h_{{R_Z^H}}(Z, Z_\bullet, Z_\bullet') =
Z_5^5 + Z_5 Z_{13} Z_7'+ Z_{14} Z_{11}'+Z_5 Z_7 Z_{13}'+Z_{11} Z_{14}',
\qquad
h_H(Z) = 5 Z_5^5.
\end{gather*}
\begin{table}[htb]\centering\vspace{-5mm}
 \begin{tabular}{|c|ccccc|}
\hline
$i$&0&1&2&3&4\\
\hline
$\fe_i $ &0 &7 & 11& 13 &14 \\
$\hfe_i$ &25 &13+5 & 14& 7+5 &11 \\
\hline
 \end{tabular}
\end{table}\label{2tb:Ex5}

\item $H=\langle 6, 13, 14, 15, 16\rangle$ (symmetric) case
$(d_h=\fe_{r-1},\, \delta_0=0, \, \ell=r-1)$:
$H = \{0, 6, 12, 13, 14, \allowbreak 15, 16, 18, 19, 20, 21, 22, 24, \dots\}$,
$H^\fc = \{1, 2, 3, 4, 5, 7, 8, 9, 10, 11, 17, 23\}$.
$R_H = \CC[Z_6, Z_{13},\allowbreak Z_{14}, Z_{15}, Z_{16}]/ {\sim} $ and noting
$Z_{13}Z_{16}=Z_{14}Z_{15}$,
\begin{gather*}
h_{{R_Z^H}}(Z, Z_\bullet, Z_\bullet')=
Z_{13}Z_{16} +Z_{13}Z_{16}'+Z_{14}Z_{15}'+Z_{15}Z_{14}'
+Z_{16}Z_{13}'+Z_{13}'Z_{14}',\\
h_H(Z) = 6 Z_{13}Z_{16}.
\end{gather*}
\begin{table}[htb]\centering
 \begin{tabular}{|c|cccccc|}
\hline
$i$&0&1&2&3&4&5\\
\hline
$\fe_i$ &0 &13 & 14& 15 &16 & 29\\
$\hfe_i$ &29 &16 & 15& 14 &13 & 0\\
\hline
 \end{tabular}
\end{table}\label{2tb:Ex6}
\end{enumerate}
\end{Example}

By comparing the semigroup $H_X$ with $\{\fe_i\}$, we have the following corollaries:

\begin{Corollary}\label{2cr:h_Hhomo}\quad
\begin{enumerate}\itemsep=0pt
\item[$(1)$] $\th_{H}(Z_r, Z_\bullet, Z'_\bullet)$ is a homogeneous polynomial whose degree is $d_h=\fe_\ell+\delta_0 r$ in \eqref{2eq:hfe0},

\item[$(2)$]
$\th_{H}(Z_r, Z_\bullet, Z'_\bullet)
=\sum_{i=0}^{r-1}
Z_r^{\delta_i} \fZ_{\fe_{\ell,i}^*} \fZ_{\fe_i}'$,

\item[$(3)$]
$Z_r^{\delta_0} \fZ_{\fe_\ell}=Z_r^{\delta_i} \fZ_{\fe_{\ell,i}^*} \fZ_{\fe_i}$,

\item[$(4)$]
$h_{H}(Z)=r Z_r^{\delta_i} \fZ_{\fe_i^*}\fZ_{\fe_i}= r Z_r^{\delta_0} \fZ_{\fe_\ell}$ for $i\in \ZZ_r$.
\end{enumerate}
\end{Corollary}

\begin{proof}
They are obvious from Definition \ref{2df:NSG1}, Lemma \ref{2lm:NSG1}\,(3),
 Lemma~\ref{2lm:Z_standardbasis} and Lemma~\ref{2lm:h_H}.
\end{proof}

\begin{Corollary}\label{2cr:h_X-h_X}
$\fe_i-r -1 \ge 0$, $\hfe_i -r -1 \ge 0$ and $d_h - r -1 \ge 0$, for $i \in \ZZ_r^\times$.
\end{Corollary}

\begin{proof}
$\fe_i\ge \min_{j=2}^{m_X} r_j = r_2\ge r+1$ because of $r+1 \le r_2$.
\end{proof}

\begin{Remark}\label{2rk:RP-moduleRX2}
Corollary \ref{2cr:h_Hhomo} determines the $R_\PP$-module structure of $R_X$ in Proposition~\ref{2pr:RP-moduleRX}.
\end{Remark}

\subsubsection{Trace structure in Weierstrass curves (W-curves)}\label{2sssc:TraceStrWcurves}

We use the $R_\PP$-module structure of $R_X$ in the previous subsubsection and the properties in $R_H$ and $R_H^Z$ noting the surjection $\varphi_\infty\colon R_X\to R_H$ to define $p_X^{(j)}$ in the quotient field $\cQ(R_X\otimes_{R_\PP} R_X)$ for every affine ring $\CC[x,y]/\big(f_X^{(j)}\big)$:
\begin{Definition}\label{2df:p_Xj}
For $f_X^{(j)}\in \CC[x, y]$, we define
\[
p_X^{(j)}(x,y,y'):=\frac{
\big(\delta_{y,y'} f_X^{(j)}\big)(x,y,y')}{\big(f_{X,y}^{(j)}\big)(x,y)}.
\]
\end{Definition}

We regard $p_X^{(j)}(x,y,y')$ as an element in $\cQ\big(\CC[x,y]/\big(f_X^{(j)}\big)\otimes_{\CC[x]} \CC[x,y]/\big(f_X^{(j)}\big)\big)$ associated with $\cQ(R_X\otimes_{R_\PP} R_X)$.
We extend the group action $G_x$ to the action on $R_X\otimes_{R_\PP} R_X$, such that $\hzeta (x, y_\bullet, y_\bullet')= \big(x, \hzeta y_\bullet, \hzeta y_\bullet'\big)$.

\begin{Lemma}\label{2lm:p_Xj}
$p_X^{(j)}(x,y,y')=\begin{cases}
1 & \text{for } y = y',\\
0 & \text{for } y \neq y',
\end{cases}$
and for a group action
$\hzeta\in G_x$,
\[
p_X^{(j)}\big(x,\hzeta y,\hzeta y'\big)=p_X^{(j)}(x,y,y').
\]
\end{Lemma}

Let us consider $p_{R_X}:=\prod_{j=2}^{m_X} p_X^{(j)}$ as an element of $\cQ(R_X\otimes_{R_\PP} R_X)$.
The following is obvious:

\begin{Proposition}\label{2pr:p_varpi0}
For $(P,Q) \in X \times_\PP X$,
\[
p_{R_X}(P,Q)=
\begin{cases}
1 & \text{for } P = Q,\\
0 & \text{for } P \neq Q.
\end{cases}
\]
\end{Proposition}

However, some parts in its numerator and denominator are canceled because they belong to~$R_\PP$.
Thus we introduce an element $h(x, y_\bullet, y_\bullet') \in R_X\otimes_{R_\PP} R_X$ such that $h(x, y_\bullet, y_\bullet') /h(x, y_\bullet, y_\bullet)$
reproduces the product.

The ring homomorphism $\varphi^X_{R_X}$ in Proposition~\ref{2pr:R_H^Z} is extended to the surjective ring homomorphism from $R_X\otimes_{R_\PP}R_X$ to $R_H^Z\otimes_{\CC[Z_r]}R_H^Z$.

\begin{Lemma}\label{2lm:h_RXe}
For a point $(P = (x, y_\bullet), P'=(x, y_\bullet'))\in X \times_\PP X$, there is a polynomial $\th_{R_X}(x, y_{\bullet}, y_{\bullet}')\in R_X \otimes_{R_\PP} R_X$ such that
\begin{enumerate}\itemsep=0pt
\item[$(1)$] by regarding the element $a$ in $R_X$ as
$a\otimes 1$ in $R_X \otimes_{R_\PP} R_X$,
$\th_{R_X}(x, y_{\bullet}, y_{\bullet})$ and
$\th_{R_X}(x, y_{\bullet}, y_{\bullet}')$ are coprime
as elements in $R_X \otimes_{R_\PP} R_X$,

\item[$(2)$] for a group action $\hzeta \in G_x$,
$ \frac{\th_{R_X}(x, \hzeta y_{\bullet}, \hzeta y_{\bullet}')}
{\th_{R_X}(x, \hzeta y_{\bullet}, \hzeta y_{\bullet})}=
 \frac{\th_{R_X}(x, y_{\bullet}, y_{\bullet}')}
{\th_{R_X}(x, y_{\bullet}, y_{\bullet})}$,

\item[$(3)$] it satisfies
$\frac{\th_{R_X}(x, y_{\bullet}, y_{\bullet}')}
{\th_{R_X}(x, y_{\bullet}, y_{\bullet})}=
p_{R_X}(x, y_{\bullet}, y_{\bullet}')$, and

\item[$(4)$]
$\varphi_H^X\big(\th_{R_X}(x, y_{\bullet}, y_{\bullet}')\big)
=\th_H(Z_r, Z_\bullet, Z'_\bullet)
$, $\wt\big(\th_{R_X}(x, y_{\bullet}, y_{\bullet}')\big)=d_h$
with respect to $(x,y)$,
for~$\th_{H}$ in Lemma~{\rm \ref{2lm:h_H}}
by letting $\varphi_H^X(y_i) =Z_i$, $\varphi_H^X(y_i') =Z_i'$.
\end{enumerate}
\end{Lemma}

\begin{proof}By considering the numerator and denominator in $\prod_{j=2}^{m_X}p_X^{(j)}\bigr(\varpi_x(P), \varpi_{r_j}(P),\varpi_{r_j}(P')\bigr)$ modulo $\big(f_i^X\big)_{i=1, \dots, k_X}$, they are reduced to $\th_{R_X}$.
Let the numerator be denoted by $\th_{R_X}(x, y_{\bullet}, y_{\bullet}')\allowbreak \in R_X\otimes_{R_\PP}R_X$.
It is obvious that $\frac{\th_{R_X}(x, y_{\bullet}, y_{\bullet}')}{\th_{R_X}(x, y_{\bullet}, y_{\bullet})}$ must be invariant for the group action~$G_x$ on $R_X\otimes_{R_\PP}R_X$.
Due to the condition~(1), we have a unique ${\th_{R_X}(x, y_{\bullet}, y_{\bullet})}$.
Then we can find $d_h$ in Lemma~\ref{2lm:h_H} such that the image $\varphi^X_H\big(\th_{R_X}(x, y_{\bullet}, y_{\bullet}')\big)$ is equal to $\th_{H}$ in Lemma~\ref{2lm:h_H} because it is invariant for the $\GG_\fm$ action; the reduction in~(3) correspond to the reduction in $R_H^Z\otimes_{\CC[Z_r]}R_H^Z$ as in~(4).
Then it is clearly that (1), (2), and (3) are satisfied.
\end{proof}

\begin{Definition}\label{2df:h_X00}
Let $h_{X}(x, y_{\bullet}):=\th_{R_X}(x, y_{\bullet}, y_{\bullet})$.
\end{Definition}

Noting Corollary \ref{2cr:h_Hhomo}, we have the expression of $\th_{R_X}(x, y_{\bullet}, y_{\bullet}')$:

\begin{Proposition}\label{2pr:h_RXe_ys}$\th_{R_X}\in R_X\otimes_{R_\PP} R_X$ is expressed by
\begin{gather*}
\th_{R_X}(x, y_{\bullet}, y_{\bullet}')= \hUpsilon_{0}\cdot 1
+ \hUpsilon_{1}\fy'_{\fe_1}+\cdots+ \hUpsilon_{r-1}\fy'_{\fe_{r-1}}\\
\hphantom{\th_{R_X}(x, y_{\bullet}, y_{\bullet}')}{}
= 1\cdot\hUpsilon_{0}'
+ \fy_{\fe_1}\hUpsilon_{1}'+\cdots+ \fy_{\fe_{r-1}}\hUpsilon_{r-1}'\\
\hphantom{\th_{R_X}(x, y_{\bullet}, y_{\bullet}')}{}
= \rfy_{\fe_0}\fy'_{\fe_0}
+ \rfy_{\fe_1}\fy'_{\fe_1}\!+\cdots+ \rfy_{\fe_{r-1}}\fy'_{\fe_{r-1}} \!
 +\text{lower weight terms with respect to $-\wt$}\!
\end{gather*}
as an $R_\PP$-module.
Here $\fy'_{\fe_0}=1$, and each $\hUpsilon_{i}$ holds the following properties:
\begin{enumerate}\itemsep=0pt
\item[$(1)$]
$\hUpsilon_{i}=\sum_{j=0}^{r-1} \fb_{i,j} \fy_{\fe_j}$, with certain $\fb_{i,j}\in \CC[x]$, and

\item[$(2)$]
$\hUpsilon_{i}= \rfy_{\fe_i}+
\text{lower weight terms with respect to $-\wt$}$,
where
$\rfy_{\fe_i}
=\tdelta_i(x) \fy_{\fe_{\ell,i}^*}$
with a monic polynomial $\tdelta_i(x) \in \CC[x]$ whose weight is
$-\delta_i r$, $($especially, $\rfy_{\fe_0}
=\tdelta_0(x) \fy_{\fe_{\ell,0}^*}=\tdelta_0(x) \fy_{\fe_{\ell}})$
such that
\begin{gather*}
\rfy_{\fe_0}=\rfy_{\fe_i}\fy_{\fe_i}+
\text{lower weight terms with respect to $-\wt$}
\end{gather*} for
$i \in \ZZ_r$, $\wt(\rfy_{\fe_i})=-\hfe_i=-(d_h - \fe_i)$ in
Lemma~{\rm \ref{2lm:h_H}}, where $b_{i, j}$ is a certain element in~$R_\PP$ for $(i,j)$.
\end{enumerate}
\end{Proposition}

\begin{proof}From the construction, $\th_{R_X}(x, y_{\bullet}, y_{\bullet}')=\th_{R_X}(x, y_{\bullet}', y_{\bullet})$; $\th_{R_X}(x, y_{\bullet}, y_{\bullet}')$ is invariant for the exchanging between $y_{\bullet}$ and $y_{\bullet}'$.
From Proposition~\ref{2pr:RP-moduleRX}, $\hUpsilon_{i}$'s are uniquely determined.
It is obvious that $\fb_{i j}$ belongs to $\CC[x]=R_\PP$.
 However by letting $\ell_i$ satisfy
$\fe_{\ell_i}=\fe_{\ell,i}^*$, due to the Sato--Weierstrass weight of $\hUpsilon_{i}$, $-\wt(\fb_{i,j} \fy_{\fe_j}) < -\wt(\fb_{i,\ell_i} \fy_{\fe_{\ell_i}})$ and $\wt\big(\hUpsilon_{i}\big)=\wt(\fb_{i,\ell_i} \fy_{\fe_{\ell_i}})$.
Thus we let $\tdelta_i:=\fb_{i,\ell_i}$, and then we obtain $\rfy_{\fe_i}=\tdelta_i \fy_{\fe_{\ell,i}^*}$, in particular $\rfy_{\fe_0}=\tdelta_0 \fy_{\fe_\ell}$. $\rfy_{\fe_i}$ is monic from Corollary~\ref{2cr:h_Hhomo}.
Hence $\th_{R_X}(x, y_{\bullet}, y_{\bullet}')$ must have the form mentioned above.

Since $\rfy_{\fe_0}$ and $\rfy_{\fe_i}\fy_{\fe_i}$ are homogeneous with respect to $-\wt_\lambda$ and correspond to the relations in Corollary \ref{2cr:h_Hhomo}, $\rfy_{\fe_0}-\rfy_{\fe_i}\fy_{\fe_i}$ is given by the lower weight terms with respect to $-\wt$.
Further~(\ref{2eq:fyfy=fy}) shows $\tdelta_0 \fy_{\fe_\ell}$,
\begin{gather*}
\rfy_{\fe_i}\fy_{\fe_i} = \tdelta_i \fy_{\fe_{\ell_i}}\fy_{\fe_i}
=\sum_{j \in \ZZ_r} \tdelta_i \fa_{i, \ell_i, j}\fy_{\fe_j}\\
\hphantom{\rfy_{\fe_i}\fy_{\fe_i}}{} =\tdelta_i \fa_{i, \ell_i, \ell} \fy_{\fe_{\ell}} +\text{lower weight terms with respect to $-\wt$}.
\end{gather*}
Unless $\tdelta_i \fa_{i, \ell_i, \ell}/\tdelta_0$ belongs to $R_\PP$, it must be expressed as $\alpha/\beta$. where $\alpha$ and $\beta$ are elements of~$R_\PP$ and their Sato--Weierstrass weights are the same, and thus we redefine
\[
\rfy_{\fe_i}:=\beta \tdelta_i \fy_{\fe_{\ell,i}^*},
\qquad \rfy_{\fe_0}:=\alpha \tdelta_0 \fy_{\fe_\ell}.
\]
We repeat such operation for each $i$ and then, for every $i$, we obtain
\[
\rfy_{\fe_i}\fy_{\fe_i}
=\rfy_{\fe_0}+
\sum_{j \in \ZZ_r, j \neq \ell} \tdelta_i \fa_{i, \ell_i, j}\fy_{\fe_j}.
\]
However if $\th_{R_X}(x, y_{\bullet}, y_{\bullet}')$ and $\th_{X}(x, y_{\bullet}, y_{\bullet})=h_X(x, y_{\bullet})$ are not coprime, it means that there is a cofactor in $\th_{R_X}(x, y_{\bullet}, y_{\bullet}')$ and we can divide $\th_{R_X}(x, y_{\bullet}, y_{\bullet})$ by the factor.
\end{proof}

Recalling Proposition~\ref{2pr:dh_symmetric}, $\rfy_{\fe_i}$ and $h_{X}(x, y_{\bullet})$ are expressed as follows.
\begin{Corollary}\label{2cr:phi_h0}\quad
\begin{enumerate}\itemsep=0pt
\item[$1.$]
If $d_h$ is symmetric, $\rfy_{\fe_i}=\tdelta_i(x)\fy_{\fe_{r-1,i}^*}$
for $i\in \ZZ_r$,
$\tdelta_0=\tdelta_{r-1}=1$,
$\rfy_{\fe_0}=\fy_{\fe_{r-1}}$, and
$\rfy_{\fe_{r-1}}=\fy_{\fe_0}=1$.

\item[$2.$]
If $d_h$ is not symmetric, $\rfy_{\fe_i}=\tdelta_i(x)\fy_{\fe_{\ell,i}^*}$ for $i\in \ZZ_r^\times$, and $\tdelta_0\neq 1$.

\item[$3.$] $h_{X}(x, y_{\bullet})=r \rfy_{\fe_0}+$ lower weight terms with respect to $-\wt$.
\end{enumerate}
\end{Corollary}

We extend the arguments for the $m_X=2$ case to general $m_X$ cases.

Obviously, due to Proposition \ref{2pr:RP-moduleRX}, the dual $\omega_{R_X/R_\PP}$ of $R_X$ as an $R_\PP$-module has the standard basis as a trace.

\begin{Lemma}\label{2lm:Delta_th}
We define
\[
\Delta_{\th}:= \hUpsilon_{0}\otimes 1
+\sum_{k\in \ZZ_r^\times}
 \hUpsilon_{k}\otimes \fy_{\fe_k}
=\sum_{i,k\in \ZZ_r}\fy_{\fe_i}\fb_{i,k}\otimes \fy_{\fe_k}
=\sum_{i,k\in \ZZ_r}\fy_{\fe_i}\otimes\fb_{i,k} \fy_{\fe_k}.
\]
Then for every $a\in R_X$, $(a\otimes 1-1\otimes a)\Delta_{\th}=0$
and $\Delta_{\th}R_X^e$ is equal to $\Ann_{R_X^e}(\Ker \mu_\cQ)$.
\end{Lemma}

\begin{proof}
First we consider
\[
\fy_{\fe_j}\th_{R_X}(x, y_{\bullet}, y_{\bullet}') =\sum_{j,k}\fb_{j,k}\fy_{\fe_i}\fy_{\fe_j}\fy_{\fe_k}'=\sum_{j,k,\ell}\fb_{j,k}\fa_{i j\ell}\fy_{\fe_\ell}\fy_{\fe_k}'
\] and
\[ \fy_{\fe_j}'\th_{R_X}(x, y_{\bullet}, y_{\bullet}') =\sum_{j,k,\ell}\fb_{j,k}\fa_{i k\ell}\fy_{\fe_j}\fy_{\fe_\ell}'.\]
The latter can be expressed by $\sum_{\ell,j,k}\fb_{\ell,j}\fa_{i j k}\fy_{\fe_\ell}\fy_{\fe_k}'$.
Lemma \ref{2lm:h_RXe} with Proposition \ref{2pr:p_varpi0} shows that both cases agree
\[
\fy_{\fe_j}\th_{R_X}(x, y_{\bullet}, y_{\bullet}')
=\fy_{\fe_j}'\th_{R_X}(x, y_{\bullet}, y_{\bullet}')=
\begin{cases}
\fy_{\fe_j}h_X & \text{for } y_{\bullet} = y_{\bullet}',\\
0 & \text{otherwise}.
\end{cases}
\]
It means that $\sum_j \fb_{k,j}\fa_{i j\ell}=\sum_j \fb_{\ell,j}\fa_{i j k}$.
This relation shows
$\sum_{i,k}\fb_{j,k}\fy_{\fe_i}\fy_{\fe_j}\otimes \fy_{\fe_k}=\sum_{i,k}\fb_{j,k}\fy_{\fe_j}\otimes\fy_{\fe_i} \fy_{\fe_k}$
and we obtain the equality $(a\otimes 1-1\otimes a)\Delta_{\th}=0$.
\end{proof}

Further from Proposition \ref{2pr:A_traceH20}, $\Delta_{\th}$ in Lemma~\ref{2lm:Delta_th} provides the trace $\tau_\th\in \omega_{R_X/R_\PP}$ with the dual basis $\big\{\hUpsilon_{i}\big\}$ of $R_X$ with respect to the trace $\tau_\th$.

\begin{Lemma}\label{2lm:tau_th_hUpsiolon}
$\big\{\hUpsilon_{i}\big\}$ is the dual basis of $R_X$ with respect to the trace $\tau_\th$, such that
\[
\tau_\th(\hUpsilon_{i}) =
\begin{cases}
1 & \text{for } i=0,\\
0 & \text{otherwise}.
\end{cases}
\]
\end{Lemma}

\begin{proof}These $\big\{\hUpsilon_{i}\big\}$ correspond to the dual basis $\{\hby_i\}$ of $R_X$ with respect to the trace $\tau_\th$. By considering
\[
1= \tau_\th\big(\hUpsilon_{0}\big) \cdot 1
+\sum_{k\in \ZZ_r^\times}
 \tau_\th\big(\hUpsilon_{k}\big) \fy'_{\fe_k},
\]
we have the relation.
\end{proof}

We introduce an $R_\PP$-module,
$R_{X,\tau_\th}^*:=\big\langle\hUpsilon_{0}, \hUpsilon_{1}, \dots, \hUpsilon_{r-1} \big\rangle_{R_\PP}$,
and consider its structure as an $R_X$-module.
Due to Proposition~\ref{2pr:A_traceH20}, $\omega_{R_X/R_\PP}=R_X\circ \tau_\th$, and thus $\omega_{R_X/R_\PP}\cong R_{X,\tau_\th}^*$ as an $R_X$-module:

\begin{Lemma}\label{2lm:ideal_hUpsilon}
$R_{X,\tau_\th}^*$ is an ideal of $R_X$, especially
 $\fy_{\fe_i} R_{X,\tau_\th}^* \subset R_{X,\tau_\th}^*$, and
\[
R_{X,\tau_\th}^*=
\big\langle \hUpsilon_{0}, \hUpsilon_{1}, \dots, \hUpsilon_{r-1}\big\rangle_{R_X}=
\big\langle \hUpsilon_{1}, \dots, \hUpsilon_{r-1}\big\rangle_{R_X}.
\]
It means that as $R_\PP$-modules,
\[
\big\langle \hUpsilon_{0}, \hUpsilon_{1}, \dots, \hUpsilon_{r-1}\big\rangle_{R_\PP}=
\big\langle \hUpsilon_{1}, \dots, \hUpsilon_{r-1}\big\rangle_{R_X}.
\]
\end{Lemma}

\begin{proof}
Let us consider $\fy_{\fe_i} \hUpsilon_{j}$, which is equal to
\[
 \fy_{\fe_i} \hUpsilon_{j}
 =\sum_{k=0}^{r-1} \fb_{j k}\fy_{\fe_j}\fy_{\fe_k}
=\sum_{k,\ell=0}^{r-1} \fb_{j k}\fa_{i k\ell}\fy_{\fe_\ell}.
\]
However from the proof in Lemma \ref{2lm:Delta_th},
$\sum_k\fb_{j k}\fa_{i k\ell}=\sum_k\fb_{\ell k}\fa_{i k j}$, and thus
\[
 \fy_{\fe_i} \hUpsilon_{j}
 =\sum_{k=0}^{r-1} \fb_{k\ell}\fa_{i k j}\fy_{\fe_\ell}
=\sum_{k=0}^{r-1}\fa_{i k j}\hUpsilon_{k}.
\]
Hence, $R_{X,\tau_\th}^*$ is closed for the action of $\fy_{\fe_i}$.

On the other hand, we take an integer $i$, $i\neq \ell$, and then
consider $\fy_{\fe_i} \hUpsilon_{i}$ which is decomposed by
$\hUpsilon_{i}$'s but has the form,
\begin{equation}
 \fy_{\fe_i} \hUpsilon_{i}=\hUpsilon_0 +\sum_{j=1}^{r-1} \fd_{i j}\hUpsilon_{j},\label{2eq:fyhUpsilon}
\end{equation}
where $\fd_{i j}$ belongs to $R_\PP$ because the leading terms of the both $\fy_{\fe_i} \hUpsilon_{i}$ and $\hUpsilon_0$ must agree in the both sides.
It means that
$\hUpsilon_0 \in \big\langle \hUpsilon_{1}, \dots, \hUpsilon_{r-1}\big\rangle_{R_X}$.
\end{proof}

This proof shows the following lemma:

{\samepage
\begin{Lemma}\label{2lm:tau_th_rfe}
$\tau_\th\big(\fy_{\fe_i} \hUpsilon_{i}\big) =1$ for every $i \in \ZZ_r^\times$, and thus $\tau_\th\big(\fy_{\fe_i} \hUpsilon_{j}\big)=\delta_{i j}$.
\end{Lemma}

\begin{proof}(\ref{2eq:fyhUpsilon}) shows this.
\end{proof}}

Noting Remark~\ref{2rmk:hetam=2} for the $m_X=2$ case, we also introduce the more convenient basis $\{\hfy_{\fe_i}\}_{i \in \ZZ_R}$ with respect to $\tau_{\th}$ instead of $\{\hUpsilon_i\}_{i \in \ZZ_R}$:

\begin{Definition} \label{2df:hfy}
For $i \in \ZZ_r$, we define a truncated polynomial $\hfy_{\fe_i}$ of $\hUpsilon_{i}$ such that the weight $-\wt$ of $\hUpsilon_i -\hfy_{\fe_i}$ is less than $-\wt(\rfy_{\fe_i})$, i.e., $\hfy_{\fe_i}= \rfy_{\fe_i}+$ certain terms, and the number of the terms is minimal satisfying the relations as $R_\PP$-modules,
\begin{gather*}
\big\langle \hUpsilon_{1}, \dots, \hUpsilon_{r-1}\big\rangle_{R_X}
=\big\langle \hfy_{\fe_1}, \dots, \hfy_{\fe_{r-1}}\big\rangle_{R_X}
=\big\langle \hfy_{\fe_0}, \dots, \hfy_{\fe_{r-1}}\big\rangle_{R_\PP},
\qquad \tau_{\th}(\hfy_{\fe_i})=
\begin{cases}
1 & \text{for } i=0,\\
0 & \text{otherwise}.
\end{cases}\!
\end{gather*}
\end{Definition}

\begin{Lemma}\label{2lm:Emp_hfy_e}\quad
\begin{enumerate}\itemsep=0pt
\item[$1.$] When $d_h$ is symmetric, $\hfy_{\fe_i}=\fy_{\fe_{r-1-i}}$, especially for $m_X=2$ case, $\hfy_{\fe_i}=\rfy_{\fe_i}=y_s^i$.

\item[$2.$] When $X$ has a cyclic symmetry of order $r$, $\hfy_{\fe_i}=\rfy_{\fe_i}$.
\end{enumerate}
\end{Lemma}

\begin{proof}The $m_X=2$ case is described in Section~\ref{2sssec:m_X=2} and the statement is obvious.
For symmetric~$d_h$ case, under which the numerical semigroup $H_X$ is symmetric due to Proposition~\ref{2pr:dh_symmetric}, $\rfy_{\fe_i}=\fy_{\fe_{r-1-i}}$, $\rfy_{\fe_{r-1}}=\hUpsilon_{r-1}=1$, and thus it is obvious
\[
\big\langle \hUpsilon_{1}, \hUpsilon_{2}, \dots, \hUpsilon_{r-1}\big\rangle_{R_X}
=\langle \fy_{\fe_0}, \fy_{\fe_1}, \dots, \fy_{\fe_{r-2}}\rangle_{R_X}=R_X.
\]
Further, when $X$ has the cyclic symmetry of the order~$r$, $\hUpsilon_{i}= \rfy_{\fe_i}$ because of the invariance for the cyclic action.
\end{proof}

Since $\mu(\Delta_{\th})=h_X(x, y_\bullet)$, the trace $\tau_{R_X/R_\PP}$ of $R_X/R_\PP$ is given by
\[
\tau_{R_X/R_\PP}= h_X(x, y_\bullet)\circ \tau_\th.
\]

Then we obviously have the following lemma:
\begin{Lemma}\label{2lm:tau_RXRP_phi}
For a monomial $\phi$ in $R_X$ as an $R_\PP$-module and $\hfy_{\fe_i}$, we have the following:
\begin{enumerate}\itemsep=0pt

\item[$1.$] When $d_h$ is symmetric $(d_h= \fe_{r-1}, \delta_0=0)$ in~\eqref{2eq:hfe0},
\[
\tau_{R_X/R_\PP}
\left(\frac{\hfy_{\fe_i}}{h_X(x, y_\bullet)}\right)
=
\begin{cases}
1 & \text{for } i=0,\\
0 & \text{otherwise}.
\end{cases}
\]

\item[$2.$] When $d_h$ is not symmetric $(d_h= \fe_\ell + \delta_0 r, \delta_0\neq0)$ in~\eqref{2eq:hfe0},
\[
\tau_{R_X/R_\PP}
\left(\frac{\phi}{h_X(x, y_\bullet)}\right)
=
\begin{cases}
\phi/\hfy_{\fe_0}
& \text{for } \hfy_{\fe_0}|\phi, \\
0 & \text{otherwise}.
\end{cases}
\]
Here we note that the case $\hfy_{\fe_0}|\phi$ means $\big(\fy_{\fe_\ell}\tdelta_0\big)|\phi$ whereas the case $\hfy_{\fe_0}\nmid \phi$ consists of two cases $1)$~$\phi = (f(x) \fy_{\fe_\ell})$, $\tdelta_0\nmid f$ and $2)$ $\phi = (f(x) \fy_{\fe_i})\ i \in \ZZ_r\setminus \{\ell\}$.
\end{enumerate}
\end{Lemma}

Let us consider elements in $R_{X}^{\fc}$ as in (\ref{2eq:RXc_tau}).
We consider $u \in \cQ(R_X)$ given by
\begin{equation*}
u = \frac{1}{h_X(x, y_\bullet)} \sum_{i=0}^{r-1} a_i
 \hfy_{\fe_i}\in \cQ(R_X),
\end{equation*}
which satisfies the condition in (\ref{2eq:RXc_tau}).
Indeed, for any element $v = \sum_{j=0}^{r-1} b_j \fy_{\fe_j}\in R_X$, $b_j \in R_\PP$, $u$~satisfies $\tau_{R_X/R_\PP}(u v) \in R_\PP$, i.e.,
\[
\tau_{R_X/R_\PP}(u v)= \sum_{i,j}
 \tau_\th(a_i b_j \hfy_{\fe_i} \fy_{\fe_j})=
\sum_{i=0}^{r-1}(a_{i} b_{i})\in R_\PP.
\]
The condition $\tau_{R_X/R_\PP}(u v) \in R_\PP$ means that $a_i$ belongs to $R_\PP$ for every $i=0, 1, \dots, r-1$, and thus we obtain $R_{X}^\fc$.

We, now, state the first theorem in this paper.

\begin{Theorem}\label{2pr:R_Xfc}
The complementary module $R_X^\fc$ is given as a fractional ideal of $R_X$,
\[
R_X^\fc=\frac{\langle\hfy_{\fe_1}, \dots, \hfy_{\fe_{r-1}}\rangle_{R_X}}{h_X(x,y_\bullet)}.
\]
\end{Theorem}

\begin{proof}
It is obvious that due to the above arguments and the identity $\big\langle \hUpsilon_{0}, \hUpsilon_{1}, \dots, \hUpsilon_{r-1}\big\rangle_{R_X}=\big\langle\hfy_{\fe_1}, \dots, \hfy_{\fe_{r-1}}\big\rangle_{R_X} $ due to Lemma~\ref{2lm:ideal_hUpsilon} and Definition~\ref{2df:hfy}.
\end{proof}

Let us identify the generator $\fh_{X,P}$ of the principal ideal $R_{X,P}^\fc$ locally.

\begin{Proposition}\label{2pr:d_hcases}\quad
\begin{enumerate}\itemsep=0pt

\item[$1.$] Symmetric $d_h$ case $(\delta_0 = 0$ in \eqref{2eq:hfe0} or $\tdelta_0=1$ and $d_h = \fe_{r-1})$:
The complementary module $($as a fractional ideal$)$ is given by
\begin{equation}
R_X^\fc=\frac{1}{h_X(x,y_\bullet)} R_X= \fh_X R_X,\label{2eq:R_X^fc1}
\end{equation}
and we define $\fh_X:=\frac{1}{h_X(x,y_\bullet)}$.

\item[$2.$] Non-symmetric $d_h$ case $(\delta_0 \neq 0$ in \eqref{2eq:hfe0}$)$ or $\tdelta_0\neq 1$ and $d_h = \fe_{\ell}+\delta_0 r$:
For complex numbers $(a_i(\neq 0))_i$, e.g., $a_i=1$, we define
\begin{equation*}
\fh_X:=\frac{\sum_{i=1}^{r-1}a_i \hfy_{\fe_i}}{h_X(x,y_\bullet)}.
\end{equation*}
\end{enumerate}
For the both cases, we obtain the local expression of the complementary module at $P\in X$,
\begin{equation}
R_{X,P}^\fc= \fh_{X,P} R_{X,P}.
\label{2eq:R_X^fc2}
\end{equation}
\end{Proposition}

\begin{proof}
Symmetric $d_h$ case: The ideal $\cI$ contains $1$ because $\hfy_{\fe_{r-1}}=1$.
It includes $m_X=2$ case in Section~\ref{2sssec:m_X=2}.
In other words, $\fh_{X,P}$ in Definition \ref{2df:2.32} is given by $\varphi_P(\fh_X)=\fh_{X,P}$ by using the homomorphism $\varphi_P\colon R_X \to R_{X,P}$ in Remark~\ref{2rm:g_Xc_X}.

Non-symmetric $d_h$ case:
By noting
\[
\deg_{P,0}\left(\sum_{i=1}^{r-1}a_i \hfy_{\fe_i}\right)=
\min\big(\deg_{P,0}\big( \hfy_{\fe_i}\big)\big),
\]
we also identify $\fh_{X,P}$ in Definition \ref{2df:2.32} with $\varphi_P(\fh_X)$.
\end{proof}

Using the complementary module for both cases, we have the Dedekind different,
\[
\diff(R_X/R_\PP)=-
\sum_{P\in X\setminus\{\infty\}}
\deg_P(\fh_{X}) P.
\]

By Dedekind's different theorem (Proposition~\ref{2pr:DedekindDiff1}), we have the following.
\begin{Proposition}\label{2pr:DedekindDiff30}
$e_{B_i}-1=-\deg_{B_i}(\fh_X)$, and the support of $\Div(\fh_X)$ is equal to $\fB_X$.
\end{Proposition}

Since some of $f_{X,y}^{(j)}(P)=0$ at $P=B_i \in \fB_X\setminus\{\infty\}$, $h_X(x, y_\bullet)\in R_X$ vanishes only at the ramification point $B_i\in X$ from Proposition \ref{2pr:DedekindDiff30} and the construction of $h_X$, we have the following corollary:

\begin{Corollary}\label{2cr:h_X}
\[
\Div(h_X(x,y_\bullet)) =
\sum_{B_i\in \fB_X\setminus\{\infty\}}d_{B_
i} B_i
-d_h\infty,
\]
where $d_{B_i}:=\deg_{B_i,0}(h_x) \ge (e_{B_i}-1)$, and $d_h = \sum_{B_i\in \fB_X\setminus\{\infty\}}d_{B_i}=-\wt (h_X)$.
\end{Corollary}

\begin{Definition}\label{2df:fKX_fkX}
The effective divisor, $\sum_{B_i\in \fB_X\setminus\{\infty\}}(d_{B_i}-e_{B_i}+1)B_i$, is denoted by
$\fK_X$, i.e., $\fK_X>0$ and let
$\fk_X:=\sum_{B_i\in \fB_X\setminus\{\infty\}}(d_{B_i}-e_{B_i}+1)=\deg(\fK_X)\ge0$.
\end{Definition}

\begin{Lemma}\label{2lm:dx_hX}
The divisor of $\frac{{\rm d} x}{h_X}$ is expressed by $(2g-2+\fk_X)\infty-\fK_X$, and $2g-2+\fk_X=d_h - r- 1$ or $\fk_X=d_h - 2g- r+1$.
\end{Lemma}

\begin{proof}
Since the meromorphic one-form in general has degree $2g-2$, we obtain $\deg({\rm d} x/h_X)=2g-2$. From Definition~\ref{2df:fKX_fkX}, its divisor is expressed by $\Div({\rm d} x/h_X)=(2g-2+\fk_X)\infty-\fK_X$.
Further since at the $\infty$, its degree $\deg_\infty({\rm d} x/h_X)=d_h - r-1$, we have $2g-2+\fk_X=d_h - r- 1$.
\end{proof}

From Corollary \ref{2cr:(dx)}, we note that these $\fK_X$ and $\fk_X$ play crucial roles in the investigation of the differentials on $X$ (e.g., Lemma \ref{2lm:nuIo1}).

\begin{Proposition}\label{2pr:fK_X}
$\fk_X$ is equal to zero if $d_h$ is symmetric whereas $\fk_X$ is not zero otherwise.
\end{Proposition}

\begin{proof}
For the symmetric case, (\ref{2eq:R_X^fc1}) and Proposition \ref{2pr:DedekindDiff30} show that $d_{B_i}=e_{B_i}-1$, whereas for the non-symmetric case, (\ref{2eq:R_X^fc2}) and Proposition~\ref{2pr:DedekindDiff30} yield non-vanishing $\fk_X$.
\end{proof}

We recall $\hfe_i$ in (\ref{2eq:hfe}), $\fe_i$ in Definition \ref{2df:NSG1} for the standard basis in Lemma~\ref{2lm:Z_standardbasis} and Proposition~\ref{2pr:RP-moduleRX}, and $\hfy_{\fe_i}$ in Definition~\ref{2df:hfy}.

By applying Propositions~\ref{2pr:R_Xfc} and~\ref{2pr:DedekindDiff30} to differentials on $X$, we consider $\frac{x^k \hfy_{\fe_i}(x,y_\bullet) {\rm d} x}{h_X(x,y_\bullet)}$, which holds the following proposition:

\begin{Proposition}\label{2pr:phi_hi}
\begin{gather*}
\Div\left(\frac{x^k \hfy_{\fe_i}(x,y_\bullet) {\rm d} x}
 {h_X(x,y_\bullet)}\right)
= \Div\left(\frac{x^k \tdelta_i(x)\fy_{\fe_{\ell,i}}(x,y_\bullet) {\rm d} x}
 {h_X(x,y_\bullet)}\right)=k\Div_0(x)+\Div_0(\hfy_{\fe_i})
\\
\hphantom{\Div\left(\frac{x^k \hfy_{\fe_i}(x,y_\bullet) {\rm d} x}
 {h_X(x,y_\bullet)}\right)=}{}
-\sum_{j=1}^{\ell_{\fB}} (d_{B_j}-e_{B_j} -1) B_j
+(d_h - \hfe_{i}-(k+1)r-1) \infty\\
\hphantom{\Div\left(\frac{x^k \hfy_{\fe_i}(x,y_\bullet) {\rm d} x}
 {h_X(x,y_\bullet)}\right)}{}
=k\Div_0(x)+\Div_0(\hfy_{\fe_i})-\fK_X+(\fe_{i}-(k+1)r-1) \infty,
\end{gather*}
where $\Div_0(\hfy_{\fe_i})-\fK_X\ge0$.
We have
\[
\left\{
\wt
\left(\frac{x^k \hfy_{\fe_i}(x,y_\bullet) {\rm d} x}{h_X(x,y_\bullet)}\right)+1 \,
\Bigr| \, i \in \ZZ_r^\times,\, k\in \NN_0,\, \fe_{i}-(k+1)r>0\right\}
=H_X^\fc,
\]
and
\[
\left\{
\wt
\left(\frac{x^k \hfy_{\fe_i}(x,y_\bullet) {\rm d} x}{h_X(x,y_\bullet)}\right)+1 \,
\Bigr| \, i \in \ZZ_r,\, k\in \NN_0\right\}
=\BH_X^\fc.
\]
\end{Proposition}

\begin{proof}The former statement is asserted from the previous lemma, whereas the latter two relations on $H_X^\fc$ and $\BH_X^\fc$ proved by Lemma~\ref{2lm:NSG1} noting $\wt(\hfy_{\fe_i})=-\hfe_i=-(d_h-\fe_i)$ due to~(\ref{2eq:hfe}) and $\fe_0 =0$.
\end{proof}

\section[W-normalized Abelian differentials H\^{}0(X,A\_X(*infty))]{W-normalized Abelian differentials $\boldsymbol{\bH^0(X, \cA_X(*\infty))}$}\label{2ssc:W-norm_nuI}

Following K.~Weierstrass \cite{WeiWIV},
H.F.~Baker \cite{Baker97}, V.M.~Buchstaber, D.V.~Leykin and V.Z.~Enolskii~\cite{BEL20}, J.C.~Eilbeck, V.Z.~Enolskii and D.V.~Leykin~\cite{EEL00}, we construct the Abelian differentials of the first kind and the second kind $\bH^0(X, \cA_X(*\infty))$ on $X$ for the hyperelliptic curves and plane Weierstrass curves (W-curves).
We extend them to more general W-curves based on Proposition~\ref{2pr:phi_hi}.

We consider the Abelian differentials of the first kind on a general W-curve.
Due to the Riemann--Roch theorem, there is the $i$-th holomorphic one-form
whose behavior at $\infty$ is given by
\begin{equation*}
 \bigr(t^{N^\fc(g-i)-1} (1+ d_{>0}(t))\bigr) {\rm d} t,
\end{equation*}
where $N^\fc(i) \in H_X^\fc$, $i=1, 2, \dots, g$, $N^\fc(i) <N^\fc(i+1)$, and $t$ is the arithmetic local parameter at~$\infty$.
We call this normalization the {\emph{W-normalization}}.
Similarly we find the differentials or the basis of $\bH^0(X, \cA_X(*\infty))$ associated with $\BH_X^\fc$.

\subsection{W-normalized Abelian differentials}

The W-normalized holomorphic one-forms are directly obtained from Proposition \ref{2pr:phi_hi}:
\begin{Lemma}\label{2lm:hphi_hi}
For $x^k \hfy_{\fe_i}$ in Proposition~{\rm \ref{2pr:phi_hi}}, we have the relation,
\[
\left\langle\frac{x^k\hfy_{\fe_i}}{h_X(x,y_\bullet)}{\rm d} x \, \Bigr|
\, i \in \ZZ_r^\times, k\in \NN_0, \fe_{i}-(k+1)r>0 \right\rangle_\CC
= \bH^0(X, \cA_X).
\]
\end{Lemma}

By re-ordering $x^k\hfy_{\fe_j}$ with respect to the weight at $\infty$, we define the ordered set $\big\{\hphi_i\big\}$:

\begin{Definition}\label{2df:nuI_hX}\quad
\begin{enumerate}\itemsep=0pt
\item
Let us define the ordered subset
$\hS_X$ of $R_X$ by
\[
\hS_X=\big\{\hphi_i\, | \, i \in \NN_0\big\}
\]
such that $\hphi_i$ is ordered by the Sato--Weierstrass weight,
i.e., $-\wt \,\hphi_i < -\wt\, \hphi_j$ for $i < j$,
and
$\hS_X$ is equal to $\big\{x ^k \hfy_{\fe_i} \, |\,i \in \ZZ_r,\, k \in \NN_0\big\}$ as a set.

\item Let $\hR_X$ be an $R_X$-module generated by $\hS_X$, i.e., $\hR_X:=\langle \hS_X\rangle_{R_X} \subset R_X$.

\item Recalling $\fK_X$ and $\fk_X$ in Definition~\ref{2df:fKX_fkX}, we let $\hN(n):= -\wt \big(\hphi_n\big)-\fk_X$, $\hH_X:=\big\{{-}\wt \big(\hphi_n\big)\, |$ $n \in \NN_0\big\}$, and we define the dual conductor $\hc_X$ as the minimal integer satisfying $\hc_X + \NN_0 \subset \hH_X-\fk_X$.

\item We define $\hS_X^{(g)}:=\big\{\hphi_0, \hphi_1, \dots, \hphi_{g-1}\big\}$ and the \emph{W-normalized holomorphic one form}, or the \emph{W-normalized Abelian differentials of the first kind} $\nuI{i}$ as the canonical basis of $X$,
\begin{equation}
\left\langle\nuI{i} := \frac{\hphi_{i-1} {\rm d} x }{h_X}\,\Bigr|\,
\hphi_{i-1} \in \hS_X^{(g)}\right\rangle_\CC = \bH^0(X, \cA_X).\label{2eq:nuI_hX}
\end{equation}
\end{enumerate}
\end{Definition}

We note that at $\infty$, $\nuI{i}$ behaves like
$
\nuI{i}=\big(t^{N^\fc(g-i-1)-1} (1+ d_{>0}(t))\big) {\rm d} t$ for the arithmetic local parameter $t$ at $\infty$, and further
$\frac{\hphi_{i-1} {\rm d} x }{h_X}=
\big(t^{N^\fc(g-i-1)-1}(1+d_{>0})\big){\rm d}t$ where $N^\fc(i)$ indicates the element in $\BH_X^\fc$ such that $N^\fc(-i)=-i$ for $i \in \NN$; they are W-normalized Abelian differentials.

Finally we state our second theorem, which is obvious.

\begin{Theorem}\label{2lm:nuI_hX}\quad
\begin{gather*}
1. \ \bH^0(X, \cA_X(*\infty)) = \bigoplus_{i=0} \CC\frac{\hphi_i {\rm d} x}{h_X}=R_X^\fc {\rm d} x = \frac{\hR_X {\rm d} x}{h_X}.\\
2. \
{\overline H}_X^{\mathfrak c}=\left\{ \wt\left(\frac{\hphi_i {\rm d} x}{h_X}\right)+1 \, |\,i \in \NN_0\right\}
\\
\hphantom{2. \ {\overline H}_X^{\mathfrak c}}{}
 =\big\{ \wt\big(\hphi_i\big)+d_h -r \, |\,i \in \NN_0\big\} =\big\{ \wt\big(\hphi_i\big)+2g-1 -\fk_X \, |\,i \in \NN_0\big\}.
\end{gather*}
\end{Theorem}

The Riemann--Roch theorem shows that
\begin{equation*}
\dim \bH^0(X, \cO_X(-n\infty))-\dim \bH^0(X, \cA_X(n\infty))-n=1-g.
\end{equation*}

\begin{Lemma}\label{2lm:nuIo1}
If $d_h=\fe_{r-1}$ or $\delta_0 = 0$ in~\eqref{2eq:hfe0} $(d_h$ is symmetric in Proposition~{\rm \ref{2pr:dh_symmetric})}, $H_X$ is symmetric, otherwise $(\delta_0 \neq 0$ or $d_h$ is not symmetric$)$ $\fk_X\neq 0$ and $H_X$ is not symmetric.
\end{Lemma}

\begin{proof}If $d_h$ is symmetric, $\fk_X=0$ from Proposition~\ref{2pr:fK_X}.
Thus if $d_h$ is symmetric, $({\rm d} x/h_X)=(2g-2)\infty$, and due to the Riemann--Roch theorem, $H_X$ is symmetric.
It corresponds to $\delta_0 = 0$ and $d_h = \fe_{r-1}$ in Proposition~\ref{2pr:dh_symmetric} and Lemma~\ref{2lm:h_H}.
On the other hand, the case $\fk_X\neq 0$ means that~$H_X$ is not symmetric and $\delta_0 \neq 0$ in Lemma~\ref{2lm:h_H}. Thus we prove it.
\end{proof}

\begin{Proposition}\label{2pr:dh_HXsymmetric}\quad
\begin{enumerate}\itemsep=0pt
\item[$1.$] Assume $d_h$ is symmetric or $d_h = \fe_{r-1}$, $\delta_0 = 0$, in~\eqref{2eq:hfe0}.
Then we have the following:
\begin{enumerate}\itemsep=0pt
\item[$(a)$] $\fk_X=0$ in Definition~{\rm \ref{2df:fKX_fkX}} and $H_X$ is symmetric,

\item[$(b)$] $\hS_X = S_X$ in Definitions~{\rm \ref{2df:t_SR_N}} and {\rm \ref{2df:nuI_hX}}, $\fK_X=0$ in Definition~{\rm \ref{2df:fKX_fkX}}, $\hc_X=c_X$, and $\hR_X = R_X$ in Definition~{\rm \ref{2df:nuI_hX}}.
\end{enumerate}
\item[$2.$] In general, $\hR_X \neq R_X$ and we have the equality if and only if $H_X$ is symmetric.
\end{enumerate}
\end{Proposition}

\begin{proof}We note Proposition \ref{2pr:N(n)}\,(5). They are obvious.
\end{proof}

By the Abel--Jacobi theorem, $\fK_X$ in Definition \ref{2df:fKX_fkX} can be divided into two pieces, which are related to the spin structure in $X$ or the half-canonical form \cite{KMP16}.

\begin{Definition}\label{2df:fK_s}
Let $\fK_\fs$ and $\fK_X^\fc$ be the effective divisors which satisfy
\[
\fK_X-\fk_X\infty \sim2\fK_\fs-2 \fk_\fs\infty,\qquad
\fK_X + \fK_X^\fc - \big(\fk_X+\fk_X^\fc\big)\infty\sim 0
\]
as the linear equivalence, where $\fk_\fs$ and $\fk_X^\fc$ are the degree of $\fK_\fs$ and $\fK_X^\fc$ respectively.
\end{Definition}

Since the W-normalized holomorphic one form is given by the basis (\ref{2eq:nuI_hX}), Definition \ref{2df:fK_s} shows the canonical divisor:

\begin{Proposition}\label{2pr:cKX}
The canonical divisor is given by
\begin{gather*}
K_X \sim(2g-2+\fk_X)\infty-\fK_X\sim (2g-2+2\fk_\fs)\infty - 2\fK_\fs
\sim \big(2g-2-\fk_X^\fc\big)\infty +\fK_X^\fc.
\end{gather*}
\end{Proposition}

This expression can be applied to the shifted Riemann constant for the non-symmetric W-curves~\cite{KMP16}.
Theorems \ref{2pr:R_Xfc} and \ref{2lm:nuI_hX} enable us to define the fundamental 2-form of the second kind algebraically and to construct the sigma functions of every W-curve as we show in a follow-up paper~\cite{KMP2022b}.
Using the results~\cite{KMP16}, we connect them with the sigma functions, which is defined as a modified Nakayashiki's sigma function~\cite{KMP2022b}.
We find the explicit relations between $R_X$ and the meromorphic functions
on the Jacobi variety $J_X$ associated with the sigma functions like Weierstrass' elliptic function theory.

As mentioned in introduction, Segal and Wilson showed that $\BH_X^\fc$ provides the embedding of the algebraic systems associated with $X$ into the UGM \cite[p.~46]{SegalWilson}.
In contrast, we find the $R_X$-module $R_X^\fc {\rm d}x$ as an algebraic system with the same Sato--Weierstrass weight as $\BH_X^\fc-1$ explicitly.
Though Nakayashiki defined his sigma functions based on the embedding of the UGM \cite{Nak16} in terms of the {\lq}wave function{\rq} with a half (spin) density, it is expected that our results might directly show the construction of the sigma functions by the embedding of the complementary module $R_X^\fc dx$ into the UGM as a natural generalization of his approach in~\cite{Nak10b} for plane W-curves.

For the application in \cite{KMP2022b}, as the end of the above discussion, we will summarize the properties of these parameters.

\begin{Lemma}\label{2lm:nuIo1_2}\quad
\begin{enumerate}\itemsep=0pt
\item[$1.$] $\big\{{-}\wt\big(\hphi_i\big)\big\}=\{\hfe_i + k r\, |\,i \in \ZZ_r,\, k \in \NN_0\}=\{d_h - \fe_i + k r\, |\,i \in \ZZ_r,\, k \in \NN_0\}$.

\item[$2.$] $\Div\big(\hphi_i\big)\ge(\fK_X-(2g-2+\fk_X)\infty)$ for every $\hphi_i
\in \hS_R^{(g)}$, $i=0, 1, 2,\dots,g-1$.

\item[$3.$] $\Div\big(\hphi_i\big)\ge(\fK_X-(g+\fk_X+i)\infty)$ for every $\hphi_i \in \hS_R$, $i\ge g$.
\end{enumerate}
\end{Lemma}

\begin{proof}They are obvious.
\end{proof}

\begin{Lemma}\label{2lm:nuIo2}\quad
\begin{enumerate}\itemsep=0pt
\item[$(1)$] $-\wt\, \hphi_{0}= \hfe_{r-1}$ $(=0$ if $H_X$ is symmetric$)$
$= \hc_X+\fk_X-c_X =d_h - r-c_X+1$,

\item[$(2)$] $-\wt\, \hphi_{g-1}= \hc_X+\fk_X-2= d_h - r- 1=
(2g-2)+\fk_X$, i.e., $\hN(g-1) = 2g-2$,

\item[$(3)$] $ \hc_X = 2g=d_h -\fk_X-r+1$,
$-\wt\, \hphi_{g}=2g+\fk_X=\hc_X+\fk_X=d_h - r+1$,

\item[$(4)$] $-\wt\, \hphi_{g-1}+\wt\, \hphi_{0}=\fe_{r-1}-r-1= c_X-2$, and
$c_X= \fe_{r-1} - r+1$.

\end{enumerate}
\end{Lemma}

\begin{proof}
Lemma \ref{2lm:nuIo1}\,(2) means that $-\wt\big(\hphi_{0}\big)=\hfe_{r-1}$ due to
the order of the weight.
Then~(1) is proved by the relations $-\wt\, \hphi_{0}=\hN(0)+\fk_X$ and $\hN(0) = \hc_X-c_X$.
Let us consider~(2). Since from the Riemann--Roch theorem, $\wt \, \nuI{g} = 0$ whereas on $\nuI{g}=\hphi_{g-1}{\rm d} x/h_X$, $\wt({\rm d} x/h_X)=(2g-2)+\fk_X=d_h - r-1$, we have $-\wt\, \hphi_{g-1}=(2g-2)+\fk_X$ or~(2).
The Riemann--Roch theorem also shows that $-\wt\, \hphi_{g}=-\wt\, \hphi_{g-1}+2 = \hc_X +\fk_X$, or~(3).
We compare them and obtain~(4).
\end{proof}

\begin{Remark}\label{rmk:Final}
As we show in a follow-up paper \cite{KMP2022b}, we mention how the results in this paper provide the connection between~$R_X$ and the sigma function shortly in this remark.
We extend $p_{R_X}(P,Q)$ for $(P,Q) \in X \times_\PP X$ in Proposition~\ref{2pr:p_varpi0} and Lemma~\ref{2lm:h_RXe} to
\[
p_{\varpi}(P,Q)
:=\frac{\th_X(x_P, y_{P\bullet}, y_{Q\bullet})}
{\th_X(x_P, y_{P\bullet})}
\]
for $(P=(x_P, y_{P\bullet}),Q=(x_Q, y_{Q\bullet})) \in X \times X$ as in \cite[Definition~12]{KMP2022b}.
This extension of the domain from $X \times_\PP X$ to $X \times X$ is not unique in general.
However, it has an excellent property
\[
p_{\varpi}(P,Q)
=
\begin{cases}
1 & \text{for } P = Q,\\
0 & \text{for } P \neq Q \text{ and }
 \varpi_x(P)=\varpi_x(Q),
\end{cases}
\]
for $X \times X$ except ramification points. Thus we introduce the one-forms \cite[Proposition~15]{KMP2022b},
\begin{equation*}
 \Sigma (P, Q )
 :=
\frac{{\rm d} x_P}{(x_P - x_Q)}p_{\varpi}(P,Q)
=
\frac{{\rm d} x_P}{(x_P - x_Q)}
\frac{\th_X(x_P, y_{\bullet P}, y_{\bullet Q})}
{h_X(x_P, y_{\bullet P})}.
\end{equation*}
By investigating the one-form $\Sigma\big(P, Q\big)$ and its derivative ${\rm d}_Q\Sigma (P, Q )$ in~$Q$, we can define the W-normalized differentials of the second kind and the third kind, and the fundamental differential of the second kind $\Omega(P, Q)$ such that \cite[Theorem~3]{KMP2022b},
\begin{enumerate}\itemsep=0pt
\item[(1)] $\Omega(P, Q)=\Omega(Q, P)$,
\item[(2)] for any $\zeta \in G_{\varpi_r(P)}$,
$\Omega(\zeta P, \zeta Q)=\Omega(P, Q)$ if $\varpi_r(P)=\varpi_r(Q)$, and
\item[(3)]
 $\Omega(P, Q)$ is holomorphic except $Q$ as a function of $P$ and behaves like
\[
\Omega(P, Q)=\frac{{\rm d} t_P {\rm d} t_Q}{(t_P-t_Q)^2}(1 + d_{>0}(t_P, t_Q)).
\]
\end{enumerate}
It turns out that for any extension of $p_\varpi$, the differential $\Omega$ is unique under the cohomological meaning so that the choice of how we select the extension does not affect the final results essentially \cite[Proposition~16]{KMP2022b}.
In the considerations, the properties of the complementary module, as the results in this paper, play a crucial role.
Further, as the differential $\Omega$ is connected with the differential of the third kind, we have its connection to the sigma function, $\sigma$ \cite[Theorem~4]{KMP2022b}:
for $(P, Q, P_i, P'_i) \in X^2 \times \big(S^g(X)\setminus S^g_1(X)\big) \times
\big(S^g(X)\setminus S^g_1(X)\big)$,
\begin{gather}
	u := \tw_\fs(P_1, \dots, P_g), \qquad
	v := \tw_\fs(P'_1, \dots, P'_g),\nonumber\\
\exp\left(
\sum_{i, j = 1}^g
 \Pi_{P_i, P'_j}^{P, Q} \right)
=
\frac{\sigma(\tw(P) - u) \sigma(\tw(Q) - v)}
 {\sigma(\tw(Q) - u) \sigma(\tw(P) - v)}.
\label{eq:5.5}
\end{gather}
where $\Pi^{P_1, P_2}_{Q_1, Q_2} := \int^{P_1}_{P_2} \int^{Q_1}_{Q_2} \Omega(P, Q)$, and $\tw$ and~$\tw_\fs$ are the ordinary and the shifted Abelian integrals.
The sigma function $\sigma$ is the modified version of Nakayashiki's one~\cite{Nak16} based on the results in~\cite{KMP16}.
In the Weierstrass elliptic function theory, as $\wp(u-v){\rm d}u{\rm d}v$ has the double order pole at $u=v$ and the elliptic sigma function, $\sigma(u-v)$, is connected with the integral $\wp(u-v)$ with respect to ${\rm d}u$ and ${\rm d}v$,
(\ref{eq:5.5})~means that based on our results in this paper, we can generalize the picture to every W-curve as mentioned in \cite[Theorem~4]{KMP2022b}.
Our results in this paper undoubtedly contribute to the significant progress in the Weierstrass sigma function theory for general algebraic curves.
\end{Remark}

\section{Examples of Weierstrass curves (W-curves)}\label{sec:Exmples}

\subsection{Special other curves: pentagonal, non-cyclic trigonal, 6-symmetric curves}\label{2ssc:SpCurve}

I. Non-cyclic trigonal curve (3,7,8):
$y^3 + a_1 k_2(x) y^2
+ a_2\tk_2(x) k_2(x) y +k_2(x)^2 k_3(x)=0$,
where $k_2(x) = (x-b_1)(x-b_2)$,
$k_3(x) = (x-b_3)(x-b_4)(x-b_5)$
$\tk_2(x) = (x-b_6)(x-b_7)$, for pairwise distinct $b_i\in \CC$ and $a_j$ generic constants.
Here (\ref{2eq:WCF3}) and (\ref{2eq:WCF4}) correspond to
\[
y^2 = -a_1k_2 y - k_2 a_2\tk_2 - k_2 w, \qquad
y w = k_2 k_3.
\]
Multiplying the first equation by $y$ and using the second equation gives the curve's equation.
Besides them, we have
\[
w^2 = -\big( a_2\tk_2 w + a_1k_2 k_3 + k_3y\big),
\]
since multiplying the first equation by $w^2$ gives
\[
w^3 + a_2\tk_2 w^2 + a_1k_2 k_3 w + k_2 k_3^2=0.
\]
This curve is trigonal with $H_X=\langle 3,7,8\rangle$ but not necessarily cyclic,
\begin{gather*}
\th_X(x, y, w, y', w')=k_2 k_3^2 + k_3(y w'+ y' w)\\
\hphantom{\th_X(x, y, w, y', w')=}{} + \frac{1}{3} k_3\big(2 a_1 y y'+a_1 k_2 (w+w')+(a_1^2k_2+a_2 \tk_2)(y+y')+2a_1a_2k_2 \tk_2\big)\\
\hphantom{\th_X(x, y, w, y', w')=}{} + \frac{1}{3}\tk_2\big(a_1a_2( y w'+y' w)+2a_2w w'+a_2^2 \tk_2(w+w')\big),\\
h_X(x,y, w) = 3k_2 k_3^2.
\end{gather*}

\begin{table}[t]\centering
	\caption{Examples of $\phi$ of special curves.}\vspace{1mm}
\renewcommand{\arraystretch}{1.3}\small
 \begin{tabular}{|r|ccccccccccccc|}
\hline
 &0&1&2&3&4&5&6&7&8&9&10&11&12\\
\hline
I&$1$&--&--&$x$&--&--&$x^2$&$y$&$w$&$x^3$&$xy$&$xw$&
 $x^4$\\
II&$1$&--&--&--&--&$x$&--&$y$&--&--&$x^2$&$w$&$xy$\\
III&$1$&--&--&--&--&--&$x$&--&--&--&--&--&$x^2$\\
\hline
\hline
&13&14&15&16& 17&18&19&20&21&22&23&24&\\
\hline
I&$x^2y$&$y^2$&$yw$&$w^2$&$xy^2$&$xyw$&$xw^2$&$x^2y^2$&
 $x^2yw$&$x^2w^2$&$x^3y$&$x^3w^2$&\\
II&--&$y^2$&$x^3$&$xw$&$x^2y$&$wy$&$xy^2$&$x^4$&$y^3$&$x^3y$&$xyw$&$x^2y^2$&\\
III&$y_{13}$&$y_{14}$&$y_{15}$&$y_{16}$
 &--&$x^3$&$xy_{13}$&$xy_{14}$&$xy_{15}$&$xy_{16}$&--&$x^4$&\\
\hline
 \end{tabular}
\end{table}

The differentials of the first kind are given as follows:
\[
\nuI{1}=\frac{y {\rm d} x}{3k_2k_3^2},\qquad
\nuI{2}=\frac{w {\rm d} x}{3k_2k_3^2},\qquad
\nuI{3}=\frac{xy {\rm d} x}{3k_2k_3^2},\qquad
\nuI{4}=\frac{xw {\rm d} x}{3k_2k_3^2}.
\]

II. Cyclic pentagonal curve (5,7,11):
$y^5 = k_2(x)^2 k_3(x)$, where
$k_2(x) = (x-b_1)(x-b_2)$,
$k_3(x) = (x-b_3)(x-b_4)(x-b_5)$ for pairwise distinct $b_i\in \CC$:
(\ref{2eq:WCF4}) corresponds to
\[
\begin{pmatrix}
-y & 0 & 1\\
k_2 & -y & 0\\
0 & 0 & - y
\end{pmatrix}
\begin{pmatrix}
w \\ y^2\\ wy
\end{pmatrix}
=\begin{pmatrix}
0 \\ 0\\ -k_2^2 k_3
\end{pmatrix}.
\]
The affine ring is
$R_X^\circ =\CC[x, y, w]/\big(y^3-k_2w, w^2-k_3y, y^2 w - k_2 k_3\big)$.
Here (\ref{2eq:WCF5}) is reduced to
\begin{gather*}
w = \frac{k_2 k_3}{y^2},\qquad
yw = \frac{k_2 k_3}{y},\qquad
y^2 = \frac{k_2 w}{y}.
\end{gather*}
This is a pentagonal cyclic curve $(X,\infty )$ with $H_X=\langle
5,7,11\rangle$.
\begin{gather*}
\th_X(x, y, w, y', w')=y^2 w + y w y' + y^2 w' + w y^{\prime\,2}+y w' y',\\
h_X(x,y, w) = 5 y^2 w =5k_2k_3.
\end{gather*}

These $\fe$'s and $d_h$ in Lemma \ref{2lm:h_H} are given as
$\fe_0 = 0$, $\fe_1 = 7$, $\fe_2 = 11$,
$\fe_3 = 14$, $\fe_4 = 18$, and $d_h = 25$.

The differentials of the first kind are given as follows:
\[
\nuI{i}=\frac{\hphi_{i-1} {\rm d} x}{5k_2k_3}.
\]

III. 6-symmetric (6,13,14,15,16) curve:
We construct a non-singular curve $X$ by giving an affine patch, an ideal in the ring $\CC[x, y_{13}, y_{14},$ $y_{15}, y_{16}]$.
For any complex numbers $\{b_i\}_{i=1, \dots, 7}$ such that each is distinct from the others, we let
\begin{gather*}
 k_3(x) := (x - b_1)(x - b_2)(x - b_3)
= x^3 + \lambda^{(3)}_{1} x^2 + \lambda^{(3)}_{2} x + \lambda^{(3)}_{3}, \\
 k_2(x) := (x - b_4) (x - b_5)
= x^2 + \lambda^{(2)}_{1} x + \lambda^{(2)}_{2} . \\
 \hk_{2}(x) := (x-b_6)(x-b_7)
 = x^2 + \hlambda^{(2)}_1 x + \hlambda^{(2)}_0, \qquad
 \hk_{5}(x) := \hk_{2}(x) k_{3}(x), \\
 k_{13}(x) := k_3(x)k_2(x)^2 \hk_2(x)^3, \qquad
 k_{14}(x) := k_7(x)^2 = k_3(x)^2k_2(x)^4 , \\
 k_{15}(x) := \hk_5(x)^3, \qquad k_{16}(x) := k_8(x)^2 = k_3(x)^4k_2(x)^2 .
\end{gather*}
The Weierstrass canonical form is given by $y_6^6=\hk_2^3k_2^2 k_3$, which is normalized as follows.

\begin{table}[t]\centering
 \caption{Examples of $\hphi$ in $\hS_X$ of special curves with respect to $\wt\big(\hphi_i\big)$.}\vspace{1mm}
\renewcommand{\arraystretch}{1.3}\small
 \begin{tabular}{|r|ccccccccccccc|}
\hline
 &0&1&2&3&4&5&6&7&8&9&10&11&12\\
\hline
I&--&--&--&--&--&--&--&$y$&$w$&--&$xy$&$xw$&
 --\\
II&--&--&--&--&--&--&--&$y$&--&--&--&$w$&$xy$\\
III&$1$&--&--&--&--&--&$x$&--&--&--&--&--&$x^2$\\
\hline
\hline
&13&14&15&16& 17&18&19&20&21&22&23&24&\\
\hline
I&$x^2y$&$y^2$&$yw$&$w^2$&$xy^2$&$xyw$&$xw^2$&$x^2y^2$&
 $x^2yw$&$x^2w^2$&$x^3y$&$x^3w^2$&\\
II&--&$y^2$&--&$xw$&$x^2y$&$wy$&$xy^2$&--&$y^3$&$x^3y$&$xyw$&$x^2y^2$&\\
III&$y$&$z$&$v$&$w$&--&$x^3$&$xy$&$xz$&$xv$&$xw$&--&$x^4$&\\
\hline
 \end{tabular}
\end{table}

Let the prime ideal $\cP$ in $\CC[x, y_{13}, y_{14}, y_{15}, y_{16}]$ be defined by
\begin{equation*}
 \cP := (f_{12, 1}, f_{12, 2}, f_{12, 3}, f_{12, 4},
f_{12, 5}, f_{12, 6}, f_{12, 7}, f_{12, 8}, f_{12, 9}),
\end{equation*}
where
\begin{gather*}
f_{12,1} := y_{13}^{2} - \hk_{2}(x) y_{14}, \qquad
f_{12,2} := y_{13} y_{14} - k_{2}(x) y_{15}, \qquad
f_{12,3} :=\hk_2(x) y_{14}^{2} - y_{13} y_{15} k_2(x),\! \\
f_{12,4} := y_{14}^{2} - k_{2}(x) y_{16} , \qquad
f_{12,5} :=y_{13}y_{16} - y_{14} y_{15}, \qquad
f_{12,6} := y_{15}^{2} - \hk_2(x)k_3(x), \\
f_{12,7} := y_{14}y_{16} - k_2(x)k_3(x), \qquad
f_{12,8} := y_{15} y_{16} - k_3(x)y_{13}, \qquad
f_{12,9} := y_{16}^{2} - k_3(x)y_{14},
\end{gather*}
which are the $2 \times 2$ minors of
\begin{equation*}
\left|\begin{matrix}
 k_2(x) & y_{14} & y_{16} \\
 y_{14} & y_{16} & k_3(x)\\
\end{matrix} \right|, \qquad
\left|\begin{matrix}
 \hk_2(x) & y_{13} & y_{14} \hk_2(x) & y_{15} \\
 y_{13} & y_{14} & y_{15} k_2(x) & y_{16} \\
\end{matrix} \right|;
\end{equation*}
again, the minor $y_{13} y_{15} - \hk_2(x) y_{16}$ is not in the list of $f_{i,j}$ and $f_{12,8}$ is not a minor, but they are compatible~-- the minor follows by combining $f_{12,8}$ with $f_{12,1}$ and $f_{12,9}$.

We define the $\GG_\fm$ action on $x$ and $y_a$ by $g^{-6} x$ and $g^{-a} y_a$, $a = 13, 14, 15, 16$, at $\infty \in X_{12}$.

Corresponding to Proposition \ref{2pr:WCF2}, the affine ring is given by
\[
R_{X}^\circ =\CC[x, y_{13}, y_{14}, y_{15}, y_{16}]/ \cP,
\]
and $\th_X$ and $h_X$ are
\begin{gather*}
\begin{split}
&\th_X(x, y_{\bullet,1}, y_{\bullet 2}) =
 y_{13,1} y_{16,1}
+ y_{13,2} y_{16,1}
+ y_{13,1} y_{16,2}
+ y_{14,2} y_{15,1}
+ y_{14,1} y_{15,2}
+ y_{14,2} y_{15,2},\\
&h_X(x, y_\bullet) = 6 y_{13} y_{16}.
\end{split}
\end{gather*}

The differentials of the first kind are given as follows:
\[
\nuI{i}=\frac{\phi_{i-1} {\rm d} x}{6 y_{13}y_{16}}=\frac{\phi_{i-1} {\rm d} x}{6 y_{14}y_{15}}.
\]

\subsection*{Acknowledgements}

The second author is grateful to Professor Yohei Komori for valuable discussions in his seminar and to Professor Takao Kato for letting him know about the paper \cite{CoppensKato}.
He has been supported by the Grant-in-Aid for Scientific Research (C) of Japan Society for the Promotion of Science Grant, no.~21K03289.

At the end of May 2022, since this paper was almost completed, the authors decided that it would be submitted at the beginning of July and would go through a month of checks by typos and others.
In the meantime, Emma Previato, the third author of this paper, passed away on June 29, 2022.
The first two authors sincerely wish that the great mathematician and their friend and collaborator Emma Previato rest in peace.

Thus the revised version was written only by the first two authors.
They thank the anonymous reviewers for their helpful and valuable comments.

\pdfbookmark[1]{References}{ref}
\LastPageEnding

\end{document}